\documentclass[a4paper,11pt]{article}
\usepackage{amsmath,amsthm,amssymb}
\usepackage{graphicx}
\usepackage[square,sort&compress,numbers]{natbib}
\usepackage{type1cm}
\usepackage{url}
\usepackage{mathtools}\mathtoolsset{showonlyrefs=false}
\usepackage{booktabs}

\usepackage{subcaption}
\newcommand{\Authornote}{\renewcommand{\thefootnote}{\fnsymbol{footnote}}}
\newcommand{\authornote}{\Authornote\footnote}
\theoremstyle{plain}

\theoremstyle{definition}

\theoremstyle{remark}

\newcommand{\reffig}[1]{Figure~\ref{#1}}
\newcommand{\reftab}[1]{Table~\ref{#1}}

\newcounter{alnum}

\newenvironment{Problem}{\begin{array}.{*{20}{l}}\}}{\end{array}}

\newcommand{\supp}{\mathop{\mathrm{supp}}\nolimits}

\newcommand{\MIN}{\mathop{\mathrm{Minimize}}}

\newcommand{\ST}{\mathop{\mathrm{subject~to}}}

\newcommand{\argmin}{\operatornamewithlimits{\mathrm{arg\,min}}}

\newcommand{\overRe}{\ensuremath{\Re\cup\{+\infty\}}}

\renewcommand{\Re}{\ensuremath{\mathbb{R}}}

\newcommand{\bi}[1]{\ensuremath{\boldsymbol{#1}}}

\newcommand{\rr}[1]{\ensuremath{\mathrm{#1}}}

\newcommand{\KC}{\ensuremath{\mathcal{K}}}
\newcommand{\LC}{\ensuremath{\mathcal{L}}}

\newcommand{\UC}{\ensuremath{\mathcal{U}}}

\begin{document}

\begin{center}
  {\Large\bfseries\sffamily%
  Alternating Direction Method of Multipliers }
  \par\medskip
  {\Large\bfseries\sffamily%
  for Truss Topology Optimization with Limited Number of Nodes: }
  \par\medskip
  {\Large\bfseries\sffamily%
  A Cardinality-Constrained Second-Order Cone Programming Approach 
  }%
  \par%
  \bigskip%
  {
  Yoshihiro Kanno~\authornote[2]{%
  Corresponding author. 
  Mathematics and Informatics Center, 
  The University of Tokyo, 
  Hongo 7-3-1, Tokyo 113-8656, Japan.
  E-mail: \texttt{kanno@mi.u-tokyo.ac.jp}. 
  },
  Shinnosuke Fujita~\authornote[3]{%
  Laboratory for Future Interdisciplinary Research of Science and Technology, 
  Institute of Innovative Research, 
  Tokyo Institute of Technology, 
  Nagatsuta 4259, Yokohama 226-8503, Japan.
  E-mail: \texttt{fujita.s.ag@m.titech.ac.jp}. 
  }
  }
\end{center}

\begin{abstract}
  This paper addresses the compliance minimization of a truss, where the 
  number of available nodes is limited. 
  It is shown that this optimization problem can be recast as a 
  second-order cone programming with a cardinality constraint. 
  We propose a simple heuristic based on the alternative direction 
  method of multipliers. 
  The efficiency of the proposed method is compared with a global 
  optimization approach based on mixed-integer second-order cone 
  programming. 
  Numerical experiments demonstrate that the proposed method often finds 
  a solution having a good objective value with small computational cost. 
\end{abstract}

\begin{quote}
  \textbf{Keywords}
  \par
  Topology optimization; 
  truss optimization; 
  manufacturability; 
  alternating direction method of multipliers; 
  cardinality-constrained second-order cone programming; 
  mixed-integer second-order cone programming. 
\end{quote}

\section{Introduction}

It is common to use the ground structure method \cite{Top83,Kir89,BBtZ94} 
for truss topology optimization, where the cross-sectional areas of the 
truss members are treated as design variables to be optimized. 
Particularly, the compliance minimization with the continuous design 
variables is convex \cite{BBtZ94,Ohs11}, and can be solved efficiently. 
An optimal solution of this problem often consists of 
too many members (including ones that are too thin) connected 
by many nodes\footnote{%
See, for example, \reffig{fig:x8_x9_nominal} in Section~\ref{sec:ex}.} 
and, hence, is regarded as too complex a design from the viewpoint of 
practical manufacturability. 
Also, the fabrication cost of a truss usually increases as the number of 
nodes increases. 
To obtain a practically acceptable truss design, 
\citet{AGV15} proposed to minimize the weighted 
sum of the structural weight and the fabrication cost related to the 
number of members. 
In this method, the number of members is approximated by using a 
regularized Heaviside step function. 
\citet{TLM16} used the same approach to take into 
account the number of nodes. 
In this paper, we consider the compliance minimization of a truss 
subjected to the explicit upper bound constraint on the number of nodes. 

The number of nodes in structural optimization has also been discussed 
in the layout optimization of trusses. 
In the classical layout optimization, we minimize the total 
weight of the members when the allowable stress is specified. 
When the potential locations of nodes of a truss are not limited,
the optimal solution becomes a so-called truss-like continuum 
with infinitely many nodes \cite{HP69,Mic04}. 
\citet{Pra78,Pra77} showed that, 
by adding the weight of the nodes to the objective function, 
we can obtain an optimal solution with a finite number of nodes. 
To avoid complex truss design, \citet{Par75} proposed to introduce 
modification of member lengths such that, at each node, a constant is 
added to the length of each member connected to the node. 
As a post-processing step for this method, 
\citet{HG15} proposed to make use of the geometry optimization. 
Similarly, \citet{MBT11} defined a so-called performance index, by using 
the member lengths and the axial forces, to assess the cost of a 
structure; see also \cite{Maz12}. 
The number of nodes in a truss is not specified explicitly in the 
methods in the literature 
\cite{AGV15,TLM16,Par75,HG15,Maz12,MBT11} cited above. 

In this paper, based on the ground structure method we deal with the 
compliance minimization problem of a truss subjected to the upper bound 
constraint on the number of nodes (i.e., the cardinality constraint on 
the set of nodes). 
This design optimization problem essentially consists of two decisions: 
We first select a set of nodes, satisfying the cardinality constraint, 
among the candidate nodes in a ground structure, and next find the 
optimal cross-sectional areas of the members connected to the selected 
nodes. 
The first decision gives combinatorial attribute to the design 
optimization problem. 
In this paper, we show that this optimization problem can be recast as 
{\em mixed-integer second-order cone programming\/} (MISOCP); 
see Section~\ref{sec:node.integer}. 
Since an SOCP problem can be solved efficiently with a primal-dual 
interior-point method \cite{AL12,BtN01}, 
we can compute a global optimal solution of an 
MISOCP problem with, e.g., a branch-and-bound method. 
Several software packages are available for this purpose 
\cite{ART03,Gurobi}. 
However, due to its large computational cost, the MISOCP approach can be 
applied only to small- to medium-size 
truss optimization problems. 
The reader may refer to \cite{BS09,MT15} for applications of MISOCP to 
variable selection in statistics, and 
\cite{Kan13,Kan16,Kan16b,Koc17} for applications in structural 
optimization. 

The number of nonzero components of a real vector is called the 
$\ell_{0}$-norm of the vector.\footnote{
Although this number is not a norm, it is common to 
call it the $\ell_{0}$-norm \cite{BDE09,BKS16,CWB08,Cha07,GTT15,LPLV15,ZSLS14}. 
} 
An upper bound constraint on the $\ell_{0}$-norm of a vector, i.e., the 
$\ell_{0}$-norm constraint, is also called the cardinality constraint 
(i.e., the upper bound constraint on the cardinality of the support 
of the vector). 
The cardinality constraint, as well as the $\ell_{0}$-norm minimization, 
frequently appears in diverse fields including variable selection in 
statistics, image processing, compressed sensing, and 
portfolio selection 
\cite{Nat95,Cha07,BDE09,CWB08,GTT15,LPLV15,ZSLS14,BKS16,BS09,CZZS13}. 
An application of the $\ell_{0}$-minimization to structural design 
generating link mechanisms can be found in \cite{OKT14}. 
In this paper, we show that the truss topology optimization with the 
limited number of nodes can be formulated as cardinality-constrained 
SOCP; see Section~\ref{sec:node.setting}. 

The {\em alternating direction method of multipliers\/} (ADMM) is an 
algorithm for convex optimization~\cite{BPCPE10}. 
For various nonconvex optimization problems, it is known that ADMM can 
often serve as a simple but powerful heuristic 
\cite{TMBB15,KT14,CW13,MRF16,Cha12,DTB17}. 
This motivates us to develop a simple heuristic based on ADMM, to find 
approximate solutions to the truss topology optimization with limited 
number of nodes. 
The proposed method might be expected to find a local optimal solution 
having the reasonable objective value with small computational cost. 
In control theory, ADMM has been used for various 
sparsity-promoting optimal control method, including design of 
sparse feedback gains \cite{LFJ13}, 
sparse output feedback \cite{ABKM15}, and 
a sparse gain matrix for the extended Kalman filter \cite{MFV12}.


The paper is organized as follows: 
Section~\ref{sec:preliminary} provides an overview of necessary 
backgrounds of ADMM. 
Section~\ref{sec:node} formulates the truss topology optimization 
problem with the limited number of nodes as cardinality-constrained SOCP, 
and recasts it as MISOCP. 
Section~\ref{sec:admm} presents a heuristic based on ADMM for the 
formulation as cardinality-constrained SOCP. 
Section~\ref{sec:overlap} is devoted to discussion on treatment of 
overlapping members in a ground structure. 
Section~\ref{sec:ex} reports the results of numerical experiments. 
Some conclusions are drawn in Section~\ref{sec:conclude}.


In our notation, 
${}^{\top}$ denotes the transpose of a vector or a matrix. 
We use $\bi{1} = (1,1,\dots,1)^{\top}$ to denote the all-ones vector. 
For vectors $\bi{x} = (x_{i}) \in \Re^{n}$ and 
$\bi{y} = (y_{i}) \in \Re^{n}$, we write $\bi{x} \ge \bi{y}$ 
if $x_{i} \ge y_{i}$ $(i=1,\dots,n)$. 
We use $\| \bi{x} \|$ to denote the Euclidean norm 
(or the $\ell_{2}$-norm) of $\bi{x}$, i.e., 
$\| \bi{x} \| = \sqrt{\bi{x}^{\top} \bi{x}}$. 
We denote by $\| \bi{x} \|_{0}$ the number of nonzero components of 
$\bi{x}$, which is the so-called $\ell_{0}$-norm of $\bi{x}$. 
For a finite set $T$, let $|T|$ denote the cardinality of $T$, 
i.e., the number of elements in $T$. 
If we define $\supp(\bi{x}) \subseteq \{ 1,\dots,n \}$ by 
$\supp(\bi{x}) 
  = \{ i \in \{ 1,\dots,n \} \mid 
  x_{i} \not= 0 \}$, 
then $\| \bi{x} \|_{0} = | \supp(\bi{x}) |$. 
Therefore, a constraint on the $\ell_{0}$-norm is also called the 
cardinality constraint. 
For a set $S \subseteq \Re^{n}$, we denote by 
$\delta_{S} : \Re^{n} \to \overRe$ the indicator function of $S$, which 
is defined by 
\begin{align*}
  \delta_{S}(\bi{x}) = 
  \begin{cases}
    0 \quad 
    & \text{if } \bi{x} \in S, \\
    +\infty \quad 
    & \text{if } \bi{x} \not\in S. \\
  \end{cases}
\end{align*}
For a closed set $S \subseteq \Re^{n}$, the projection of a point 
$\bi{z} \in \Re^{n}$ onto $S$, denoted $\Pi_{S}(\bi{z}) \in \Re^{n}$, is 
defined by 
\begin{align*}
  \Pi_{S}(\bi{z})  \in S, 
  \quad
  \| \bi{z} - \Pi_{S}(\bi{z}) \|
  = \min \{ \| \bi{z} - \bi{x} \| \mid \bi{x} \in S \} .
\end{align*}
If $S$ is closed and convex, then $\Pi_{S}(\bi{z})$ exists uniquely for 
any point $\bi{z} \in \Re^{n}$. 
The $n$-dimensional second-order cone, denoted $\LC^{n}$, is 
defined by 
\begin{align*}
  \LC^{n} = \{
  (s_{0},\bi{s}_{1}) \in \Re \times \Re^{n-1}
  \mid
  \| \bi{s}_{1} \| \le s_{0} 
  \} . 
\end{align*}
The $n$-dimensional rotated second-order cone, denoted $\KC^{n}$, is 
defined by 
\begin{align*}
  \KC^{n} = \{
  (\bi{x},y,z) \in \Re^{n-2} \times \Re \times \Re
  \mid
  \bi{x}^{\top} \bi{x} \le y z , \
  y \ge 0 , \
  z \ge 0 
  \} . 
\end{align*}
We have that $(\bi{x},y,z) \in \KC^{n}$ if and only if 
$(y+z, y-z, 2\bi{x}) \in \LC^{n}$. 
We use $\UC(a,b)$ to denote the continuous uniform 
distribution on the interval $(a,b) \subset \Re$.

\section{Fundamentals of alternating direction method of multipliers}
\label{sec:preliminary}

In this section, we briefly outline the 
{\em alternating direction method of multipliers\/} (ADMM) 
for solving convex optimization; see \cite{BPCPE10} for more accounts. 

Let $f : \Re^{n} \to \Re \cup \{ +\infty \}$ and
$g : \Re^{m} \to \Re \cup \{ +\infty \}$ be closed proper convex 
functions. 
Consider the following convex optimization problem in variables 
$\bi{x} \in \Re^{n}$ and $\bi{z} \in \Re^{m}$: 
\begin{subequations}\label{P.convex.1}%
  \begin{alignat}{3}
    & \MIN_{\bi{x},\bi{z}}  &{\quad}& 
    f(\bi{x}) + g(\bi{z}) \\
    & \ST && 
    A \bi{x} + B \bi{z} = \bi{c} . 
  \end{alignat}
\end{subequations}
Here, $A \in \Re^{l \times n}$ and $B \in \Re^{l \times m}$ are constant 
matrices, and $\bi{c} \in \Re^{l}$ is a constant vector. 

The augmented Lagrangian of problem \eqref{P.convex.1} is defined as 
\begin{align}
  L_{\rho}(\bi{x},\bi{z},\bi{y}) 
  = f(\bi{x}) + g(\bi{z}) 
  + \bi{y}^{\top} (A \bi{x} + B \bi{z} - \bi{c}) 
  + \frac{\rho}{2} \| A \bi{x} + B \bi{z} - \bi{c} \|^{2} , 
  \label{eq.ADMM.Lagrangian.1}
\end{align}
where $\rho > 0$ is the penalty parameter, and 
$\bi{y} \in \Re^{l}$ is the Lagrange multiplier (also called the dual 
variable). 
At each iteration of ADMM, we update $\bi{x}^{k}$, $\bi{z}^{k}$, and 
$\bi{y}^{k}$ as 
\begin{align}
  \bi{x}^{k+1} 
  &:= \argmin_{\bi{x}} L_{\rho}(\bi{x},\bi{z}^{k},\bi{y}^{k}) , 
  \label{eq.ADMM.update.1.1} \\
  \bi{z}^{k+1} 
  &:= \argmin_{\bi{z}} L_{\rho}(\bi{x}^{k+1},\bi{z},\bi{y}^{k}) , 
  \label{eq.ADMM.update.1.2} \\
  \bi{y}^{k+1} 
  &:= \bi{y}^{k} 
  + \rho (A \bi{x}^{k+1} + B \bi{z}^{k+1} - \bi{c}) . 
  \label{eq.ADMM.update.1.3}
\end{align}

The so-called {\em scaled form\/} of ADMM is defined below. 
Letting $\bi{v} = \bi{y} / \rho$, we see that 
\eqref{eq.ADMM.Lagrangian.1} is reduced to 
\begin{align}
  \tilde{L}_{\rho}(\bi{x},\bi{z},\bi{v}) 
  = f(\bi{x}) + g(\bi{z}) 
  + \frac{\rho}{2} \| A \bi{x} + B \bi{z} - \bi{c} + \bi{v} \|^{2} 
  - \frac{\rho}{2} \| \bi{v} \|^{2} . 
  \label{eq.ADMM.Lagrangian.2}
\end{align}
By using $\tilde{L}_{\rho}$ in \eqref{eq.ADMM.Lagrangian.2}, the 
iteration of ADMM given by \eqref{eq.ADMM.update.1.1}, 
\eqref{eq.ADMM.update.1.2}, and \eqref{eq.ADMM.update.1.3} is written as 
\begin{align}
  \bi{x}^{k+1} 
  &:= \argmin_{\bi{x}} \tilde{L}_{\rho}(\bi{x},\bi{z}^{k},\bi{v}^{k}) , 
  \label{eq.ADMM.update.2.1} \\
  \bi{z}^{k+1} 
  &:= \argmin_{\bi{z}} \tilde{L}_{\rho}(\bi{x}^{k+1},\bi{z},\bi{v}^{k}) , 
  \label{eq.ADMM.update.2.2} \\
  \bi{v}^{k+1} 
  &:= \bi{v}^{k} + A \bi{x}^{k+1} + B \bi{z}^{k+1} - \bi{c} .
  \label{eq.ADMM.update.2.3}
\end{align}
The form given in \eqref{eq.ADMM.update.2.1}, \eqref{eq.ADMM.update.2.2}, 
and \eqref{eq.ADMM.update.2.3} is called the scaled form of ADMM, and 
$\bi{v}$ is called the scaled dual variable. 

Primarily, ADMM is an algorithm for solving convex optimization. 
It is known that ADMM can often serve as an efficient heuristic for 
diverse nonconvex optimization problems; see, e.g., 
\cite{TMBB15,KT14,CW13,MRF16,Cha12}, and \cite[Section~9]{BPCPE10}. 
For nonconvex problems, ADMM does not necessarily converge. 
Also, when it converges, the obtained solution is not necessarily 
optimal. 
Furthermore, the obtained solution can depend on the penalty parameter 
and the initial point. 
Nevertheless, ADMM can be a simple algorithm, and can be efficient in 
the sense that it often converges to a solution with a good objective 
value.

\section{Design optimization with limited number of nodes}
\label{sec:node}

In Section~\ref{sec:node.setting}, we define truss topology optimization 
under the upper bound constraint on the number of nodes. 
In Section~\ref{sec:node.integer}, we show that this problem can be 
recast as MISOCP. 

\subsection{Problem setting}
\label{sec:node.setting}

\begin{figure}[tp]
  \centering
  \includegraphics[scale=1.20]{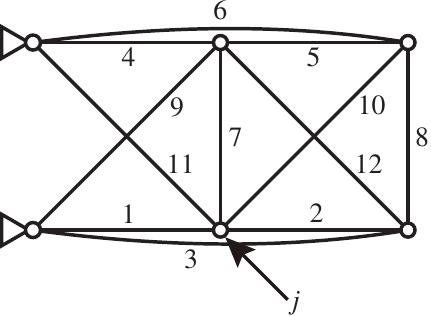}
  \caption{An example of ground structure with overlapping members.}
  \label{fig:gs2x1}
\end{figure}

Following the ground structure approach, consider an initial truss 
consisting of many candidate members that are connected by nodes with 
the given locations. 
The cross-sectional areas of the members are treated as the design 
variables to be optimized. 
It is worth noting that the ground structure may involve some 
overlapping members, as an example shown in \reffig{fig:gs2x1}. 
The necessity, as well as the treatment, of overlapping members in a 
ground structure is thoroughly discussed in Section~\ref{sec:overlap}. 
We use $m$, $l$, and $d$ to denote the number of members, the number of 
nodes, and the number of degrees of freedom of the nodal displacements, 
respectively. 

Let $x_{i}$ $(i=1,\dots,m)$ denote the member cross-sectional areas. 
We use $K(\bi{x}) \in \Re^{d \times d}$ to denote the stiffness matrix, 
which can be written as 
\begin{align*}
  K(\bi{x}) 
  = \sum_{i=1}^{m} \frac{E}{c_{i}} x_{i} \bi{b}_{i} \bi{b}_{i}^{\top} .
\end{align*}
Here, $c_{i}$ is the undeformed member length, 
$E$ is the Young modulus, and 
$\bi{b}_{i} \in \Re^{d}$ is a constant vector 
reflecting the member connectivity and the direction cosine of member $i$. 
For a given external load vector $\bi{p} \in \Re^{d}$, the compliance of 
the truss, denoted $\pi(\bi{x})$, is defined by 
\begin{align}
  \pi(\bi{x}) 
  = \sup\{ 2\bi{p}^{\top} \bi{u} - \bi{u}^{\top} K(\bi{x}) \bi{u} 
  \mid \bi{u} \in \Re^{d} \} . 
  \label{eq:def.compliance}
\end{align}
Let $V$ $(>0)$ denote the specified upper bound for the structural 
volume. 
The conventional compliance minimization problem is formulated as 
follows: 
\begin{subequations}\label{P.compliance.1}%
  \begin{alignat}{3}
    & \MIN_{\bi{x}}
    &{\quad}& 
    \pi(\bi{x}) \\
    & \ST && 
    \bi{x} \ge  \bi{0} , \\
    & && 
    \bi{c}^{\top} \bi{x} \le V . 
  \end{alignat}
\end{subequations}
This problem is convex, and can be recast as follows 
\citep[Section~3.4.3]{BtN01}: 
\begin{subequations}\label{P.compliance.SOCP.1}%
  \begin{alignat}{3}
    & \MIN_{\bi{x},\bi{q},\bi{w}}
    &{\quad}& 
    \sum_{i=1}^{m} w_{i} \\
    & \ST && 
    w_{i} x_{i} \ge \frac{c_{i}}{E} q_{i}^{2} , 
    \quad i=1,\dots,m , 
    \label{P.compliance.SOCP.1.2} \\
    & && 
    \bi{x} \ge \bi{0} , 
    \label{P.compliance.SOCP.1.4} \\
    & && 
    \sum_{i=1}^{m} q_{i} \bi{b}_{i} = \bi{p} , \\
    & && 
    \bi{c}^{\top} \bi{x} \le V . 
  \end{alignat}
\end{subequations}
Constraints \eqref{P.compliance.SOCP.1.2} and 
\eqref{P.compliance.SOCP.1.4} can be rewritten equivalently as the 
rotated second-order cone constraints 
\begin{align*}
  (\sqrt{c_{i}/E}q_{i} , w_{i}, x_{i}) \in \KC^{3} , 
  \quad  i=1,\dots,m . 
\end{align*}
These constraints also can be rewritten equivalently as the 
second-order cone constraints 
\begin{align*}
  w_{i} + x_{i} \ge 
  \begin{Vmatrix}
    \begin{bmatrix}
      w_{i} - x_{i} \\
      2\sqrt{c_{i}/E} q_{i} \\
    \end{bmatrix}
  \end{Vmatrix}
  , 
  \quad i=1,\dots,m . 
\end{align*}
Thus, the conventional compliance minimization, 
\eqref{P.compliance.SOCP.1}, can be recast as SOCP 
\citep[Section~3.4.3]{BtN01}; see also \cite{Koc17,Kan16}. 

We are now in position to consider the upper bound constraint on the 
number of nodes in a truss design. 
Let $n$ denote the specified upper bound. 
For the $j$th node $(j=1,\dots,l)$, define 
$I(j) \subseteq \{ 1,\dots,m\}$ as the set of indices of the members 
connected to node $j$. 
For example, in the case of \reffig{fig:gs2x1} we have  
$I(j) = \{ 1,2,7,10,11 \}$. 
Define $z_{j}$ $(j=1,\dots,l)$ by 
\begin{align}
  z_{j} = \sum_{i \in I(j)} x_{i} 
  \label{eq.node.1}
\end{align}
to see that the number of nodes becomes equal to $\| \bi{z} \|_{0}$. 
For notational simplicity, we write \eqref{eq.node.1} as 
\begin{align*}
  \bi{z} = Z \bi{x}
\end{align*}
with a constant matrix $Z \in \Re^{l \times m}$. 
The upshot is that the compliance minimization under the upper bound 
constraint for  the number of existing nodes is formulated as follows: 
\begin{subequations}\label{P.compliance.node.def}%
  \begin{alignat}{3}
    & \MIN_{\bi{x},\bi{z}}
    &{\quad}& 
    \pi(\bi{x}) \\
    & \ST && 
    \bi{x} \ge  \bi{0} , \\
    & && 
    \bi{c}^{\top} \bi{x} \le V , \\
    & && 
    \bi{z} = Z \bi{x} , \\
    & && 
    \| \bi{z} \|_{0} \le n . 
    \label{P.compliance.node.def.5}
  \end{alignat}
\end{subequations}
As mentioned above, the conventional compliance minimization in 
\eqref{P.compliance.1} can be recast as SOCP. 
Therefore, problem \eqref{P.compliance.node.def} can be reduced to 
cardinality-constrained SOCP. 
In Section~\ref{sec:node.integer}, we present its MISOCP reformulation.

\subsection{MISOCP formulation}
\label{sec:node.integer}

In this section, we show that problem \eqref{P.compliance.node.def} can 
be recast as MISOCP. 

For node $j$ $(j=1,\dots,l)$, we introduce a new variable, 
$s_{j} \in \{ 0,1 \}$, to indicate whether the node vanishes ($s_{j}=0$) 
or exists ($s_{j}=1$). 
The relation between $s_{j}$ and $z_{j}$ can be given as 
\begin{align*}
  0 \le z_{j} \le M s_{j} , 
\end{align*}
where $M>0$ is a sufficiently large constant. 
The upper bound constraint for the number of existing nodes is written 
in terms of $s_{1},\dots,s_{l}$ as 
\begin{align*}
  \sum_{j=1}^{l} s_{j} \le n . 
\end{align*}
This observation, in conjunction with the SOCP reformulation of problem 
\eqref{P.compliance.1}, concludes that problem 
\eqref{P.compliance.node.def} is reduced to the following MISOCP: 
\begin{subequations}\label{P.mixed-integer.SOCP.1}%
  \begin{alignat}{3}
    & \MIN_{\bi{x}, \bi{q}, \bi{w}, \bi{z}, \bi{s}}
    &{\quad}& 
    \sum_{i=1}^{m} w_{i} \\
    & \ST && 
    w_{i} + x_{i} \ge 
    \begin{Vmatrix}
      \begin{bmatrix}
        w_{i} - x_{i} \\
        2\sqrt{c_{i} / E} q_{i} \\
      \end{bmatrix}
    \end{Vmatrix}
    ,  \quad i=1,\dots,m , \\
    & && 
    \sum_{i=1}^{m} q_{i} \bi{b}_{i} = \bi{p} , \\
    & && 
    \bi{c}^{\top} \bi{x} \le V  , \\
    & && 
    \bi{z} = Z \bi{x} , \\
    & && 
    \bi{z} \le M \bi{s} , 
    \label{P.mixed-integer.SOCP.1.M} \\
    & && 
    \sum_{j=1}^{l} s_{j} \le n , \\
    & && 
    \bi{s} \in \{ 0,1 \}^{l} . 
  \end{alignat}
\end{subequations}
Here, optimization variables are 
$\bi{x} \in \Re^{m}$, $\bi{q}\in \Re^{m}$, $\bi{w} \in \Re^{m}$, 
$\bi{z} \in \Re^{l}$, and $\bi{s} \in \Re^{l}$. 
Although problem \eqref{P.mixed-integer.SOCP.1} is a fairly 
straightforward extension of the existing SOCP formulation for 
problem \eqref{P.compliance.1}, it cannot be found in literature to the 
best of the authors' knowledge.

\section{Simple heuristic based on alternating direction method of multipliers}
\label{sec:admm}

In this section, we present an ADMM as a heuristic for problem 
\eqref{P.compliance.node.def}. 

For notational simplicity, define $F \subseteq \Re^{m}$ and 
$G \subseteq \Re^{l}$ by 
\begin{align*}
  F &= 
  \{ \bi{x} \in \Re^{m} 
  \mid
  \bi{x} \ge  \bi{0} , \
  \bi{c}^{\top} \bi{x} \le V 
  \} , \\
  G &= \{ \bi{z} \in \Re^{l}
  \mid \| \bi{z} \|_{0} \le n \} . 
\end{align*}
We see that problem \eqref{P.compliance.node.def} can be written as 
follows: 
\begin{subequations}\label{P.node.1}%
  \begin{alignat}{3}
    & \MIN
    &{\quad}& 
    \pi(\bi{x}) + \delta_{F}(\bi{x}) + \delta_{G}(\bi{z}) \\
    & \ST && 
    Z \bi{x} - \bi{z} = \bi{0} . 
  \end{alignat}
\end{subequations}

The augmented Lagrangian for problem \eqref{P.node.1} is formulated as 
\begin{align}
  L_{\rho}(\bi{x},\bi{z},\bi{y})
  = \pi(\bi{x}) + \delta_{F}(\bi{x}) + \delta_{G}(\bi{z}) 
  + \bi{y}^{\top} (Z \bi{x} - \bi{z}) 
  + \frac{\rho}{2} \| Z \bi{x} - \bi{z} \|^{2} , 
  \label{eq.node.Lagrange.y}
\end{align}
where $\rho > 0$ is the penalty parameter, and $\bi{y} \in \Re^{l}$ is 
the Lagrange multiplier. 
Let $\bi{v}=\bi{y}/\rho$ to see that \eqref{eq.node.Lagrange.y} is 
reduced to 
\begin{align*}
  \tilde{L}_{\rho}(\bi{x},\bi{z},\bi{v})
  = \pi(\bi{x}) + \delta_{F}(\bi{x}) + \delta_{G}(\bi{z}) 
  + \frac{\rho}{2} \| Z \bi{x} - \bi{z} + \bi{v} \|^{2}
  - \frac{\rho}{2} \| \bi{v} \|^{2} . 
\end{align*}
Using $\tilde{L}_{\rho}$, we can write the iterations of ADMM in the 
scaled form as 
\begin{align}
  \bi{x}^{k+1} 
  &:= \argmin_{\bi{x}} 
  \Bigl\{ \pi(\bi{x}) + \delta_{F}(\bi{x}) 
  + \frac{\rho}{2} \| Z \bi{x} - \bi{z}^{k} + \bi{v}^{k} \|^{2} \Bigr\} , 
  \label{eq.node.ADMM.1.1} \\
  \bi{z}^{k+1} 
  &:= \argmin_{\bi{z}} 
  \Bigl\{
  \delta_{G}(\bi{z}) 
  + \frac{\rho}{2} \| Z \bi{x}^{k+1} - \bi{z} + \bi{v}^{k} \|^{2} \Bigr\} ,
  \label{eq.node.ADMM.1.2} \\
  \bi{v}^{k+1}
  &:= \bi{v}^{k} + Z \bi{x}^{k+1} - \bi{z}^{k+1} . 
  \label{eq.node.ADMM.1.3}
\end{align}

The first step of ADMM in \eqref{eq.node.ADMM.1.1} means that we let 
$\bi{x}^{k+1}$ be an optimal solution of the following convex 
optimization problem: 
\begin{subequations}\label{P.node.2}%
  \begin{alignat}{3}
    & \MIN
    &{\quad}& 
    \pi(\bi{x}) 
    + \frac{\rho}{2} \| Z \bi{x} - \bi{z}^{k} + \bi{v}^{k} \|^{2}  \\
    & \ST && 
    \bi{x} \ge  \bi{0} , \\
    & && 
    \bi{c}^{\top} \bi{x} \le V . 
  \end{alignat}
\end{subequations}
This problem can be recast as SOCP. 
To see this, using the SOCP formulation of problem 
\eqref{P.compliance.1}, we rewrite problem \eqref{P.node.2} as follows: 
\begin{subequations}\label{P.node.SOCP}%
  \begin{alignat}{3}
    & \MIN
    &{\quad}& 
    \sum_{i=1}^{m} w_{i} + \frac{\rho}{2} t \\
    & \ST && 
    t \ge \| Z \bi{x} - \bi{z}^{k} + \bi{v}^{k} \|^{2} , 
    \label{P.node.SOCP.2} \\
    &  && 
    w_{i} + x_{i} \ge 
    \begin{Vmatrix}
      \begin{bmatrix}
        w_{i} - x_{i} \\
        2\sqrt{c_{i} / E} q_{i} \\
      \end{bmatrix}
    \end{Vmatrix}
    ,   \quad i=1,\dots,m , \\
    & && 
    \sum_{i=1}^{m} q_{i} \bi{b}_{i} = \bi{p} , \\
    & && 
    \bi{c}^{\top} \bi{x} \le V , 
  \end{alignat}
\end{subequations}
where $t \in \Re$ is an auxiliary variable. 
Since constraint \eqref{P.node.SOCP.2} is a rotated second-order cone 
constraint\footnote{ 
It can also be rewritten as 
\begin{align*}
  t + 1 \ge 
  \begin{Vmatrix}
    \begin{bmatrix}
      t - 1 \\
      2(Z \bi{x} - \bi{z}^{k} + \bi{v}^{k}) \\
    \end{bmatrix}
  \end{Vmatrix}
  , 
\end{align*}
which is a second-order cone constraint. }
\begin{align*}
  (Z \bi{x} - \bi{z}^{k} + \bi{v}^{k}, t, 1) \in \KC^{l+2} , 
\end{align*}
problem \eqref{P.node.SOCP} is an SOCP. 
We adopt a primal-dual interior-point method for solving this problem. 
Next, the second step of ADMM in \eqref{eq.node.ADMM.1.2} can be 
written as 
\begin{align}
  \bi{z}^{k+1} 
  \in \Pi_{G}(Z \bi{x}^{k+1} + \bi{v}^{k}) , 
  \label{eq.projection.node.1}
\end{align}
where $\Pi_{G}$ is the projection onto $G$.\footnote{
Since $G$ is nonconvex, the projection of a point onto $G$ is not 
necessarily unique. 
} 
We can compute \eqref{eq.projection.node.1} easily \cite[Chap.~9]{BPCPE10}; 
for a point $\bi{z} \in \Re^{l}$, $\Pi_{G}(\bi{z})$ keeps the $n$ 
largest magnitude components of $\bi{z}$ and zeros out the other components. 
In this way, each step of ADMM in 
\eqref{eq.node.ADMM.1.1}, \eqref{eq.node.ADMM.1.2}, and 
\eqref{eq.node.ADMM.1.3} can be carried out very easily.

\section{On overlapping members}
\label{sec:overlap}

Unlike the conventional compliance minimization of a truss, overlapping 
members in a ground structure are not redundant for the optimization 
problem considered in this paper. 
This section explains the treatment of overlapping members. 

\begin{figure}[tp]
  \centering
  \begin{subfigure}[b]{0.40\textwidth}
    \centering
    \includegraphics[scale=0.40]{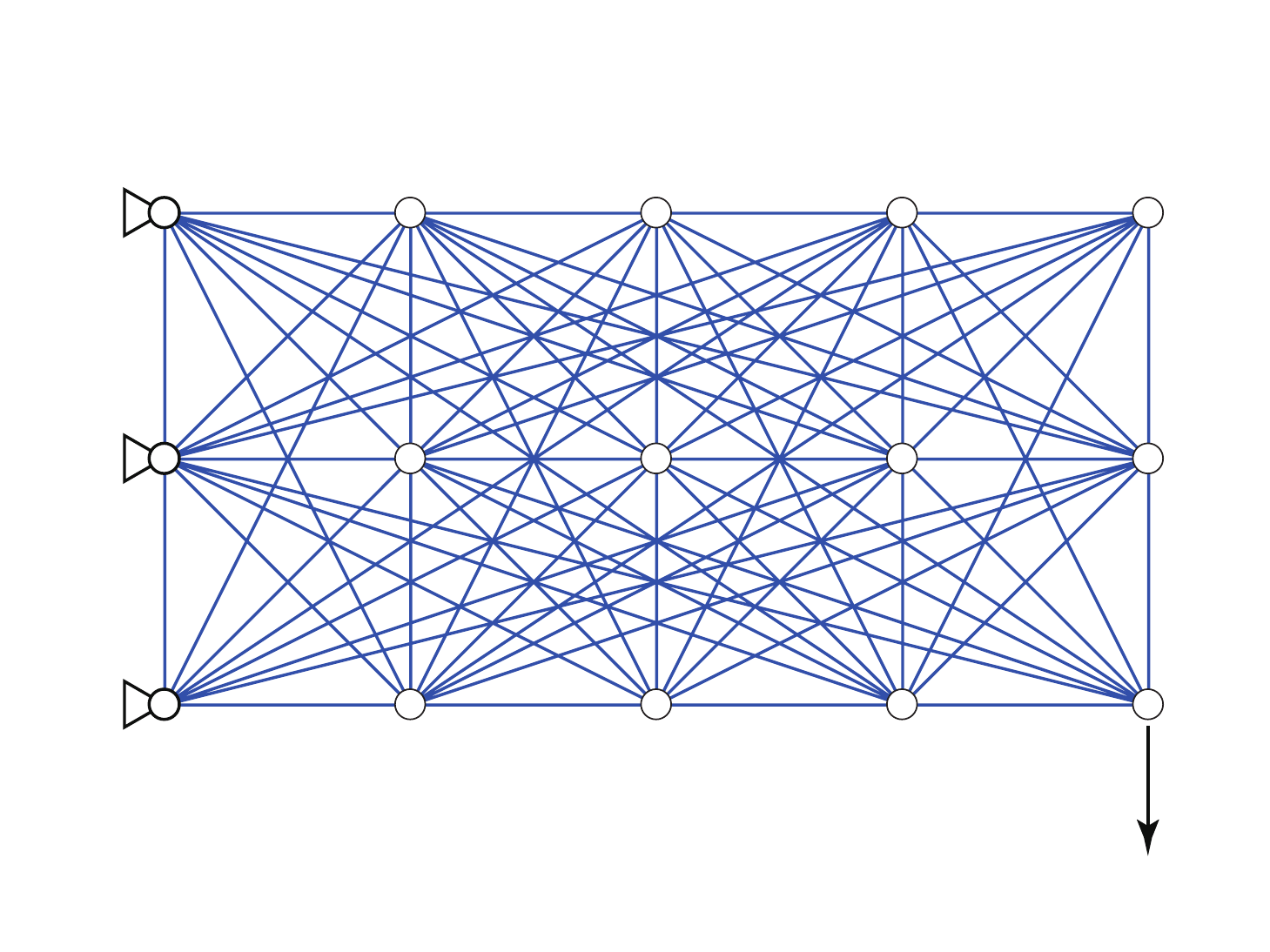}
    \caption{}
    \label{fig:x4_y2_initial}
  \end{subfigure}
  \par\medskip
  \begin{subfigure}[b]{0.40\textwidth}
    \centering
    \includegraphics[scale=0.40]{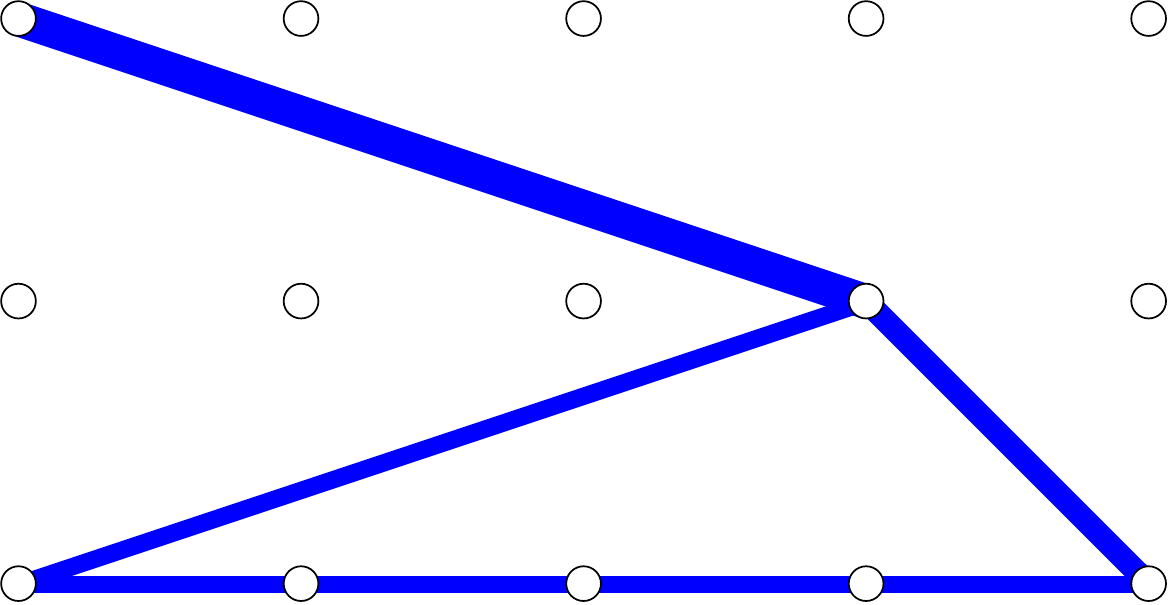}
    \caption{}
    \label{fig:x4_y2_all_node}
  \end{subfigure}
  \hfill
  \begin{subfigure}[b]{0.40\textwidth}
    \centering
    \includegraphics[scale=0.40]{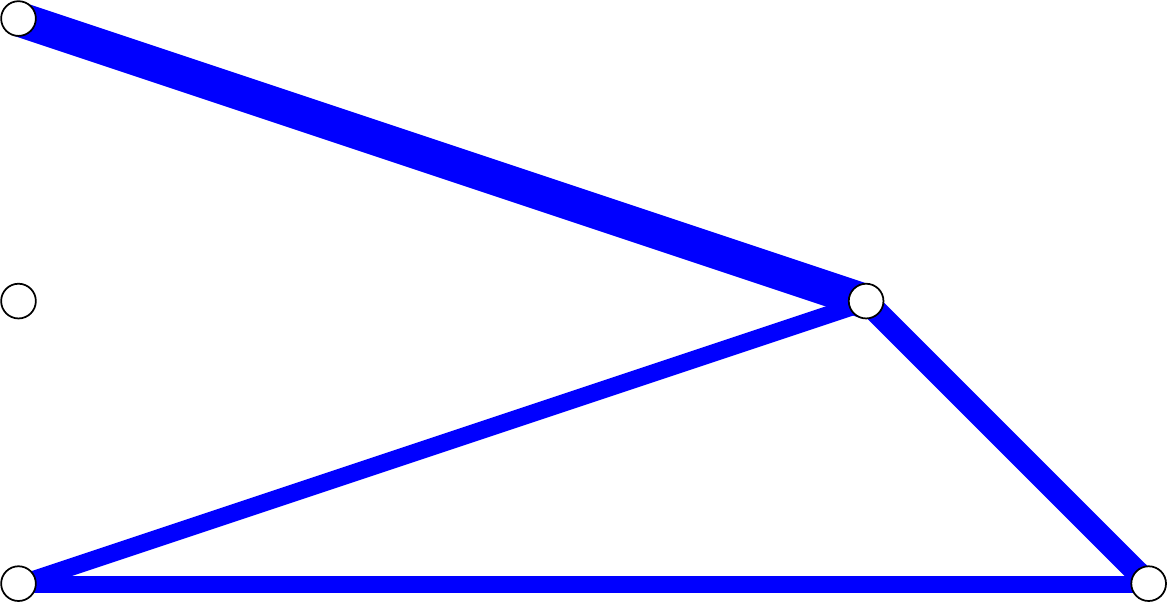}
    \caption{}
    \label{fig:x4_y2_post}
  \end{subfigure}
  \caption{An example of truss topology optimization and hinge cancellation. 
  \subref{fig:x4_y2_initial} The problem setting; 
  \subref{fig:x4_y2_all_node} the optimal solution; and 
  \subref{fig:x4_y2_post} the final design after hinge cancellation. }
  \label{fig:x4_y2}
\end{figure}

We begin by reviewing that overlapping members in a ground structure is 
redundant for the conventional compliance minimization of a truss. 
For example, consider the ground structure shown in 
\reffig{fig:x4_y2_initial}. 
Here, any two nodes are connected by a member, but overlapping of 
members is avoided by removing the longer member when two members overlap. 
The leftmost nodes are pin-supported. 
The vertical external force is applied to the bottom rightmost node. 
\reffig{fig:x4_y2_all_node} shows the optimal solution of the 
compliance minimization, i.e., problem \eqref{P.compliance.1}. 
This solution has four horizontal consecutive members 
that are connected by nodes supported only in the direction of those 
members. 
A sequence of such members is sometimes called a {\em chain\/} \cite{Ach99}. 
In this example, without changing the objective value, 
we can remove three intermediate nodes to 
replace the chain with a single longer member. 
This procedure is called the {\em hinge cancellation\/} \cite{Ach99,Roz96}. 
As a result of hinge cancellation, we obtain the final truss design 
shown in \reffig{fig:x4_y2_post}. 
Thus, longer overlapping members, like the horizontal member in 
\reffig{fig:x4_y2_post}, are unnecessary to a ground structure. 
In contrast, when we consider a constraint on the number 
of nodes, the optimal solution depends on existence of overlapping 
members in a ground structure. 
For example, the truss in \reffig{fig:x4_y2_all_node} has five free 
nodes, while the one in \reffig{fig:x4_y2_post} has two free nodes. 
Thus, the hinge cancellation can possibly change the feasibility of the 
cardinality constraint and, hence, overlapping members in a ground 
structure are not redundant.\footnote{
Such non-redundancy of overlapping members is also known for 
truss topology optimization considering, e.g., the self-weight load 
\cite{BBtZ94,KY17} and the member buckling constraints \cite{Mel14,GCO05}. 
} 

When we consider a ground structure with some overlapping members, 
existence of overlapping members in an obtained solution is not allowed 
from a practical point of view. 
The method proposed in Section~\ref{sec:admm} does not consider 
explicitly the constraint prohibiting presence of overlapping members. 
Nevertheless, in practice, a solution obtained by the proposed method 
often has no overlapping members, as illustrated through numerical 
experiments in Section~\ref{sec:admm}. 

Within the framework of MISOCP, we can explicitly 
incorporate the constraints prohibiting the presence of mutually 
overlapping members in a truss design. 
To do this, besides $\bi{s} \in \{0,1 \}^{l}$ in 
Section~\ref{sec:node.integer}, we use extra 0-1 variables 
$\bi{t}\in\{ 0,1 \}^{m}$ to indicate whether each member vanishes or 
exists. 
Namely, $t_{i}=0$ means that member $i$ is removed, while $t_{i}=1$ 
means that member $i$ exists. 
The relation between $t_{i}$ and $x_{i}$ is given by 
\begin{align*}
  0 \le x_{i} \le M t_{i} 
\end{align*}
where $M>0$ is a sufficiently large constant. 
Recall that $I(j)$ denotes the set of indices of the members connected to 
node $j$; see Section~\ref{sec:admm}. 
The relation between $t_{i}$ $(i \in I(j))$ and $s_{j}$ is given by 
\begin{align*}
  t_{i} \le s_{j} , 
  \quad \forall i \in I(j) . 
\end{align*}
Let $D$ denote the set of pairs of indices of the members that mutually 
overlap. 
Namely, $(i_{1},i_{2}) \in D$ means that member $i_{1}$ and member 
$i_{2}$ cannot exist simultaneously. 
This constraint is written as 
\begin{align*}
  t_{i_{1}} + t_{i_{2}} \le 1 , 
  \quad \forall (i_{1},i_{2}) \in D . 
\end{align*}
The upshot is that the truss topology 
optimization problem can be formulated as 
the following MISOCP: 
\begin{subequations}\label{P.mixed-integer.SOCP.2}%
  \begin{alignat}{3}
    & \MIN
    &{\quad}& 
    \sum_{i=1}^{m} w_{i} \\
    & \ST && 
    w_{i} + x_{i} \ge 
    \begin{Vmatrix}
      \begin{bmatrix}
        w_{i} - x_{i} \\
        2\sqrt{c_{i} / E} q_{i} \\
      \end{bmatrix}
    \end{Vmatrix}
    , 
    \quad i=1,\dots,m , \\
    & && 
    \sum_{i=1}^{m} q_{i} \bi{b}_{i} = \bi{p} , \\
    & && 
    \bi{c}^{\top} \bi{x} \le V  , \\
    & && 
    \bi{x} \le M \bi{t} , 
    \label{P.mixed-integer.SOCP.2.5} \\
    & && 
    t_{i} \le s_{j} \ (\forall i \in I(j)), 
    \quad j=1,\dots,l, \\
    & && 
    \bi{s} \ge \bi{0} , \\
    & && 
    \sum_{j=1}^{l} s_{j} \le n , \\
    & && 
    t_{i_{1}} + t_{i_{2}} \le 1 , 
    \quad \forall (i_{1},i_{2}) \in D , 
    \label{P.mixed-integer.SOCP.2.9} \\
    & && 
    \bi{t} \in \{ 0,1 \}^{m} . 
  \end{alignat}
\end{subequations}
It is worth noting that the $0$-$1$ constraints on 
$s_{1},\dots,s_{k}$ can be omitted.

\section{Numerical experiments}
\label{sec:ex}

In this section, we report numerical experiments on the method presented 
in Section~\ref{sec:admm}. 
In Section~\ref{sec:ex.implementation}, we describe the details of 
implementation of the algorithm and the problem settings of the 
numerical experiments. 
The computational results of the proposed ADMM approach, together with 
the comparison with the MISOCP approach, are presented in 
Sections~\ref{sec:ex.n4}, \ref{sec:ex.n5}, and \ref{sec:ex.7m}. 
Empirical evidences of our stopping criterion and selection of initial 
points are presented in Sections~\ref{sec:ex.stopping} and 
\ref{sec:ex.initial}, respectively. 
Section~\ref{sec:ex.robust} presents application of the proposed method 
to robust truss optimization, which is recast as mixed-integer 
semidefinite programming. 

\subsection{Implementation and problem settings}
\label{sec:ex.implementation}

At each iteration of the proposed method, we solved problem 
\eqref{P.node.SOCP} by using CVX ver.~2.1, 
a MATLAB package for specifying and solving convex optimization 
problems \cite{GB08,CVX}. 
As a solver, we used SDPT3 ver.~4.0 \cite{TTT03} on 
MATLAB ver.~9.1.0. 
The \texttt{cvx\_precision} of CVX is set to \texttt{best}, which means 
that the solver continues as far as it can make progress~\cite{CVX}. 
For comparison, we solved the MISOCP problem in 
\eqref{P.mixed-integer.SOCP.1} with a global optimization approach. 
The value of $M$ in constraint 
\eqref{P.mixed-integer.SOCP.1.M} is set to $1.0\times 10^5$ in $\mathrm{m}$.\footnote{%
Through our preliminary numerical experiments it was found that the 
computational cost required by MOSEK does not change drastically 
depending on the value of $M$. }  
We used PICOS ver.~1.1.2, a Python interface to diverse optimization 
solvers \cite{PICOS}. 
MOSEK ver.~8.0.1 \cite{ART03} was used as the solver. 
Computation was carried out on two $3.2\,\mathrm{GHz}$ Intel Xeon E5-2667 v4 
processors with $256\,\mathrm{GB}$ RAM. 

In practice, we slightly modify the original version of ADMM introduced 
in Section~\ref{sec:preliminary} so that the penalty parameter in the 
augmented Lagrangian is gradually increased. 
Specifically, $\rho$ in subproblem \eqref{P.node.SOCP} is given by 
\begin{align*}
  \rho_{k+1}:= \min\{ \mu \rho_{k} , \rho_{\rr{max}} \} ,  
\end{align*}
where $\mu$ ($>1$) and $\rho_{\rr{max}}$ $(>\rho_{0})$ are constants. 
In the following, we set $\mu=1.5$, 
$\rho_{0}=1$, and $\rho_{\rr{max}}=10^{6}$. 
Define $J_{0}^{k} \subseteq \{ 1,\dots,l \}$ by 
\begin{align*}
  J_{0}^{k} = \{ j \in \{ 1,\dots,l \} \mid z_{j}^{k} \le \epsilon \} , 
\end{align*}
where we set $\epsilon = 0.1\,\mathrm{mm^{2}}$. 
We terminate the ADMM when 
\begin{align*}
  l - |J_{0}^{k}| \le n 
\end{align*}
is satisfied. 
Then we solve problem \eqref{P.compliance.1} with the additional 
constraints 
\begin{align*}
  \sum_{i \in I(j)} x_{i} = 0 , 
  \quad \forall j \in J_{0}^{k}
\end{align*}
to generate the final output. 
As for the initial point for the ADMM, we examine two cases: 
\begin{itemize}
  \item Initial point (A): 
        $\bi{z}^{0} := Z \bi{x}^{0}$ and 
        $\bi{v}^{0} := \bi{0}$, where 
        $\bi{x}^{0}$ is an optimal solution of problem \eqref{P.compliance.1}. 
  \item Initial point (B): 
        $\bi{z}^{0} := Z \bi{x}^{0}$ and 
        $\bi{v}^{0} := \bi{0}$ with 
        $\bi{x}^{(0)} := (V / \bi{c}^{\top} \bi{1}) \bi{1}$. 
\end{itemize}
It should be clear that only $\bi{z}^{0}$ and $\bi{v}^{0}$ are used as 
input data of the ADMM; $\bi{x}^{0}$ is not required as input. 

Consider the problem setting shown in \reffig{fig:gs5x2}. 
The nodes are aligned on a $1\,\mathrm{m} \times 1\,\mathrm{m}$ grid. 
We vary the values of $N_{X}$ and $N_{Y}$ to generate problem instances 
with diverse sizes. 
The number of free nodes in this ground structure is $N_{X}(N_{Y}+1)$. 
The members in a ground structure are generated as follows: 
We first consider all possible members such that any two nodes are 
connected by a member. 
Then we remove members that are longer than a specified value, 
$5\,\mathrm{m}$ in Sections~\ref{sec:ex.n4} and \ref{sec:ex.n5} 
while $7\,\mathrm{m}$ in Section~\ref{sec:ex.7m}. 
It is worth noting that the ground structure retains overlapping members. 

In the following examples, the Young modulus is $E=20\,\mathrm{GPa}$, 
and the specified upper bound for the structural volume is 
$V=2 N_{X} N_{Y} \times 10^{5}\,\mathrm{mm^{3}}$. 
As for $\bi{p}$, the external vertical force of $100\,\mathrm{kN}$ is 
applied to the bottom rightmost node. 
We consider $n$ as the upper bound for the number of free nodes. 
In other words, the number of supports is not restricted in the 
following examples, and $l$ in the previous sections denotes the number 
of free nodes of a ground structure.

\subsection{Example (I)}
\label{sec:ex.n4}

\begin{table}[bp]
  \centering
  \caption{Characteristics of the problem instances for the numerical experiments.}
  \label{tab:ex.I.data}
  \begin{tabular}{lrrrr}
    \toprule
    $(N_{X},N_{Y})$ & $m$ & $d$ & $\hat{w}$ (J) & {\#}free nodes \\
    \midrule
    (5,2) & 147 & 30 & $12100.00$ & 9 \\
    (5,3) & 264 & 40 & $5007.41$  & 7 \\
    (5,4) & 411 & 50 & $2812.50$  & 5 \\
    \midrule
    (8,2) & 273 & 48 & $34515.63$ & 15 \\
    (9,2) & 315 & 54 & $45125.00$ & 8 \\
    (8,4) & 750 & 80 & $6937.81$ & 10 \\
    (9,4) & 863 & 90 & $8900.28$ & 10 \\
    (8,6) & 1296 & 112 & $3080.39$ & 10 \\
    (9,6) & 1489 & 126 & $3847.25$ & 12 \\
    \bottomrule
  \end{tabular}
\end{table}

\begin{table}[bp]
  \centering
  \caption{Computational results of example (I).}
  \label{tab:ex.I.result}
  \begin{tabular}{lrrrrrrr}
    \toprule
    & \multicolumn{5}{c}{ADMM} & \multicolumn{2}{c}{MISOCP} \\
    \cmidrule(lr){2-6}
    \cmidrule(l){7-8}
    $(N_{X},N_{Y})$ &  Init.\ sol.\ & $w^{*}$ (J) & $w^{*}/\hat{w}$ & {\#}iter & Time (s) 
    & $\bar{w}$ (J) & Time (s) \\
    \midrule
    (5,2) 
    & $*$ (A) & $12100.00$ & $1.000$ & $5$ & $3.6$ &$12100.00$ &$0.72$ \\
    & $*$ (B) & $12100.00$ & $1.000$ & $3$ & $2.2$ \\
    \midrule
    (5,3) 
    & $*$ (A) & $5007.41$ & $1.000$ & $5$ & $3.2$ &$5007.41$ &$1.69$ \\
    & (B) & $5052.45$ & --- & $3$ & $2.0$ \\
    \midrule
    (5,4) 
    & $*$ (A) & $2812.50$ & $1.000$ & $2$ & $2.6$ &$2812.50$ &$0.52$ \\
    & $*$ (B) & $2812.50$ & $1.000$ & $3$ & $3.8$ \\
    \bottomrule
  \end{tabular}
\end{table}

\begin{figure}[tp]
  \centering
  \includegraphics[scale=0.60]{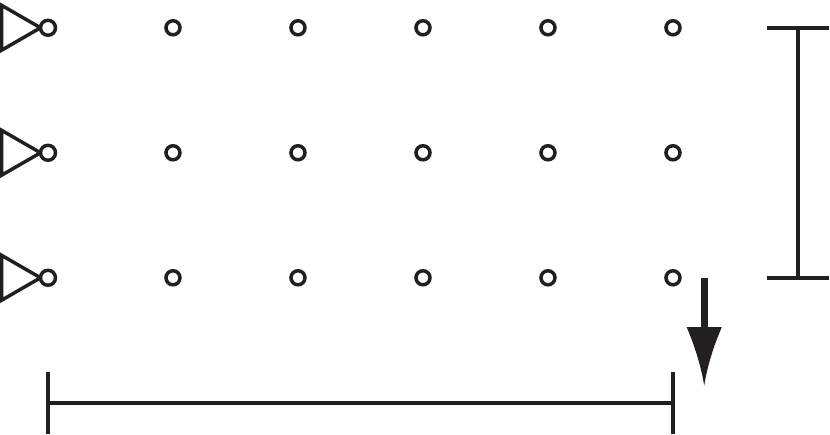}
  \begin{picture}(0,0)
    \put(-150,-50){
    \put(50,44){{\footnotesize $N_{X}${\,}@$1\,\mathrm{m}$}}
    \put(144,94){{\footnotesize $N_{Y}${\,}@$1\,\mathrm{m}$}}
    \put(127,60){{\footnotesize $\bi{p}$}}
    }
  \end{picture}
  \medskip
  \caption{The problem setting for numerical experiments with 
  $(N_{X},N_{Y})=(5,2)$. }
  \label{fig:gs5x2}
\end{figure}
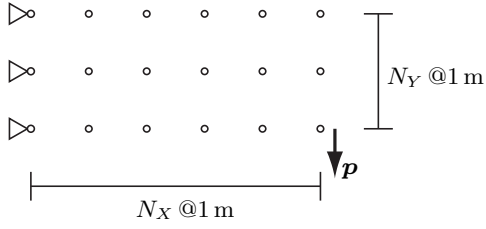

\begin{figure}[tp]
  \centering
  \begin{subfigure}[b]{0.45\textwidth}
    \centering
    \includegraphics[scale=0.40]{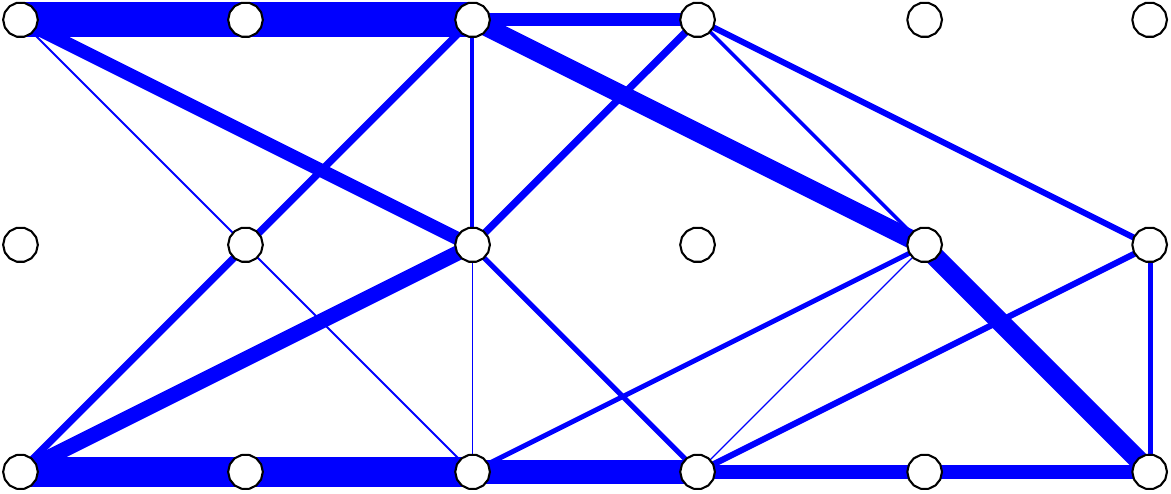}
    \caption{}
    \label{fig:x5_y2_nominal}
  \end{subfigure}
  \hfill
  \begin{subfigure}[b]{0.45\textwidth}
    \centering
    \includegraphics[scale=0.40]{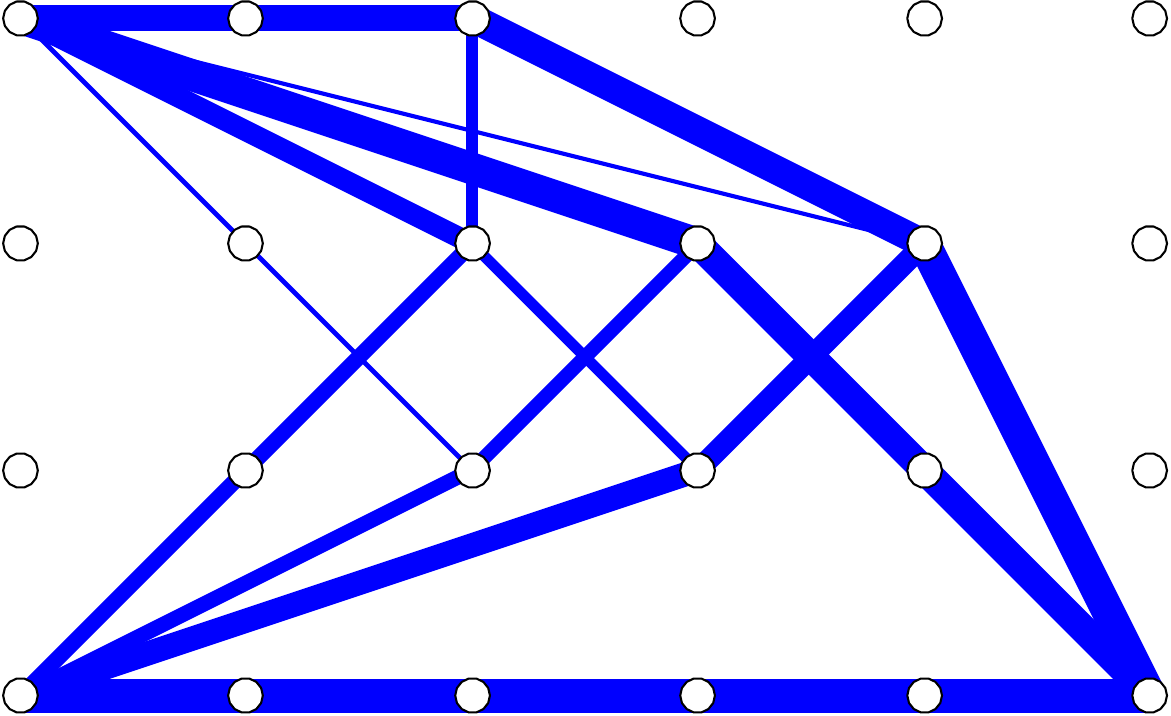}
    \caption{}
    \label{fig:x5_y3_nominal}
  \end{subfigure}
  \par\medskip
  \begin{subfigure}[b]{0.45\textwidth}
    \centering
    \includegraphics[scale=0.40]{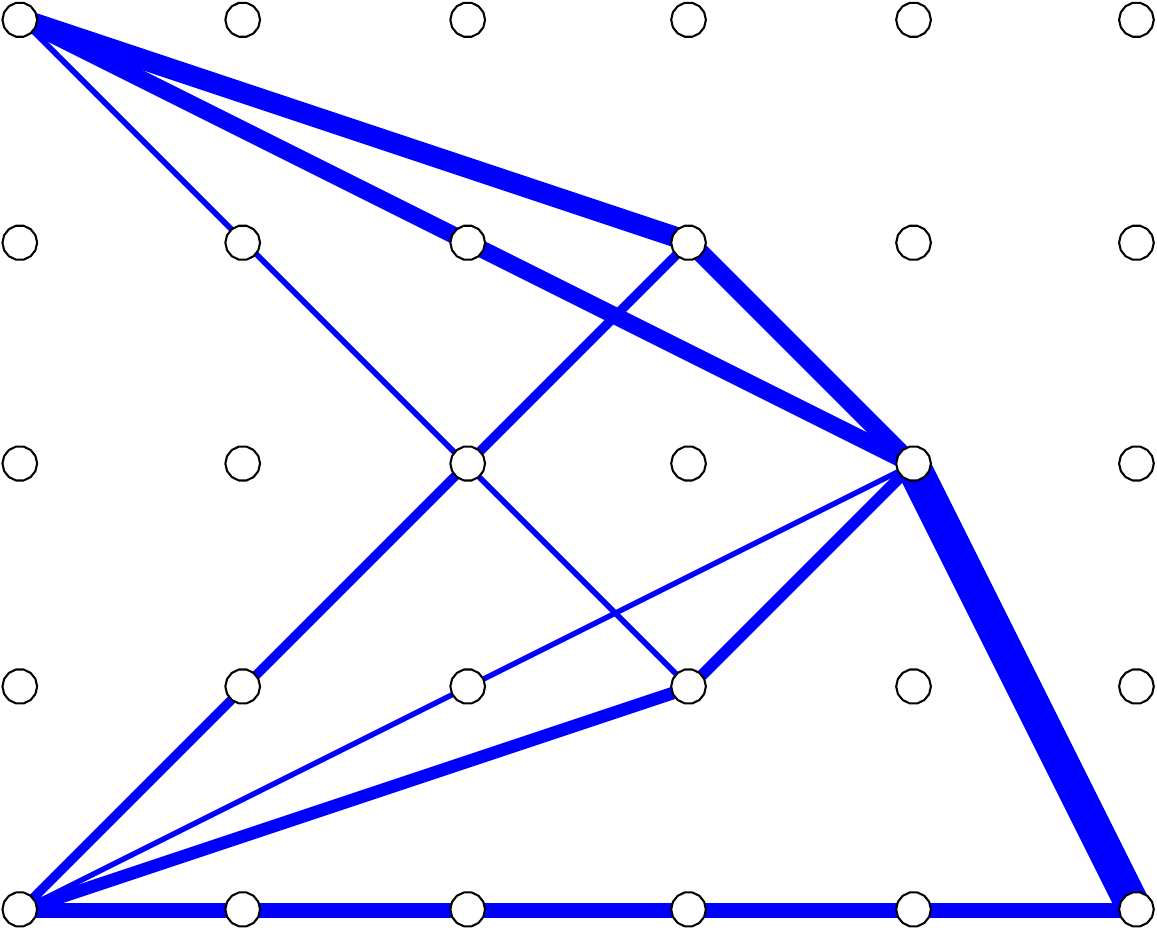}
    \caption{}
    \label{fig:x5_y4_nominal}
  \end{subfigure}
  \caption{Example (I). The optimal solutions of the compliance 
  minimization (without the cardinality constraint). 
  \subref{fig:x5_y2_nominal} $(N_{X},N_{Y})=(5,2)$; 
  \subref{fig:x5_y3_nominal} $(N_{X},N_{Y})=(5,3)$; and 
  \subref{fig:x5_y4_nominal} $(N_{X},N_{Y})=(5,4)$. 
  }
  \label{fig:x5_nominal}
\end{figure}

\begin{figure}[tp]
  \centering
  \begin{subfigure}[b]{0.45\textwidth}
    \centering
    \includegraphics[scale=0.40]{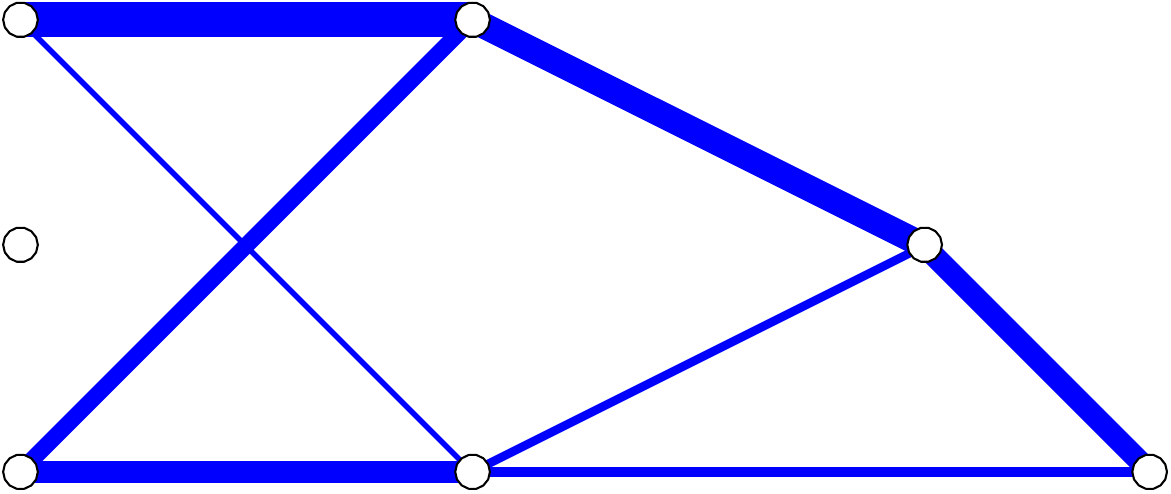}
    \caption{}
    \label{fig:x5_y2_post}
  \end{subfigure}
  \hfill
  \begin{subfigure}[b]{0.45\textwidth}
    \centering
    \includegraphics[scale=0.40]{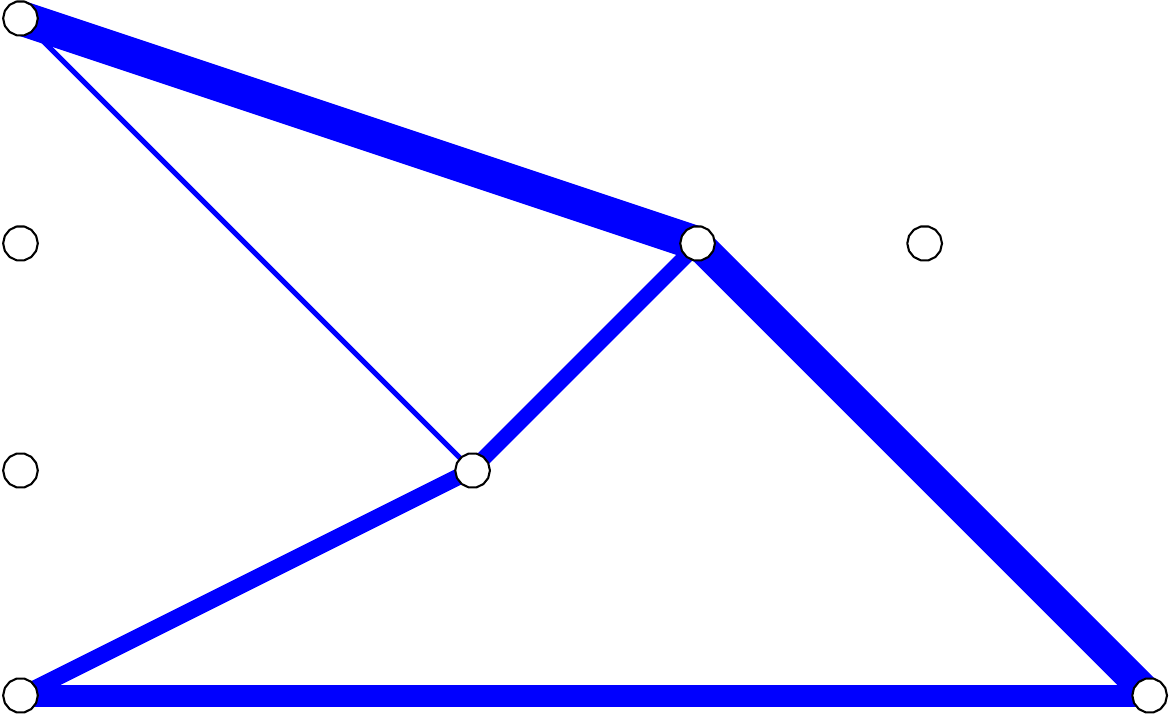}
    \caption{}
    \label{fig:x5_y3_post}
  \end{subfigure}
  \par\medskip
  \begin{subfigure}[b]{0.45\textwidth}
    \centering
    \includegraphics[scale=0.40]{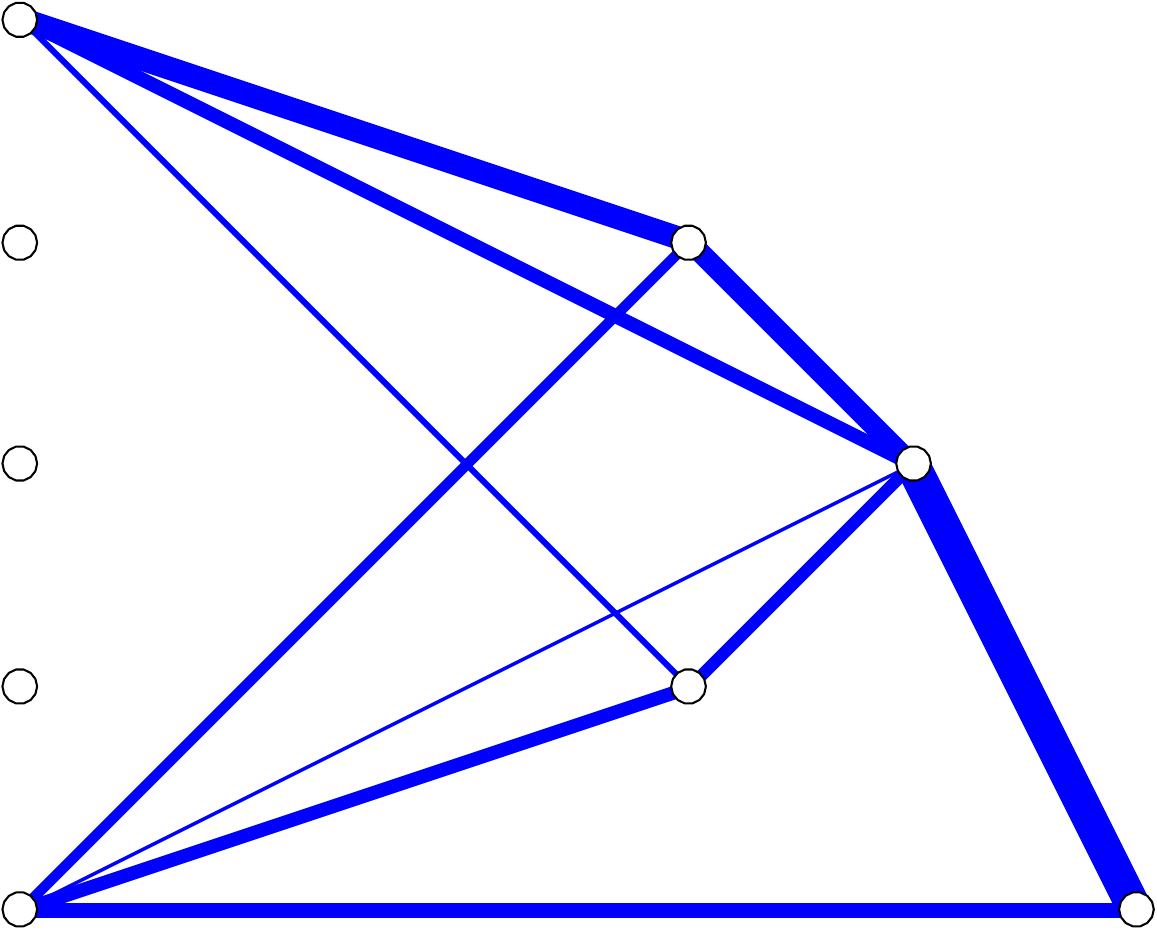}
    \caption{}
    \label{fig:x5_y4_post}
  \end{subfigure}
  \caption{Example (I). The solutions obtained by the proposed method 
  for the compliance minimization with the cardinality constraint ($n=4$). 
  \subref{fig:x5_y2_post} $(N_{X},N_{Y})=(5,2)$; 
  \subref{fig:x5_y3_post} $(N_{X},N_{Y})=(5,3)$; and 
  \subref{fig:x5_y4_post} $(N_{X},N_{Y})=(5,4)$. 
  }
  \label{fig:x5_post}
\end{figure}

\begin{figure}[tp]
  \centering
  \begin{subfigure}[b]{0.45\textwidth}
    \centering
    \includegraphics[scale=0.40]{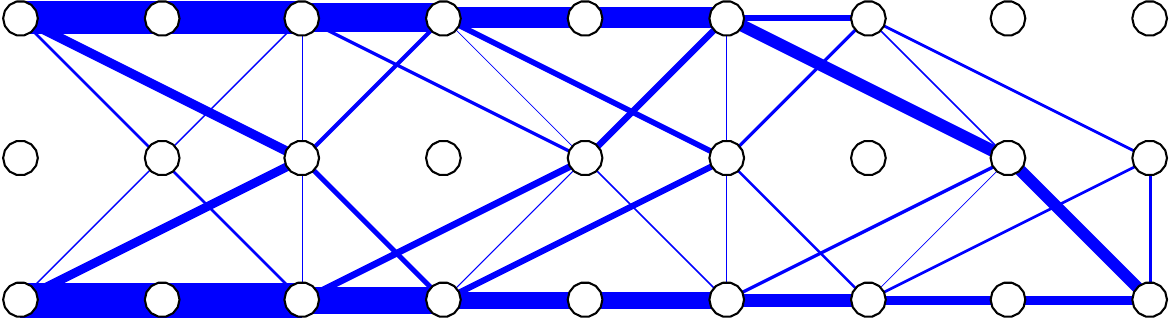}
    \caption{}
    \label{fig:x8_y2_nominal}
  \end{subfigure}
  \hfill
  \begin{subfigure}[b]{0.45\textwidth}
    \centering
    \includegraphics[scale=0.40]{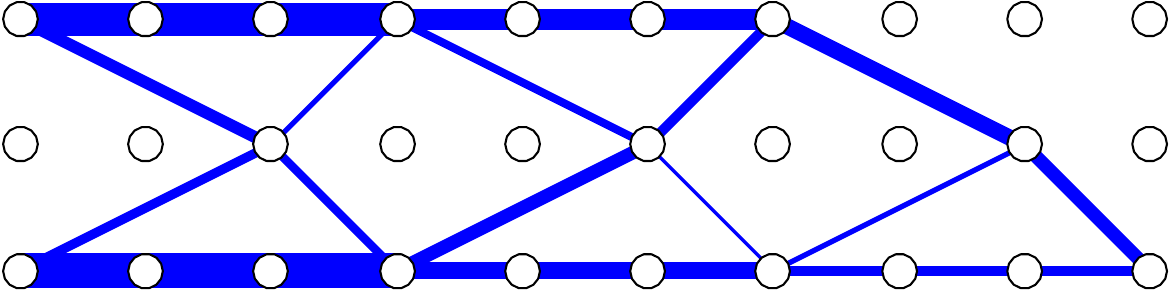}
    \caption{}
    \label{fig:x9_y2_nominal}
  \end{subfigure}
  \par\medskip
  \begin{subfigure}[b]{0.45\textwidth}
    \centering
    \includegraphics[scale=0.40]{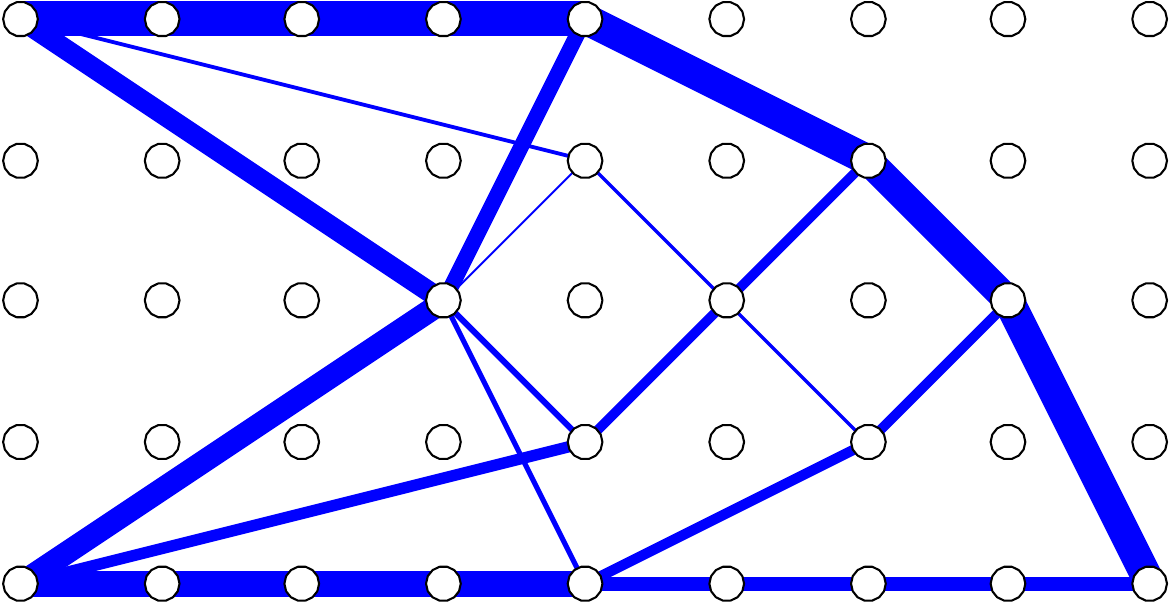}
    \caption{}
    \label{fig:x8_y4_nominal}
  \end{subfigure}
  \hfill
  \begin{subfigure}[b]{0.45\textwidth}
    \centering
    \includegraphics[scale=0.40]{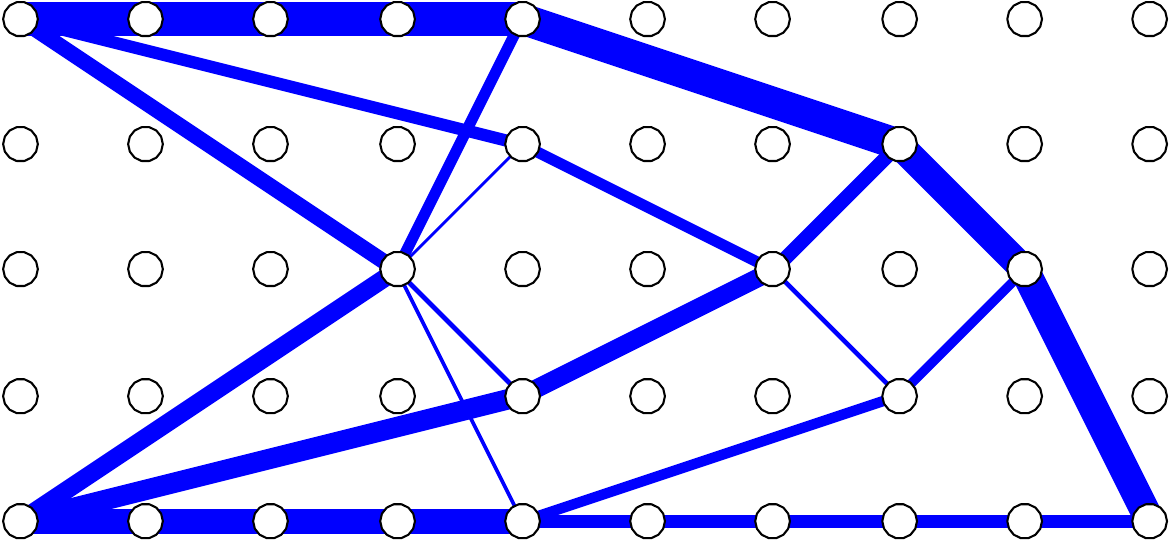}
    \caption{}
    \label{fig:x9_y4_nominal}
  \end{subfigure}
  \par\medskip
  \begin{subfigure}[b]{0.45\textwidth}
    \centering
    \includegraphics[scale=0.40]{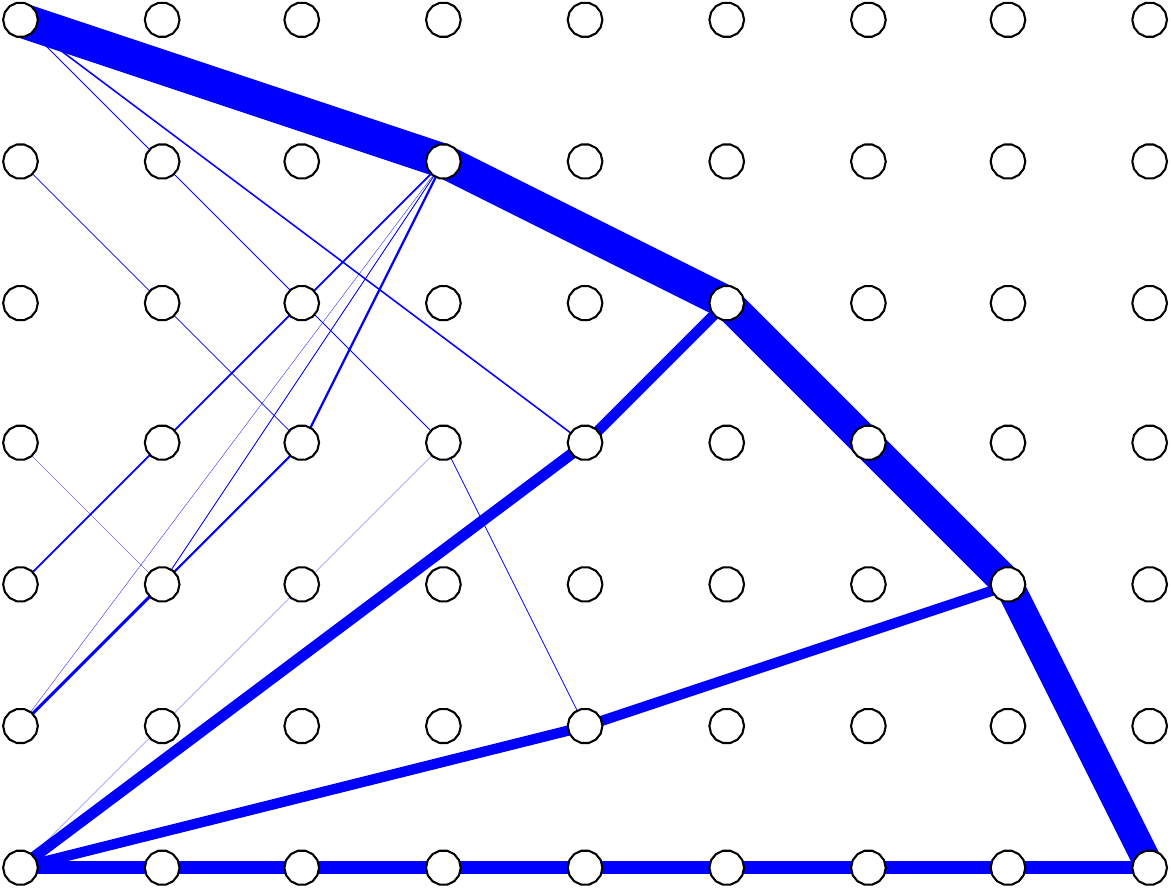}
    \caption{}
    \label{fig:x8_y6_nominal}
  \end{subfigure}
  \hfill
  \begin{subfigure}[b]{0.45\textwidth}
    \centering
    \includegraphics[scale=0.40]{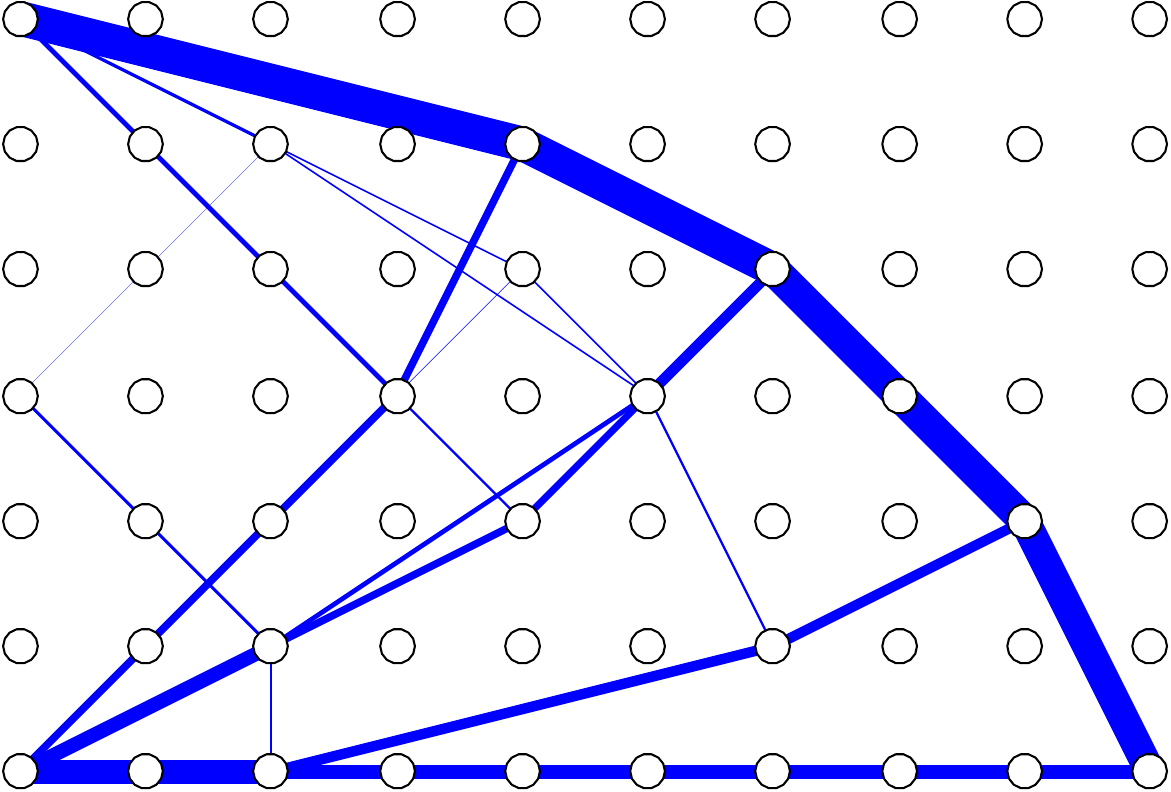}
    \caption{}
    \label{fig:x9_y6_nominal}
  \end{subfigure}
  \caption{Example (II). 
  The optimal solutions of the compliance minimization (without the 
  cardinality constraint). 
  \subref{fig:x8_y2_nominal} $(N_{X},N_{Y})=(8,2)$; 
  \subref{fig:x9_y2_nominal} $(N_{X},N_{Y})=(9,2)$; 
  \subref{fig:x8_y4_nominal} $(N_{X},N_{Y})=(8,4)$; 
  \subref{fig:x9_y4_nominal} $(N_{X},N_{Y})=(9,4)$; 
  \subref{fig:x8_y6_nominal} $(N_{X},N_{Y})=(8,6)$; and 
  \subref{fig:x9_y6_nominal} $(N_{X},N_{Y})=(9,6)$. 
  }
  \label{fig:x8_x9_nominal}
\end{figure}

\begin{figure}[tp]
  \centering
  \begin{subfigure}[b]{0.45\textwidth}
    \centering
    \includegraphics[scale=0.40]{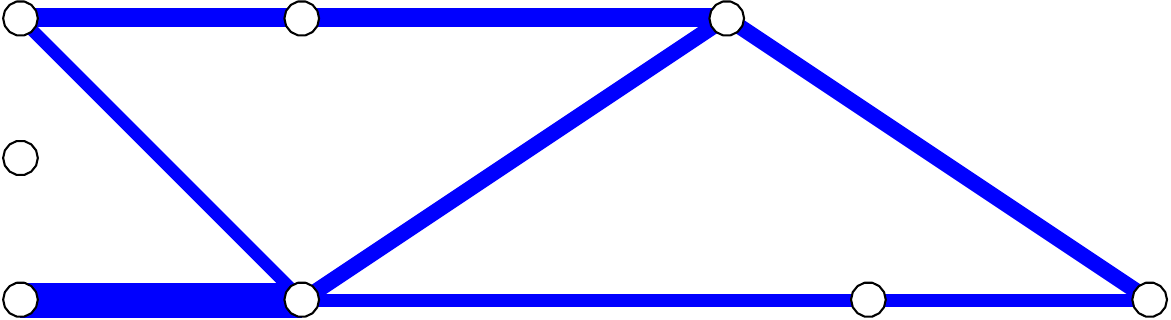}
    \caption{}
    \label{fig:x8_y2_post}
  \end{subfigure}
  \hfill
  \begin{subfigure}[b]{0.45\textwidth}
    \centering
    \includegraphics[scale=0.40]{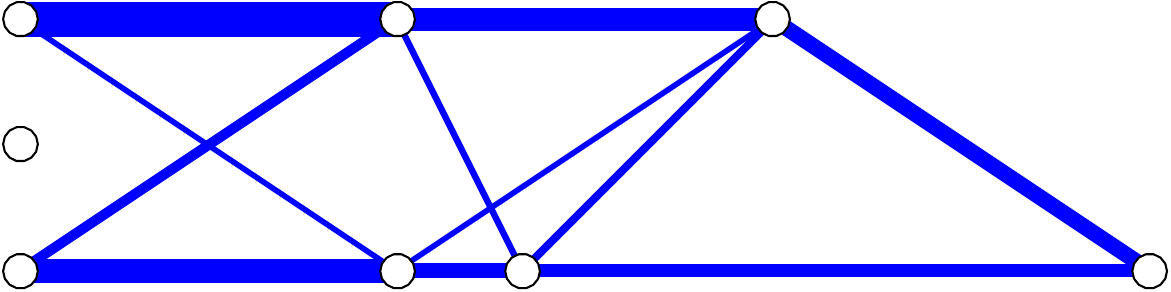}
    \caption{}
    \label{fig:x9_y2_post}
  \end{subfigure}
  \par\medskip
  \begin{subfigure}[b]{0.45\textwidth}
    \centering
    \includegraphics[scale=0.40]{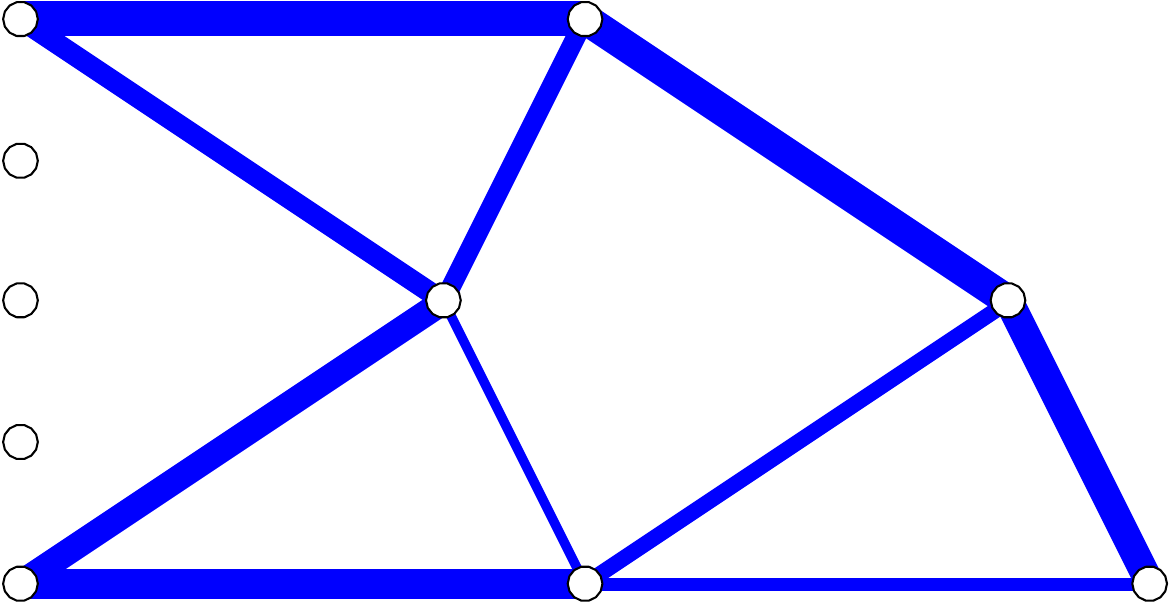}
    \caption{}
    \label{fig:x8_y4_post}
  \end{subfigure}
  \hfill
  \begin{subfigure}[b]{0.45\textwidth}
    \centering
    \includegraphics[scale=0.40]{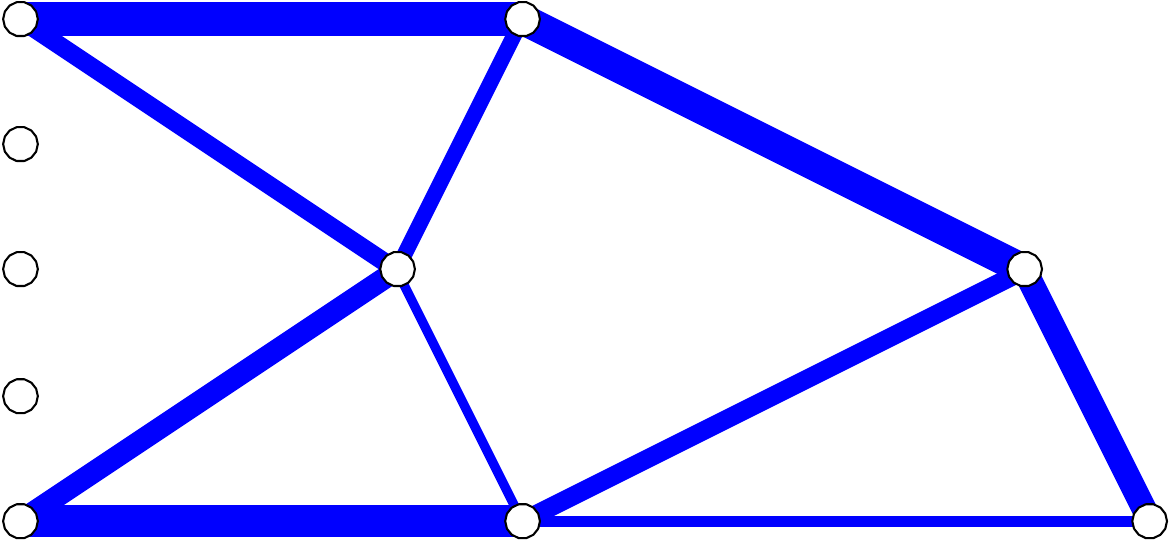}
    \caption{}
    \label{fig:x9_y4_post}
  \end{subfigure}
  \par\medskip
  \begin{subfigure}[b]{0.45\textwidth}
    \centering
    \includegraphics[scale=0.40]{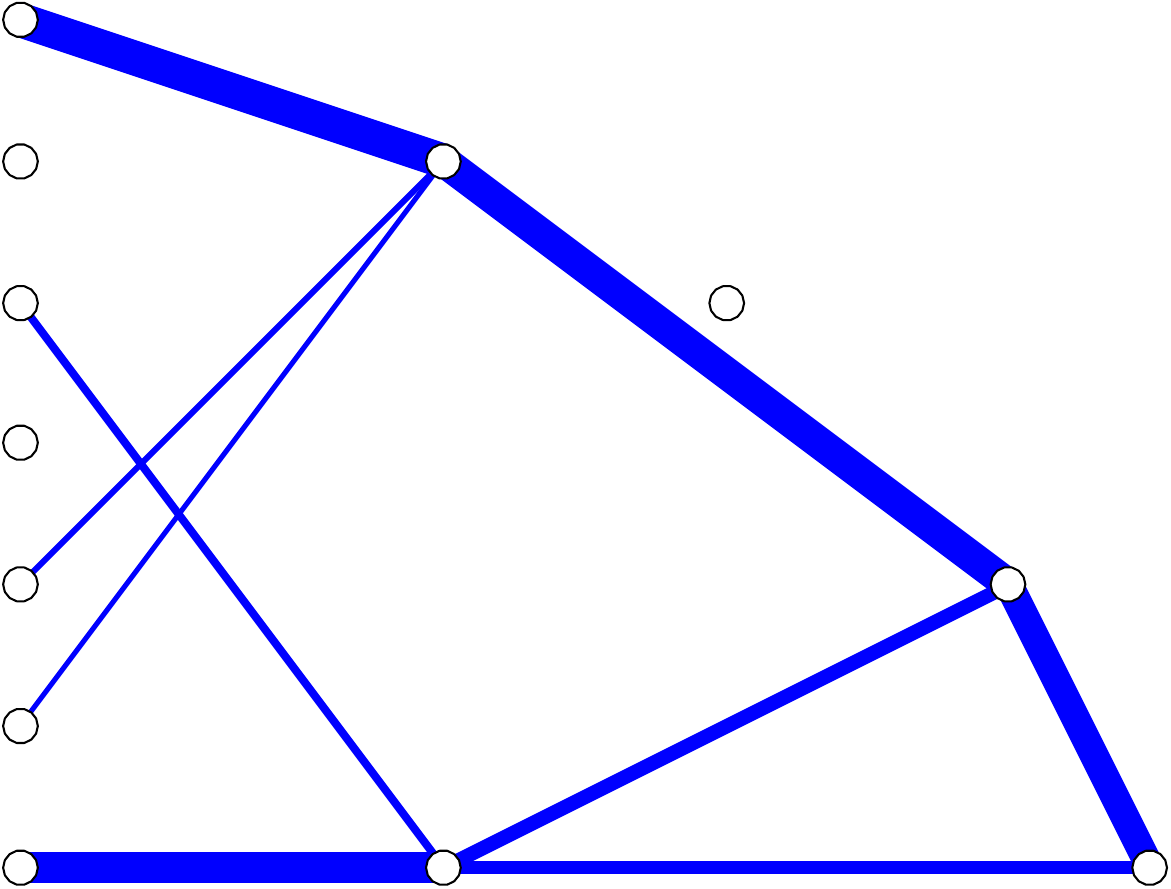}
    \caption{}
    \label{fig:x8_y6_post}
  \end{subfigure}
  \hfill
  \begin{subfigure}[b]{0.45\textwidth}
    \centering
    \includegraphics[scale=0.40]{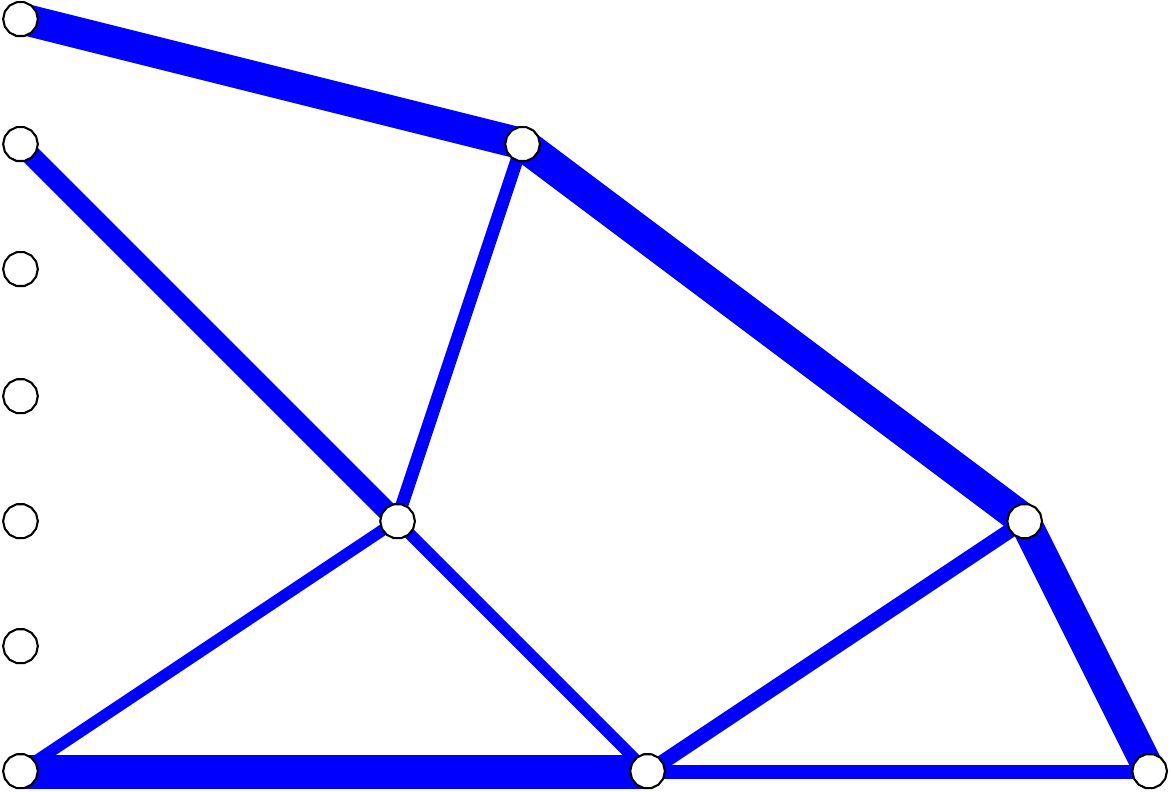}
    \caption{}
    \label{fig:x9_y6_post}
  \end{subfigure}
  \caption{Example (II). 
  The solutions obtained by the ADMM for the compliance minimization 
  with the cardinality constraint ($n=5$). 
  \subref{fig:x8_y2_post} $(N_{X},N_{Y})=(8,2)$; 
  \subref{fig:x9_y2_post} $(N_{X},N_{Y})=(9,2)$; 
  \subref{fig:x8_y4_post} $(N_{X},N_{Y})=(8,4)$; 
  \subref{fig:x9_y4_post} $(N_{X},N_{Y})=(9,4)$; 
  \subref{fig:x8_y6_post} $(N_{X},N_{Y})=(8,6)$; and 
  \subref{fig:x9_y6_post} $(N_{X},N_{Y})=(9,6)$. 
  }
  \label{fig:x8_x9_post}
\end{figure}

\begin{figure}[tp]
  \centering
  \begin{subfigure}[b]{0.45\textwidth}
    \centering
    \includegraphics[scale=0.40]{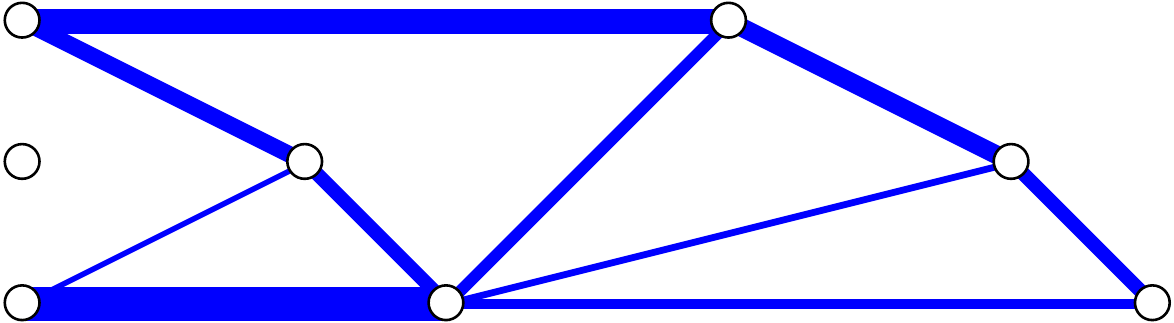}
    \caption{}
    \label{fig:x8_y2_misocp}
  \end{subfigure}
  \hfill
  \begin{subfigure}[b]{0.45\textwidth}
    \centering
    \includegraphics[scale=0.40]{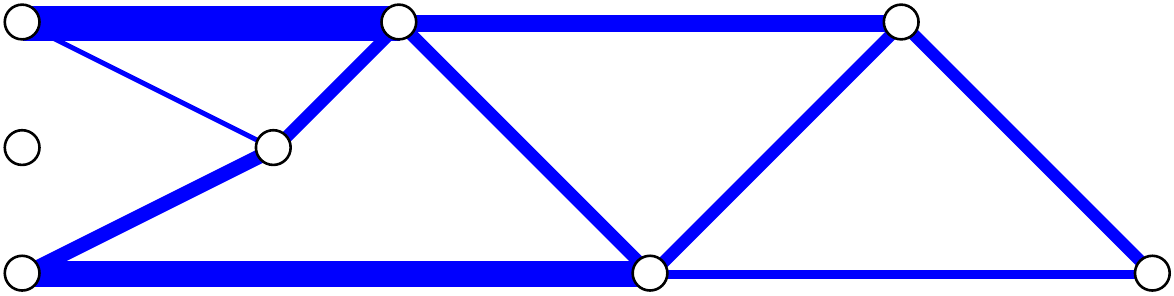}
    \caption{}
    \label{fig:x9_y2_misocp}
  \end{subfigure}
  \par\medskip
  \begin{subfigure}[b]{0.45\textwidth}
    \centering
    \includegraphics[scale=0.40]{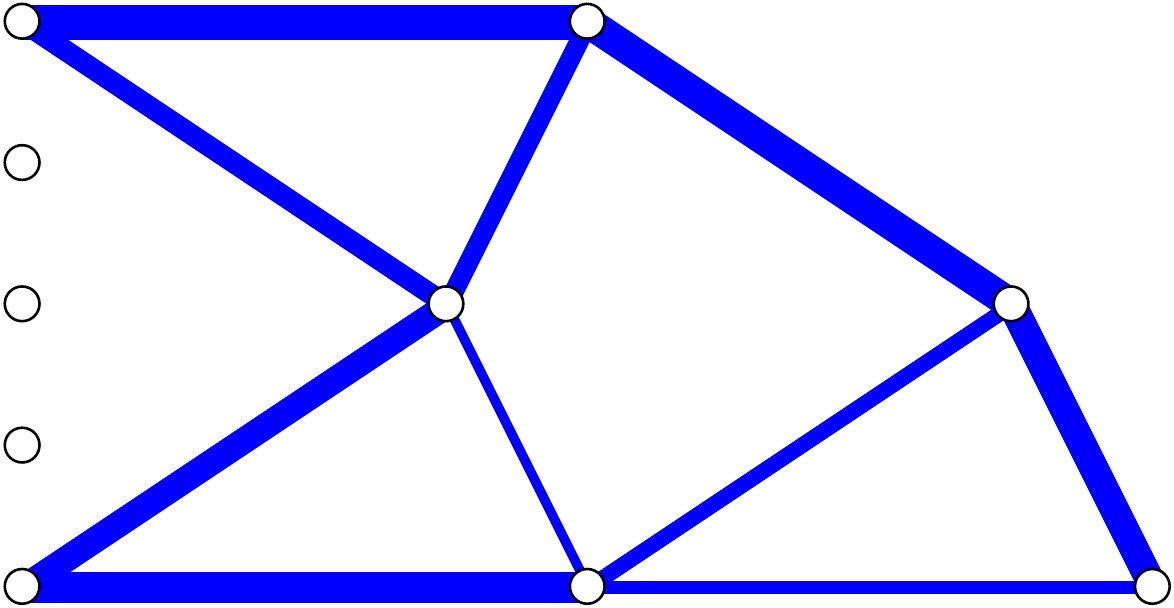}
    \caption{}
    \label{fig:x8_y4_misocp}
  \end{subfigure}
  \hfill
  \begin{subfigure}[b]{0.45\textwidth}
    \centering
    \includegraphics[scale=0.40]{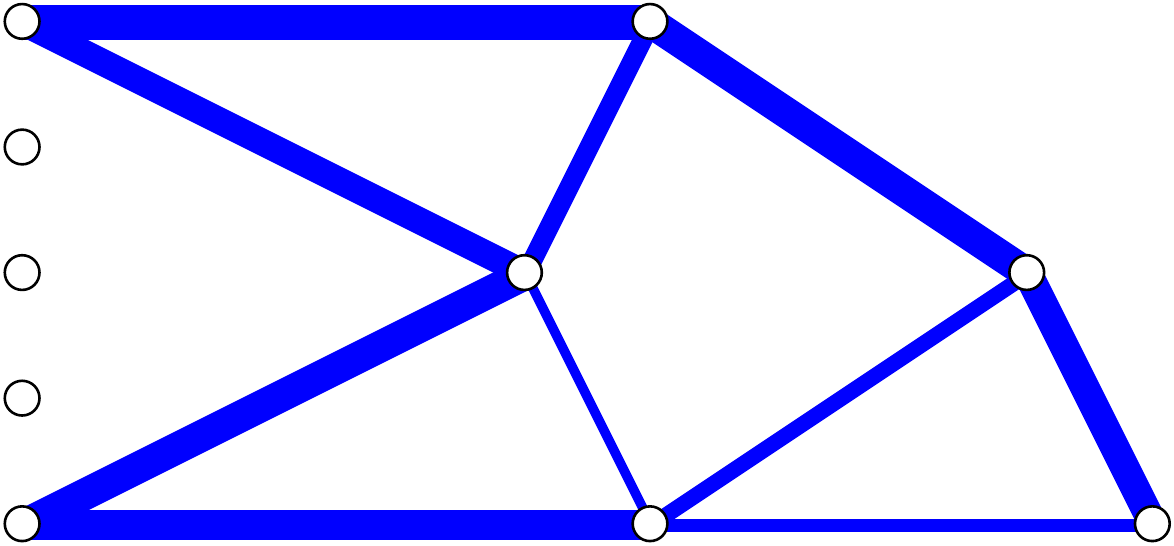}
    \caption{}
    \label{fig:x9_y4_misocp}
  \end{subfigure}
  \par\medskip
  \begin{subfigure}[b]{0.45\textwidth}
    \centering
    \includegraphics[scale=0.40]{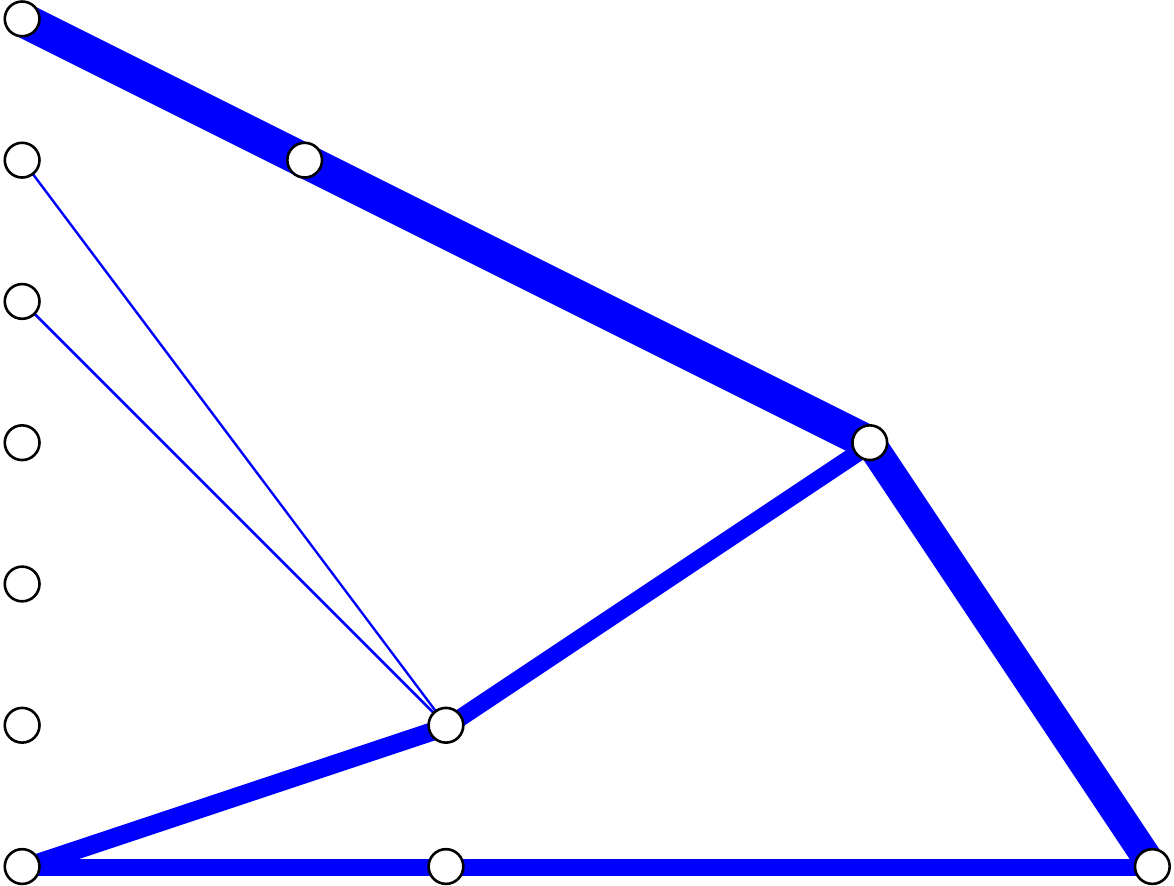}
    \caption{}
    \label{fig:x8_y6_misocp}
  \end{subfigure}
  \hfill
  \begin{subfigure}[b]{0.45\textwidth}
    \centering
    \includegraphics[scale=0.40]{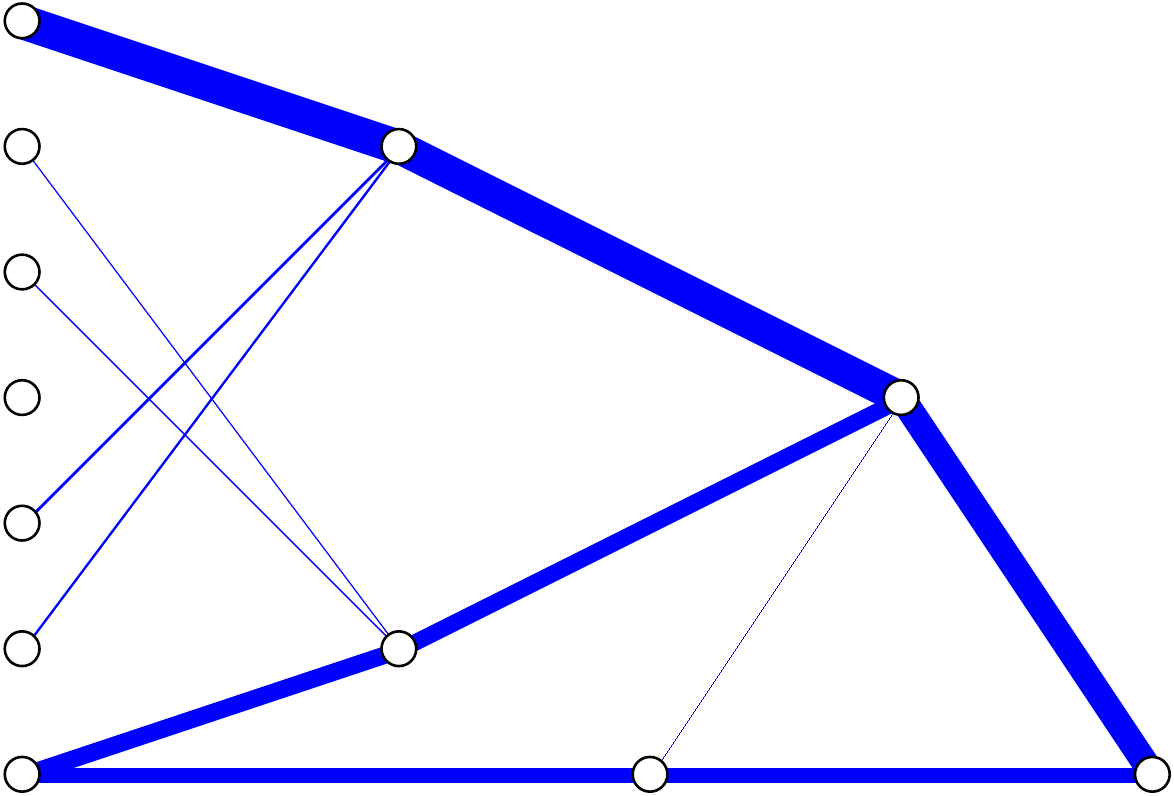}
    \caption{}
    \label{fig:x9_y6_misocp}
  \end{subfigure}
  \caption{Example (II). 
  The optimal solutions obtained by MISOCP for the compliance 
  minimization with the cardinality constraint ($n=5$). 
  \subref{fig:x8_y2_misocp} $(N_{X},N_{Y})=(8,2)$; 
  \subref{fig:x9_y2_misocp} $(N_{X},N_{Y})=(9,2)$; 
  \subref{fig:x8_y4_misocp} $(N_{X},N_{Y})=(8,4)$; 
  \subref{fig:x9_y4_misocp} $(N_{X},N_{Y})=(9,4)$; 
  \subref{fig:x8_y6_misocp} $(N_{X},N_{Y})=(8,6)$; and 
  \subref{fig:x9_y6_misocp} $(N_{X},N_{Y})=(9,6)$. 
  }
  \label{fig:x8_x9_MISOCP}
\end{figure}

In this section, we set the upper bound for the existing free nodes to $n=4$. 
As for problem instances, consider $(N_{X},N_{Y})=(5,2)$, $(5,3)$, and 
$(5,4)$ in \reffig{fig:gs5x2}. 
\reffig{fig:x5_nominal} shows the optimal solutions\footnote{
To obtain these solutions, we used the ground structures without 
overlapping members. } of the conventional 
compliance minimization without the constraint on the number of nodes, 
i.e., problem \eqref{P.compliance.1}, 
where the width of each member is proportional 
to its cross-sectional area. 
\reftab{tab:ex.I.data} reports the optimal values, denoted $\hat{w}$. 
It also lists the number of members ($m$) and the number of degrees of 
freedom of the nodal displacements ($d$). 
As mentioned in Section~\ref{sec:overlap}, the intermediate nodes on a 
chain in \reffig{fig:x5_nominal} can be removed without changing the 
objective value. 
After this hinge cancellation procedure, the numbers of free nodes in 
Figures~\subref{fig:x5_y2_nominal}, \subref{fig:x5_y3_nominal}, and 
\subref{fig:x5_y4_nominal} become $9$, $7$, and $5$, respectively, as 
listed in \reftab{tab:ex.I.data}. 

\reffig{fig:x5_post} shows the solutions obtained by the proposed ADMM 
for the problem with the limited number of free nodes. 
In \reffig{fig:x5_y3_post} we see that the number of free nodes is 
$3$ $(<n)$. 
It should be clear that $\bi{x}^{0}$ used  to generate initial point 
(A) for the ADMM is in general different from the one in \reffig{fig:x5_nominal}, 
because $\bi{x}^{0}$ is computed from the ground structure involving the 
overlapping members. 
Indeed, $\bi{x}^{0}$ for $(N_{X},N_{Y})=(5,2)$, $(5,3)$, and $(5,4)$ 
have $12$, $15$, and $13$ free nodes, respectively. 
Thus, the number of nodes is decreased successfully by the proposed method. 
It is observed that the solutions in 
\reffig{fig:x5_y2_nominal} and \reffig{fig:x5_y3_nominal} have too many 
members from a practical point of view. 
In contrast, we can see in \reffig{fig:x5_y2_post} and 
\reffig{fig:x5_y3_post} that the number of members is decreased as a 
result of optimization with the limitation of the number of nodes. 
The computational results of the ADMM are listed in 
\reftab{tab:ex.I.result}, where $w^{*}$ is the objective value of the 
obtained solution, ``{\#}iter.''\ is the number of iterations required 
before convergence, and ``time'' is the computational time. 
As mentioned before, we examine two different values, denoted (A) and (B), 
for $\bi{z}^{0}$ and $\bi{v}^{0}$. 
The one which yields the better objective value is indicated by ``$*$.'' 
It is observed in \reftab{tab:ex.I.result} that, for every instance, the 
objective value of the solution obtained by the ADMM approach 
is identical to the 
optimal value of the problem without the limitation of the number of 
nodes (i.e., problem \eqref{P.compliance.1}). 
Since problem \eqref{P.compliance.1} can be regarded as a relaxation 
problem, the solutions obtained by the proposed ADMM are globally 
optimal. 
This also illustrates that, in general, the compliance minimization of a 
truss has more than one optimal solution, and the optimal solutions may 
have different numbers of nodes. 

For comparison, we also solved MISOCP \eqref{P.mixed-integer.SOCP.1} 
with a global optimization approach. 
\reftab{tab:ex.I.result} lists the obtained results,\footnote{
It should be clear that no initial point was assigned 
for the MISOCP approach, although in \reftab{tab:ex.I.result}, for 
convenience of presentation, the results of MISOCP are placed in the 
rows concerning the results of the ADMM with initial point (A).} 
 where $\bar{w}$ is the objective value.
The solutions obtained by the MISOCP solver are identical to the ones 
obtained by the ADMM approach.

\begin{table}[bp]
  \centering
  \caption{Computational results of example (II).}
  \label{tab:ex.II.result}
  \begin{tabular}{lrrrrrrrr}
    \toprule
    & \multicolumn{6}{c}{ADMM} & \multicolumn{2}{c}{MISOCP} \\
    \cmidrule(lr){2-7}
    \cmidrule(l){8-9}
    $(N_{X},N_{Y})$ & Init.\ sol.\ & $w^{*}$ (J) 
    & $w^{*}/\hat{w}$ & $w^{*}/\bar{w}$ & {\#}iter & Time (s) 
    & $\bar{w}$ (J)  & Time (s) \\
    \midrule
    (8,2) 
    & $*$ (A) & $37515.63$ & $1.087$ & $1.072$ & $11$ & $9.0$  &$35006.93$&$13.42$ \\
    & $*$ (B) & $37515.63$ & $1.087$ & $1.072$ & $13$ & $9.8$  \\
    \midrule
    (9,2) 
    & $*$ (A) & $50000.00$ & $1.108$ & $1.070$ & $20$ & $18.8$  &$46722.20$&$19.92$ \\
    & $*$ (B) & $50000.00$ & $1.108$ & $1.070$ & $18$ & $16.5$ \\
    \midrule
    (8,4) 
    & $*$ (A) & $7031.25$ & $1.013$ & $1.000$ & $11$ & $46.4$  &$7031.25$&$4.30$ \\
    & $*$ (B) & $7031.25$ & $1.013$ & $1.000$ &  $6$ & $24.5$ \\
    \midrule
    (9,4) 
    & (A) & $9167.44$ & --- & --- & $13$ & $77.2$  \\
    & $*$ (B) & $9000.00$ & $1.011$ & $1.000$ & $8$ & $46.7$ &$8999.91$&$12.20$ \\
    \midrule
    (8,6) 
    & $*$ (A) & $3287.84$ & $1.067$ & $1.038$  & $10$ & $80.3$ &$3168.98$&$76.22$ \\
    & (B) & $3440.05$  & --- & --- & $7$ & $54.1$ \\
    \midrule
    (9,6) 
    & (A) & $4353.91$ & --- & --- & $9$ & $96.5$ \\
    & $*$ (B) & $4221.65$ & $1.097$ & $1.060$ & $12$ & $122.1$  &$3983.34$&$112.08$ \\
    \bottomrule
  \end{tabular}
\end{table}


\subsection{Example (II)}
\label{sec:ex.n5}

As for instances with larger sizes, consider 
$(N_{X},N_{Y})=(8,2)$, $(9,2)$, $(8,4)$, $(9,4)$, $(8,6)$, and $(9,6)$. 
In this section, we set the upper bound for the number of free nodes to $n=5$. 

\subsubsection{Results}

\reffig{fig:x8_x9_nominal} collects the optimal solutions without 
limiting the number of nodes. 
The number of free nodes after applying the hinge cancellation is 
reported in \reftab{tab:ex.I.data}. 
\reffig{fig:x8_x9_post} shows the solutions obtained by the ADMM approach. 
The number of free nodes in \reffig{fig:x8_y6_post} is 
$4$ ($<n$). 
Two nodes can be removed from the solution in 
\reffig{fig:x8_y2_post}, which results in a truss design with three free 
nodes. 
It is observed in \reffig{fig:x8_x9_nominal} and \reffig{fig:x8_x9_post} 
that the limitation of the number of nodes often yields a solution with 
a fewer members. 
Also, too thin members are observed in \reffig{fig:x8_x9_nominal}, while 
such thin members do not appear in \reffig{fig:x8_x9_post}. 
These two features of the solutions in \reffig{fig:x8_x9_post} are 
considered practically preferable. 
The initial design $\bi{x}^{0}$ used for generating initial points (A) 
for the ADMM to solve $(N_{X},N_{Y})=(8,2)$, $(9,2),\dots,(9,6)$ 
have $20$, $18$, $19$, $20$, $24$, and $26$ free nodes, respectively. 

It is observed in \reftab{tab:ex.II.result} that the ADMM terminates 
after at most $20$ iterations. 
Increase of the objective value from the optimal value of the problem 
without the cardinality constraint is quite small, i.e., increase by at 
most about 10\%. 
Particularly, for $(N_{X},N_{Y})=(8,4)$ and $(9,4)$ we have only about 
1\% increase. 
Thus, it is often that the number of nodes can be reduced at the expense 
of only small increase of the compliance. 

The computational results of the MISOCP approach are listed in 
\reftab{tab:ex.II.result}. 
\reffig{fig:x8_x9_MISOCP} collects the obtained solutions. 
For $(N_{X},N_{Y})=(8,4)$, the solution obtained by MISOCP is identical 
to the one obtained by the ADMM; i.e., the ADMM found a global optimal solution. 
For $(N_{X},N_{Y})=(9,4)$, it is observed in \reftab{tab:ex.II.result} 
that the objective values obtained by the two methods are almost same, 
but the two solutions are slightly different as seen in 
\reffig{fig:x9_y4_post} and \reffig{fig:x9_y4_misocp}. 
The largest value of $w^{*}/\bar{w}$ is $1.072$ in the case of 
$(N_{X},N_{Y})=(8,2)$. 
It is also worth noting that, for $(N_{X},N_{Y})=(8,6)$ and $(9,6)$, 
although the global optimal solutions in \reffig{fig:x8_y6_misocp} and 
\reffig{fig:x9_y6_misocp} involve very thin members, the solutions 
obtained by the ADMM shown in \reffig{fig:x8_y6_post} and 
\reffig{fig:x9_y6_post} do not have such a thin member. 

It is observed from \reftab{tab:ex.II.result} 
that the proposed ADMM often converges more 
quickly than the MISOCP solver; exceptions are 
$(N_{X},N_{Y})=(8,4)$ and $(9,4)$. 
The computational time required by the MISOCP solver varies drastically 
depending on problem instances. 
In contrast, the number of iterations required by the ADMM is almost 
independent of problem instances. 
Since the computational time required for solving the SOCP subproblem of 
the ADMM depends on the problem size, it is possible to roughly estimate 
the total computational cost of the ADMM from the problem size. 
This might be considered one of advantages of ADMM over MISOCP. 

As mentioned in Section~\ref{sec:overlap}, the proposed method does not 
incorporate the constraint prohibiting overlapping members. 
Nevertheless, for the problems with limitation of the number of nodes, 
all the solutions obtained in Section~\ref{sec:ex} do not involve 
overlapping members. 

\subsubsection{MISOCP with slenderness constraints}

The constraints preventing the presence of very thin members observed in 
\reffig{fig:x8_y6_misocp} and \reffig{fig:x9_y6_misocp} can be handled 
within the framework of mixed-integer programming (MIP) \cite{KY17}. 
Recall problem \eqref{P.mixed-integer.SOCP.2} in 
Section~\ref{sec:overlap}, where $t_{i}$ is a binary variable indicating 
whether member $i$ exists or vanishes. 
Let $x_{\rr{min}} > 0$ denote the specified lower bound for the member 
cross-sectional area. 
The constraint avoiding existence of too thin members can be formulated as 
\begin{align}
  x_{\rr{min}} \bi{t} \le \bi{x} \le M \bi{t} . 
  \label{eq:slender.constraint.0}
\end{align}

\begin{figure}[tp]
  \centering
  \begin{subfigure}[b]{0.45\textwidth}
    \centering
    \includegraphics[scale=0.40]{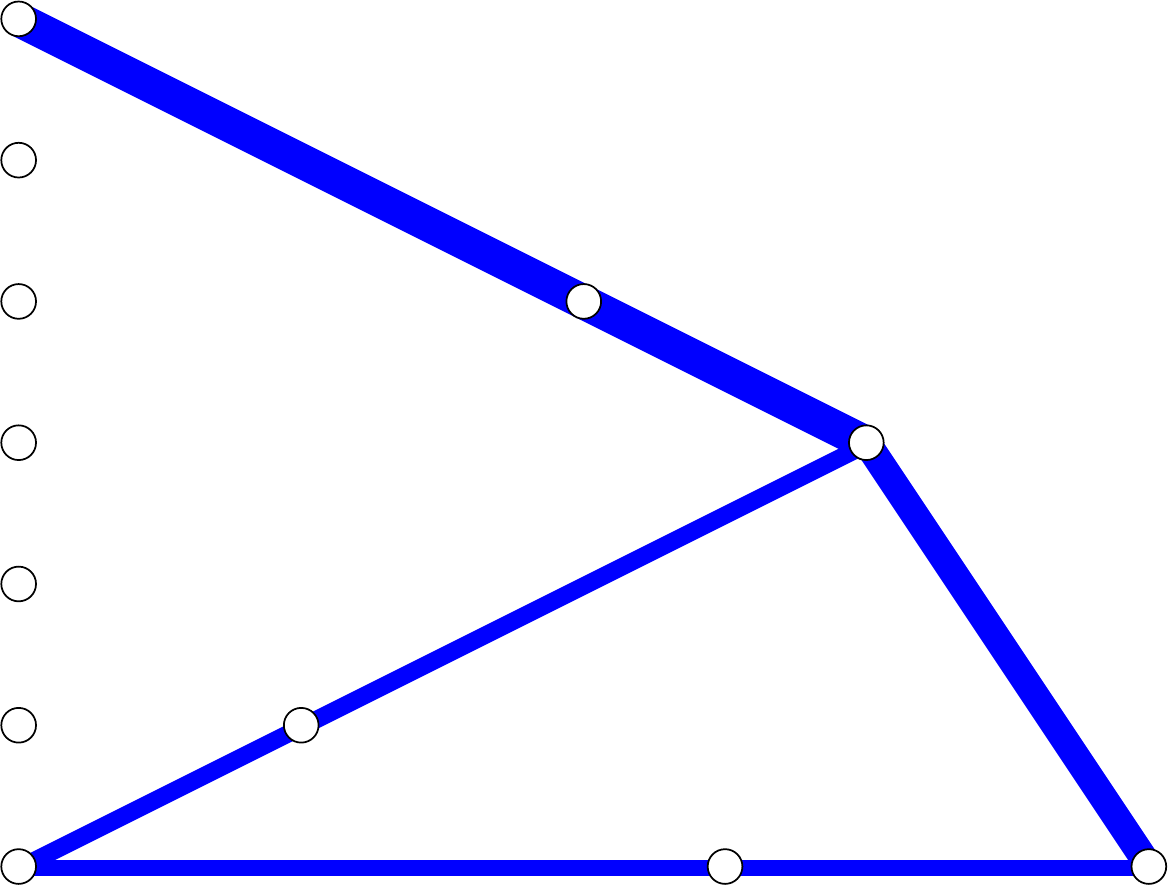}
    \caption{}
    \label{fig:8x6under5_misocp_xmin=200.0_vol_0.0096_nod_limit=5}
  \end{subfigure}
  \hfill
  \begin{subfigure}[b]{0.45\textwidth}
    \centering
    \includegraphics[scale=0.40]{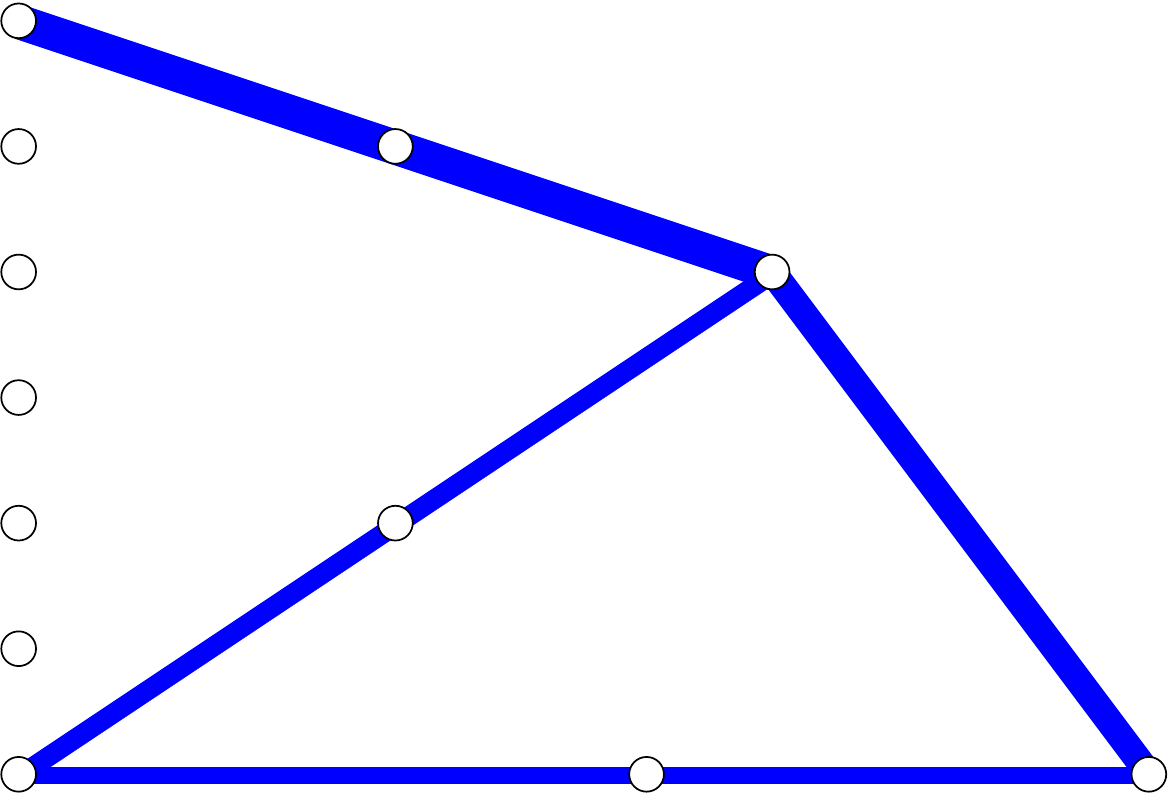}
    \caption{}
    \label{fig:9x6under5_misocp_xmin=200.0_vol_0.0108_nod_limit=5}
  \end{subfigure}
  \caption{
  The optimal solutions of example (II) 
  with the constraints avoiding the presence of thin members. 
  \subref{fig:8x6under5_misocp_xmin=200.0_vol_0.0096_nod_limit=5} $(N_{X},N_{Y})=(8,6)$; and 
  \subref{fig:9x6under5_misocp_xmin=200.0_vol_0.0108_nod_limit=5} $(N_{X},N_{Y})=(9,6)$. 
  }
  \label{fig:misocp_xmin=200.0_vol_0.0096_nod_limit=5}
\end{figure}

In problem \eqref{P.mixed-integer.SOCP.2}, we replace constraint 
\eqref{P.mixed-integer.SOCP.2.5} with \eqref{eq:slender.constraint.0}. 
The constraint avoiding the presence of overlapping members, 
\eqref{P.mixed-integer.SOCP.2.9}, is not considered. 
We solve this MISOCP for the instances $(N_{X},N_{Y})=(8,6)$ and $(9,6)$ 
with $x_{\rr{min}}=200\,\mathrm{mm^{2}}$. 
The obtained solutions are shown in 
\reffig{fig:misocp_xmin=200.0_vol_0.0096_nod_limit=5}. 
Both solutions have parallel consecutive members that are connected by 
nodes supported only in the direction of those members. 
The intermediate nodes can be removed without changing the optimal value. 
Hence, the number of free nodes of these solutions is essentially two. 
The objective value of the solution for $(N_{X},N_{Y})=(8,6)$ is 
$3168.97\,\mathrm{J}$, which is slightly less than that for the case without the 
slenderness constraints in \reftab{tab:ex.II.result}. 
This is due to the computational error in computing the objective value 
with the finite element method. 
In the solution shown in \reffig{fig:x8_y6_post} (i.e., the solution 
obtained by the ADMM without the slenderness constraints), the 
cross-sectional area of the thinest member is $75.4\,\mathrm{mm^{2}}$. 
Hence, this solution is not globally optimal under the slenderness 
constraint. 
The computational time required by MOSEK is $365.0\,\mathrm{s}$. 
The objective value of the solution for $(N_{X},N_{Y})=(9,6)$ is 
$4028.94\mathrm{J}$. 
This is larger than that for the case without the 
slenderness constraints as expected, and is less than that of the 
solution obtained by the ADMM. 
Since the cross-sectional area of the thinnest member of the 
solution shown in \reffig{fig:x9_y6_post} is $146.7\,\mathrm{mm^{2}}$, 
Hence, the solution shown in \reffig{fig:x9_y6_post} is not globally 
optimal under the slenderness constraints. 
The computational time required by MOSEK to find the solution in 
\reffig{fig:9x6under5_misocp_xmin=200.0_vol_0.0108_nod_limit=5} was 
$1676.0\,\mathrm{s}$, which is much larger than the computational time 
of the ADMM. 


\subsection{Example (III)}
\label{sec:ex.7m}

\begin{figure}[tp]
  \centering
  \begin{subfigure}[b]{0.45\textwidth}
    \centering
    \includegraphics[scale=0.40]{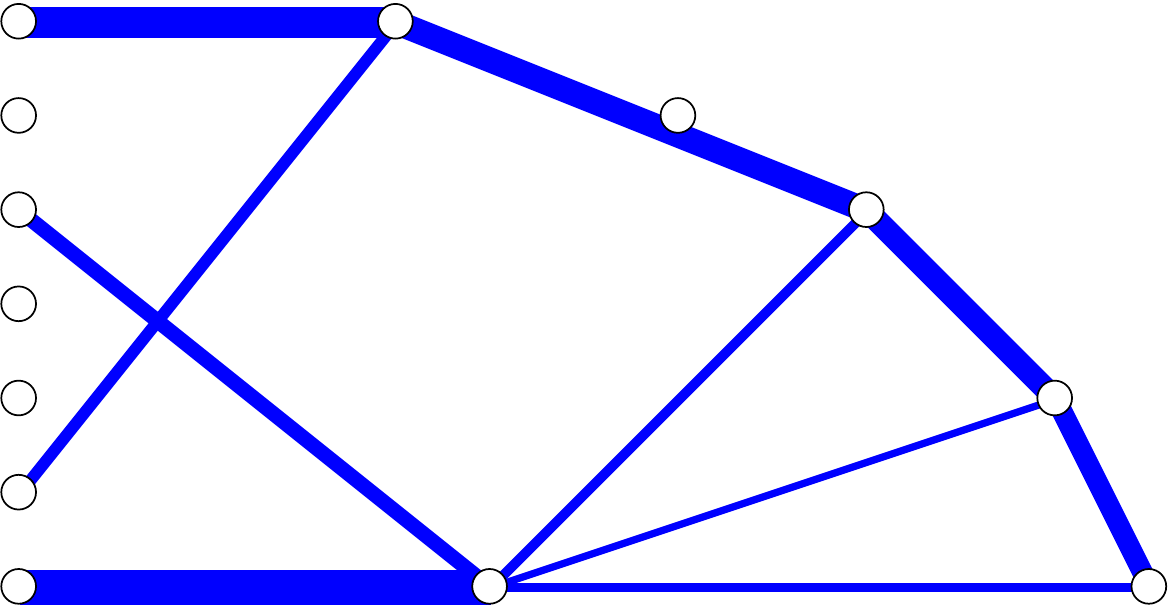}
    \caption{}
    \label{fig:x12_y6_post}
  \end{subfigure}
  \hfill
  \begin{subfigure}[b]{0.45\textwidth}
    \centering
    \includegraphics[scale=0.40]{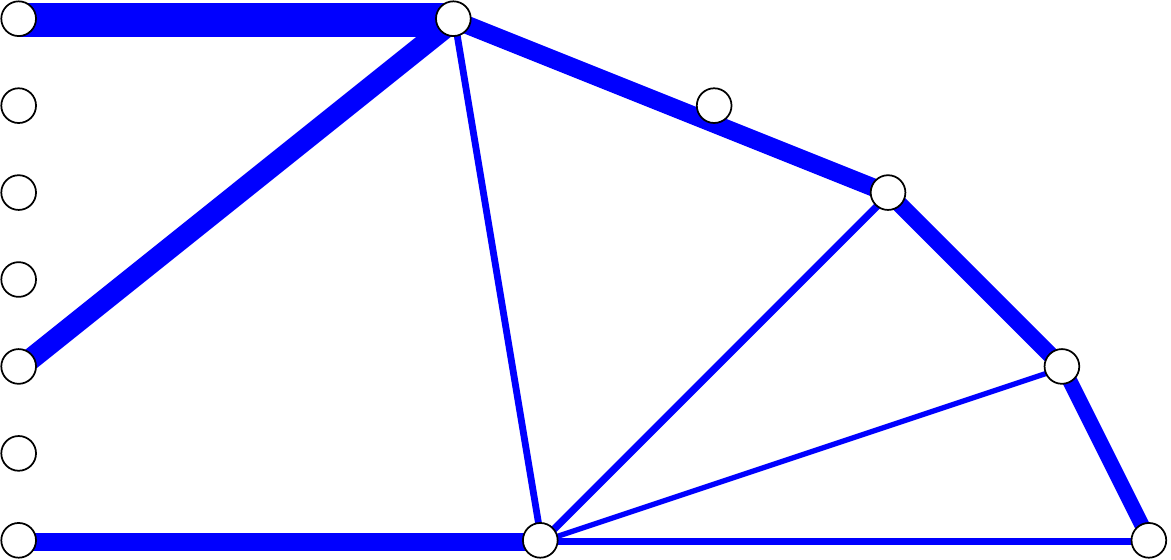}
    \caption{}
    \label{fig:x13_y6_post}
  \end{subfigure}
  \par\medskip
  \begin{subfigure}[b]{0.45\textwidth}
    \centering
    \includegraphics[scale=0.40]{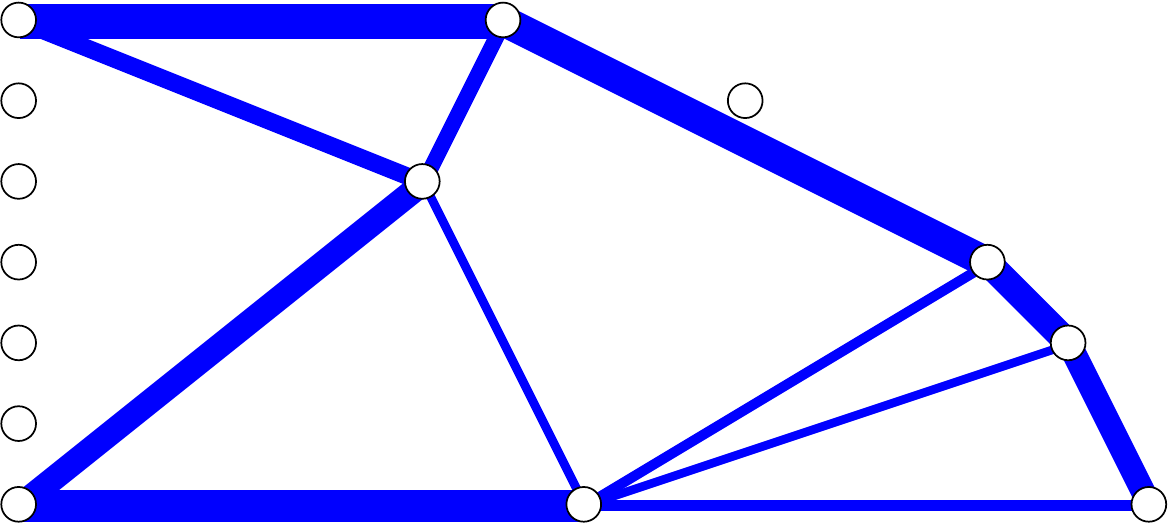}
    \caption{}
    \label{fig:x14_y6_post}
  \end{subfigure}
  \caption{Example (III). 
  The solutions obtained by ADMM for the compliance minimization 
  with the cardinality constraint. 
  \subref{fig:x12_y6_post} $(N_{X},N_{Y},n)=(12,6,6)$; 
  \subref{fig:x13_y6_post} $(N_{X},N_{Y},n)=(13,6,6)$; and 
  \subref{fig:x14_y6_post} $(N_{X},N_{Y},n)=(14,6,7)$, 
  }
  \label{fig:x12_x_13_x14_post}
\end{figure}

\begin{figure}[tp]
  \centering
  \begin{subfigure}[b]{0.45\textwidth}
    \centering
    \includegraphics[scale=0.40]{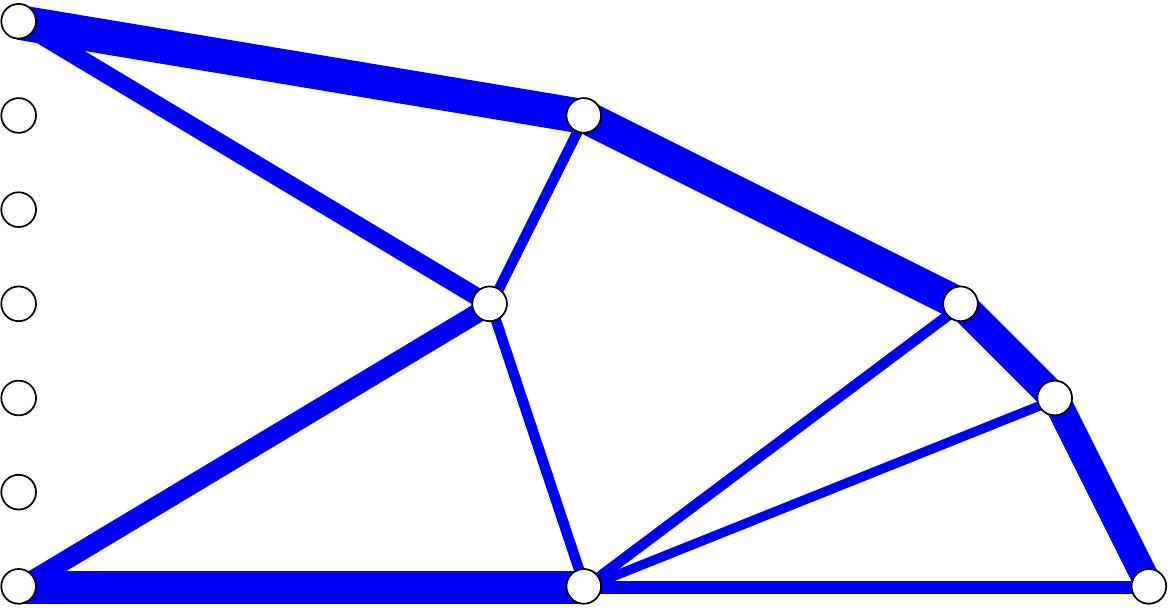}
    \caption{}
    \label{fig:x12_y6_mip}
  \end{subfigure}
  \hfill
  \begin{subfigure}[b]{0.45\textwidth}
    \centering
    \includegraphics[scale=0.40]{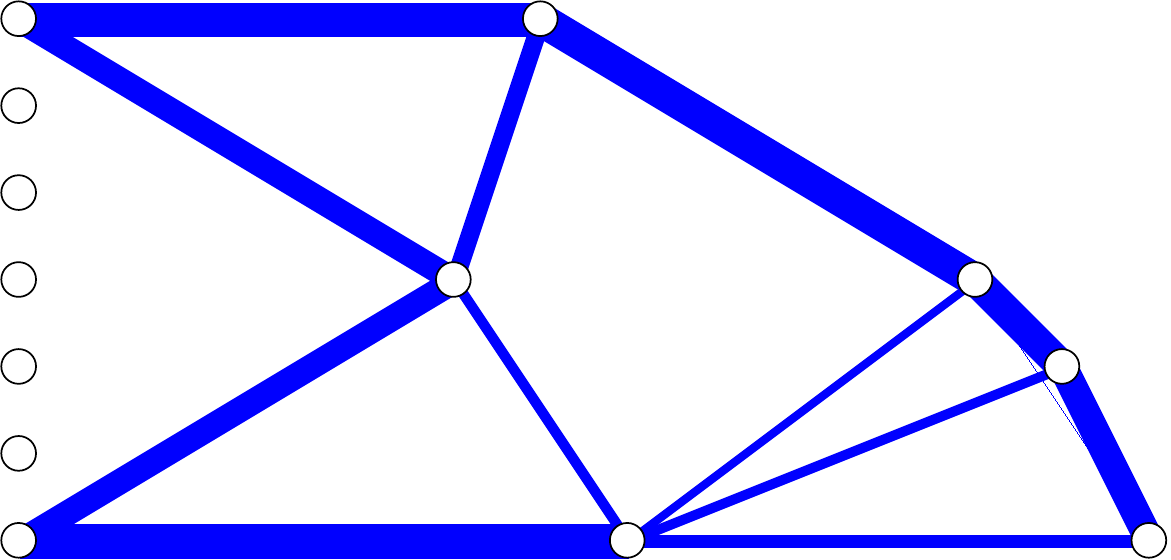}
    \caption{}
    \label{fig:x13_y6_mip}
  \end{subfigure}
  \par\medskip
  \begin{subfigure}[b]{0.45\textwidth}
    \centering
    \includegraphics[scale=0.40]{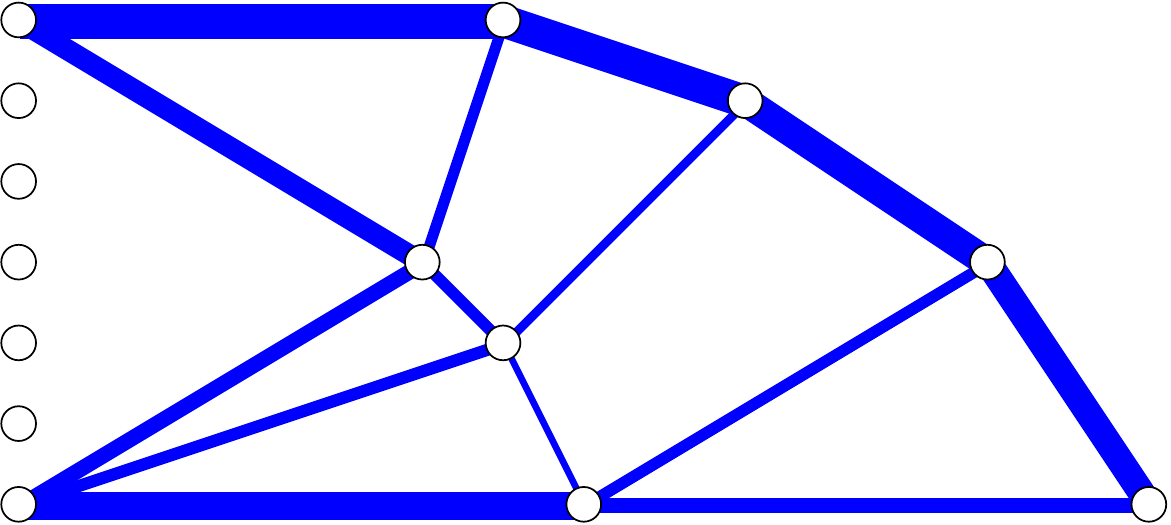}
    \caption{}
    \label{fig:x14_y6_mip}
  \end{subfigure}
  \caption{Example (III). 
  The optimal solutions obtained by MISOCP for the compliance minimization 
  with the cardinality constraint. 
  \subref{fig:x12_y6_mip} $(N_{X},N_{Y},n)=(12,6,6)$; 
  \subref{fig:x13_y6_mip} $(N_{X},N_{Y},n)=(13,6,6)$; and 
  \subref{fig:x14_y6_mip} $(N_{X},N_{Y},n)=(14,6,7)$. 
  }
  \label{fig:x12_x_13_x14_mip}
\end{figure}

Consider problem instances 
$(N_{X},N_{Y},n)=(12,6,6)$, $(13,6,6)$, and $(14,6,7)$. 
The maximum member length in a ground structure is set to $7\,\mathrm{m}$. 

\reffig{fig:x12_x_13_x14_post} shows the solutions obtained by the 
proposed ADMM approach. 
The ADMM terminates with a solution having $n$ free nodes. 
One of these nodes vanishes in the post-processing. 
The objective value as well as the computational cost is reported in 
\reftab{tab:ex.III.result}. 

\reffig{fig:x12_x_13_x14_mip} collects the optimal solutions found by 
the MISOCP approach. 
These solutions use exactly $n$ free nodes. 
It is observed from \reftab{tab:ex.III.result} that the objective value 
obtained by the ADMM for the largest instance, 
$(N_{X},N_{Y},n)=(14,6,7)$, is very close to the optimal value. 
In contrast, for $(N_{X},N_{Y},n)=(13,6,6)$ the objective value obtained 
by the ADMM is more than 20\% larger than the optimal value. 
However, the computational cost of the ADMM is much less than the MISOCP 
approach (which requires more than four hours). 
Thus, the quality of the solution obtained by the ADMM approach can 
possibly be very good, although in general it depends on problem 
instances. 
As the problem size increases, the computational cost of the ADMM 
approach becomes much smaller compared with the MISOCP approach. 

\begin{table}[tbp]
  \centering
  \caption{Computational results of example (III). }
  \label{tab:ex.III.result}
  \begin{tabular}{lrrrrrrr}
    \toprule
    & \multicolumn{5}{c}{ADMM} & \multicolumn{2}{c}{MISOCP} \\
    \cmidrule(lr){2-6}
    \cmidrule(l){7-8}
    $(N_{X},N_{Y},n)$ & Init.\ sol.\ & $w^{*}$ (J) 
    & $w^{*}/\bar{w}$ & {\#}iter & Time (s) 
    & $\bar{w}$ (J)  & Time (s) \\
    \midrule
    $(12,6,6)$ 
    & $*$ (A) & $7817.72$ & $1.115$ & 10 & 531.7 &$7012.95$ & 3492.34 \\
    &     (B) & $9527.47$ &  & 10 & 516.7 & & \\
    \midrule
    $(13,6,6)$ 
    & $*$ (A) & $10142.88$ & $1.228$ & 10 & 642.6 &$8258.85$ & 14675.17 \\
    &     (B) & $13631.69$ &  & 11 & 793.7 & & \\
    \midrule
    $(14,6,7)$ 
    & $*$ (A) & $9720.09$ & $1.018$ & 9 & 724.2 &$9550.97$ &14657.54  \\
    &     (B) & $14928.49$ &  & 14 & 1148.6 & &  \\
    \bottomrule
  \end{tabular}
\end{table}

\subsection{On heuristic for stopping ADMM}
\label{sec:ex.stopping}

\begin{table}[bp]
  \centering
  \caption{The computational results when the ADMM is run until the 
  convergence to evaluate effectiveness of the heuristic stopping criterion. }
  \label{tab:ex.post.heuristic}
  \begin{tabular}{llrrr}
    \toprule
    $(N_{X},N_{Y},n)$ & Init.\ sol.\ 
    & $\tilde{K}$ & $\tilde{K}_{\rr{p}}$ & $K^{*}$ \\
    \midrule
    $(8,2,5)$ & (A) & $35$ & $18$ & $11$ \\
    $(8,2,5)$ & (B) & $35$ & $18$ & $13$ \\
    $(9,2,5)$ & (A) & $100^{\dag}$ & $77^{\dag}$ & $20$ \\
    $(9,2,5)$ & (B) & $81$ & $58$ & $18$ \\
    $(8,4,5)$ & (A) & $27$ & $15$ & $11$ \\
    $(8,4,5)$ & (B) & $21$ & $15$ & $6$ \\
    $(9,4,5)$ & (A) & $29$ & $15$ & $13$ \\
    $(9,4,5)$ & (B) & $25$ & $15$ & $8$ \\
    $(8,6,5)$ & (A) & $26$ & $12$ & $10$ \\
    $(8,6,5)$ & (B) & $21$ & $15$ & $7$ \\
    $(9,6,5)$ & (A) & $24$ & $16$ & $9$ \\
    $(9,6,5)$ & (B) & $36$ & $25$ & $12$ \\
    \bottomrule
  \end{tabular}
\end{table}

As mentioned in section~\ref{sec:ex.implementation}, we use a heuristic 
criterion for stopping the ADMM. 
Namely, we stop the ADMM when the cardinality constraint is satisfied 
with $\epsilon$ tolerance. 
Then, as for a post-processing, we solve the compliance minimization 
problem, \eqref{P.compliance.1}, with specifying the set of vanishing nodes. 
This section presents some empirical justification for this procedure. 
Namely, it is illustrated through numerical experiments that with this 
heuristic procedure the number of subproblems to be solved is 
drastically reduced, without missing out better solutions in the sense 
of the objective valued. 
We use the problem instances in Section~\ref{sec:ex.n5}. 

We performed the following experiment. 
The ADMM is run until it terminates with a small tolerance, namely, 
$\| \bi{x}^{k+1} - \bi{x}^{k} \| \le 10^{-1}$ (in $\mathrm{mm}^{2}$) 
is satisfied. 
This requires much more iterations compared with the procedure described above. 
In the iteration history, we select every iterate that 
satisfies the cardinality constraint approximately, i.e., that satisfies 
$l-|J_{0}^{k}| \le n$. 
For every selected iterate, we run the post-processing, i.e., we solve 
problem \eqref{P.compliance.1} with specifying the set of vanishing nodes. 

The computational results are listed in \reftab{tab:ex.post.heuristic}, 
where $\tilde{K}$ and $\tilde{K}_{\rr{p}}$ are the number of iterations 
required  before convergence and the number of iterates that 
approximately satisfy the cardinality constraint, respectively. 
For all the $\tilde{K}_{\rr{p}}$ iterates, the post-processing yields 
the same solution as the one reported in \reftab{tab:ex.II.result}. 
In the case $(N_{X},N_{Y},n)=(9,2,5)$ with initial point (A), the ADMM 
does not converge after $100$ iterates. 
In this iteration history, there exist $77$ iterates satisfying the 
cardinality constraint approximately. 
From all of them, the post-processing generates the solution reported in 
\reftab{tab:ex.II.result}. 
In a nutshell, the set of vanishing nodes does not change, even if the 
ADMM iterations are continued after the iterate at which our heuristic 
stopping criterion is satisfied. 

For ease of comparison, the number of iterations reported in 
\reftab{tab:ex.II.result} is listed again as $K^{*}$ in 
\reftab{tab:ex.post.heuristic}. 
Namely, our stopping criterion reduces the number of iterations from 
$\tilde{K}$ to $K^{*}$, without changing the final output. 
It is worth noting that the solutions found at between the $(K^{*}+1)$ 
iterate and $\tilde{K}$ iterate do not necessarily satisfy the 
cardinality constraint with $\epsilon$ tolerance. 


\subsection{On choice of initial points}
\label{sec:ex.initial}

\begin{table}[bp]
  \centering
  \caption{Computational results of ADMM from randomly generated initial points. }
  \label{tab:random.experiment}
  \begin{tabular}{llrrrr}
    \toprule
    $(N_{X},N_{Y},n)$ & Init.\ sol.\ 
    & Min.\ ($\mathrm{J}$) & Max.\ ($\mathrm{J}$) & Mean ($\mathrm{J}$) & Var.\ ($\mathrm{J}^{2}$) \\
    \midrule
    $(5,3,4)$ & (C) & $5052.45$ & $5052.45$ & $5052.45$ & $0.00$ \\
    $(5,3,4)$ & (D) & $5007.41$ & $5052.45$ & $5049.74$ & $115.56$ \\
    \midrule
    $(9,4,5)$ & (C) & $9000.00$ & $9143.42$ & $9039.78$ & $3344.71$ \\
    $(9,4,5)$ & (D) & $9000.00$ & $9768.04$ & $9076.47$ & $13916.47$ \\
    \bottomrule
  \end{tabular}
\end{table}

Since we apply ADMM to a nonconvex problem, the obtained solution may 
depend on the choice of initial points. 
In Section~\ref{sec:ex.implementation}, we suggest to use the two 
initial points, (A) and (B), and adopt the better solution as the final 
output. 
In this section, we perform comparison with the results obtained by 
using randomly generated initial points to empirically justify our selection. 
As for two representative instances for which initial points (A) and (B) 
lead to different solution, we use 
$(N_{X},N_{Y},n)=(5,3,4)$ in section~\ref{sec:ex.implementation} and 
$(N_{X},N_{Y},n)=(9,4,5)$ in section~\ref{sec:ex.n4} in the following 
numerical experiments. 

As for randomly generated initial points, we examine two cases: 
\begin{itemize}
  \item Initial point (C): 
        $\bi{x}^{0} := (V / \bi{c}^{\top} \bi{\xi}) \bi{\xi}$, 
        $\bi{z}^{0} := Z \bi{x}^{0}$, and 
        $\bi{v}^{0} := Z \bi{x}^{0} - \bi{z}^{0} = \bi{0}$, where 
        $\bi{\xi}\in\Re^{m}$ is a random vector with the entries drawn from 
        $\UC(0,1)$. 
  \item Initial point (D): 
        $\bi{x}^{0} := (V / \bi{c}^{\top} \bi{\xi}) \bi{\xi}$, 
        $\bi{z}^{0} := Z \bi{x}^{0}$, and 
        $\bi{v}^{0} := \max\{ z_{1}^{0},\dots,z_{l}^{0} \} \bi{\zeta}$, 
        where $\bi{\xi}\in\Re^{m}$ and $\bi{\zeta} \in \Re^{l}$ are 
        random vectors with the entries drawn from $\UC(0,1)$. 
\end{itemize}
We generate 100 samples of each of these initial points, and run our 
ADMM approach from every sample. 
\reftab{tab:random.experiment} reports the minimum value, maximum value, 
mean, and variance of the objective value. 

For $(N_{X},N_{Y},n)=(5,3,4)$ with initial point (C), in all the 
cases the ADMM converges to the same solution. 
This solution is the one obtained by using initial point (B), as shown 
in \reftab{tab:ex.I.result}. 
Therefore, using initial point (A) yielded a better solution (which is 
globally optimal) than using $100$ samples of (C). 
In contrast, when initial point (D) was adopted, the global optimal 
solution is obtained from $7$ sampled initial points, among $100$ trials. 
From the other $93$ samples, the ADMM converges to the solution obtained 
with initial point (B). 
The mean and the variance of the objective value are listed in 
\reftab{tab:random.experiment}. 
In this manner, it is demonstrated that the global optimal solution, 
easily obtained by carrying out our ADMM procedure with initial point 
(A), is rarely obtained from randomly generated initial points. 

For $(N_{X},N_{Y},n)=(9,4,5)$ with initial point (C), the best solution 
is same as the one obtained from initial point (B) in 
\reftab{tab:ex.II.result}. 
This is not globally optimal. 
Among $100$ trials, $67$ sampled initial points yield this solution. 
In contrast, the objective value of the worst solution is larger than 
the one obtained from initial point (A). 
By using initial point (D), the variation of the objective value 
increased, but the global optimal solution was not obtained. 
Thus, the ADMM with randomly generated initial points could not find a 
solution better than the one obtained from initial point (B). 

In short, for these two problem instances, using many randomly generated 
initial points does not yield a better solution. 
Therefore, using initial points (A) and (B) might be considered a 
reasonable selection.

\subsection{Application to robust optimization against uncertainty in external load}
\label{sec:ex.robust}

\begin{table}[bp]
  \centering
  \caption{Computational results of the robust optimization with the 
  cardinality constraint.}
  \label{tab:ex.I.SDP}
  \begin{tabular}{lrrrrrrrr}
    \toprule
    & \multicolumn{5}{c}{ADMM} & \multicolumn{3}{c}{MISDP} \\
    \cmidrule(lr){2-6}
    \cmidrule(l){7-9}
    $(N_{X},N_{Y},n)$ &  Init.\ sol.\ & $w^{*}$ (J) & $w^{*}/\bar{w}$ & {\#}iter & Time (s) 
    & $\bar{w}$ (J) & {\#}iter & Time (s) \\
    \midrule
    $(5,2,4)$ 
    & $*$ (A) & $12156.25$ & $1.000$ & $3$ & $7.8$ & $12156.25$ & $22$ & $11.8$ \\
    & $*$ (B) & $12156.25$ & $1.000$ & $3$ & $6.7$ \\
    \midrule
    $(5,3,4)$ 
    & $*$ (A) & $5044.91$ & $1.000$ & $5$ & $6.6$ & $5044.91$ & $9$ & $7.6$ \\
    &     (B) & $5089.95$ & ---    & $3$ & $4.7$ \\
    \midrule
    $(5,4,4)$ 
    & $*$ (A) & $2840.63$ & $1.000$ & $2$ & $6.7$ & $2840.63$ & $5$ & $10.9$ \\
    & $*$ (B) & $2840.63$ & $1.000$ & $3$ & $7.4$ \\
    \midrule
    $(8,2,5)$ 
    & $*$ (A) & $37605.63$ & $1.071$ & $11$ & $16.5$ & $35096.94$ & $238$ & $293.1$ \\
    & $*$ (B) & $37605.63$ & $1.071$ & $13$ & $18.2$ \\
    \midrule
    $(9,2,5)$ 
    & $*$ (A) & $50101.25$ & $1.070$ & $20$ & $32.0$ & $46823.47$ & $1484$ & $2396.2$ \\
    & $*$ (B) & $50101.25$ & $1.070$ & $14$ & $20.6$ \\
    \midrule
    $(8,4,5)$ 
    & $*$ (A) & $7076.25$ & $1.000$ & $11$ & $67.0$ & $7076.25$ & $37$ & $334.1$ \\
    & $*$ (B) & $7076.25$ & $1.000$ & $6$  & $34.3$ \\
    \midrule
    $(9,4,5)$ 
    &     (A) & $9218.06$ & --- & $13$ & $110.6$ \\
    & $*$ (B) & $9050.63$ & $1.000$ & $6$  & $47.0$ & $9050.63$ & $69$ & $1026.2$ \\
    \midrule
    $(8,6,5)$ 
    & $*$ (A) & $3317.84$  & $1.037$ & $10$ & $118.8$ & $3198.98$ & $573$ & $21992.7$ \\
    &     (B) & $26377.74$ & --- & $19$ & $358.5$ \\
    \midrule
    $(9,6,5)$ 
    &     (A) & $4387.66$ & --- & $9$  & $148.6$ \\
    & $*$ (B) & $4255.40$ & ${\ddag}$ & $12$ & $181.9$ 
    & $4387.66^{\ddag}$ & $300^{\ddag}$ & $16468.0^{\ddag}$ \\
    \bottomrule
  \end{tabular}
\end{table}

The ADMM approach presented in this paper can be easily extended to the 
case in which the external load possesses uncertainty. 
The set of nodes at which the external forces can possibly be applied is 
supposed to be specified. 
Then we consider the robust optimization against the uncertainty, under 
the upper bound constraint on the number of nodes. 
In this section we examine efficiency of the ADMM applied to this 
problem, as an example of optimization problems that are not handled with 
current mainstream MIP solvers. 
The computation of this section was carried out on a $2.2\,\mathrm{GHz}$ 
Intel Core i5 processor with $8\,\mathrm{GB}$ RAM. 

As a concrete instance, consider the problem setting shown in 
\reffig{fig:gs5x2}. 
The external force is applied at the bottom right node, but this time 
its direction and magnitude are assumed to be uncertain. 
Without loss of generality, let $p_{1}$ and $p_{2}$ denote the 
horizontal and vertical components, respectively, of this external force. 
The set of possible realizations of the external load is defined by 
\begin{align*}
  P = \{ 
  (p_{1},p_{2},0,\dots,0)^{\top}
  \mid
  p_{1} = p^{0}_{1} \psi_{1} , 
  \
  p_{2} = p^{0}_{2} \psi_{2} , 
  \
  \| (\psi_{1},\psi_{2}) \| \le 1
  \} , 
\end{align*}
where $p^{0}_{1}=30\,\mathrm{kN}$ and $p^{0}_{2}=100\,\mathrm{kN}$. 
With referring to \eqref{eq:def.compliance}, we see that the compliance 
in the worst case is given by 
\begin{align*}
  \hat{\pi}(\bi{x}) = \sup \{ 
  2 \bi{p}^{\top} \bi{u} - \bi{u}^{\top} K(\bi{x}) \bi{u}
  \mid
  \bi{u} \in \Re^{d} , \
  \bi{p} \in P
  \} . 
\end{align*}
In the following, we consider the minimization problem of this function. 

When the constraint on the number of nodes is not considered, it is known that 
this optimization problem can be recast as semidefinite programming 
(SDP) \cite{BtN97}. 
Since the upper bound constraint on the number of nodes is treated as 
presented in Section~\ref{sec:node.integer}, the optimization problem 
under this constraint can be recast as 
mixed-integer semidefinite programming (MISDP). 
For comparison, we solve this MISDP problem with YALMIP~\cite{Lof04}, 
which finds a global optimal solution with a branch-and-bound 
method \cite{Lof04}. 
We used YALMIP with the default setting, where SDP subproblems are 
solved with SeDuMi ver.~1.3 \cite{Pol05,Stu99}. 
Alternatively, consider the problem obtained by replacing the objective 
function of \eqref{P.compliance.node.def} by $\hat{\pi}$. 
It is fairly straightforward to apply the ADMM in Section~\ref{sec:admm} 
to this optimization problem. 
The subproblem solved to update the variable $\bi{x}$ at each iteration 
is formulated as SDP. 

\begin{figure}[tp]
  \centering
  \begin{subfigure}[b]{0.32\textwidth}
    \centering
    \includegraphics[scale=0.25]{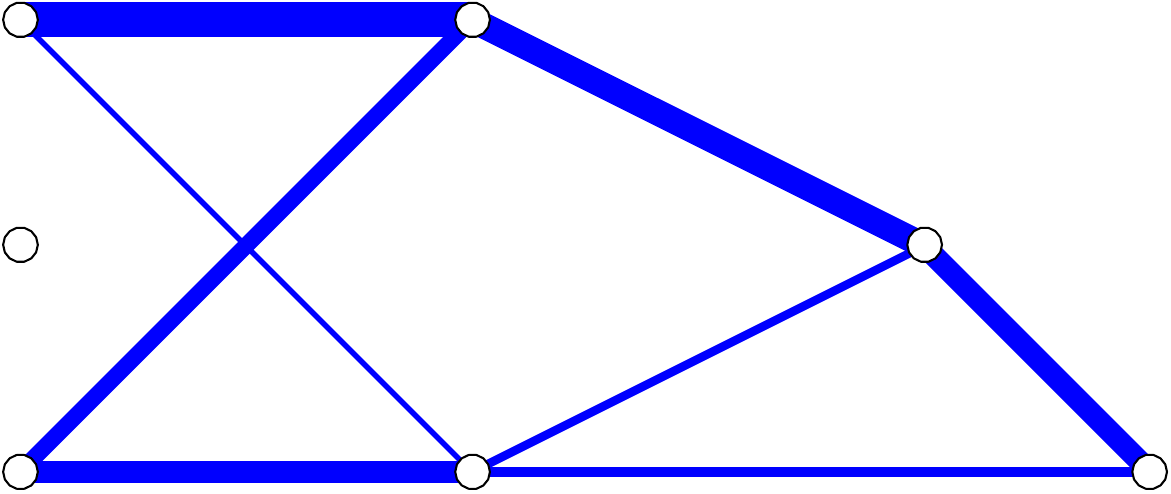}
    \caption{}
    \label{fig:x5_y2_robust}
  \end{subfigure}
  \hfill
  \begin{subfigure}[b]{0.32\textwidth}
    \centering
    \includegraphics[scale=0.25]{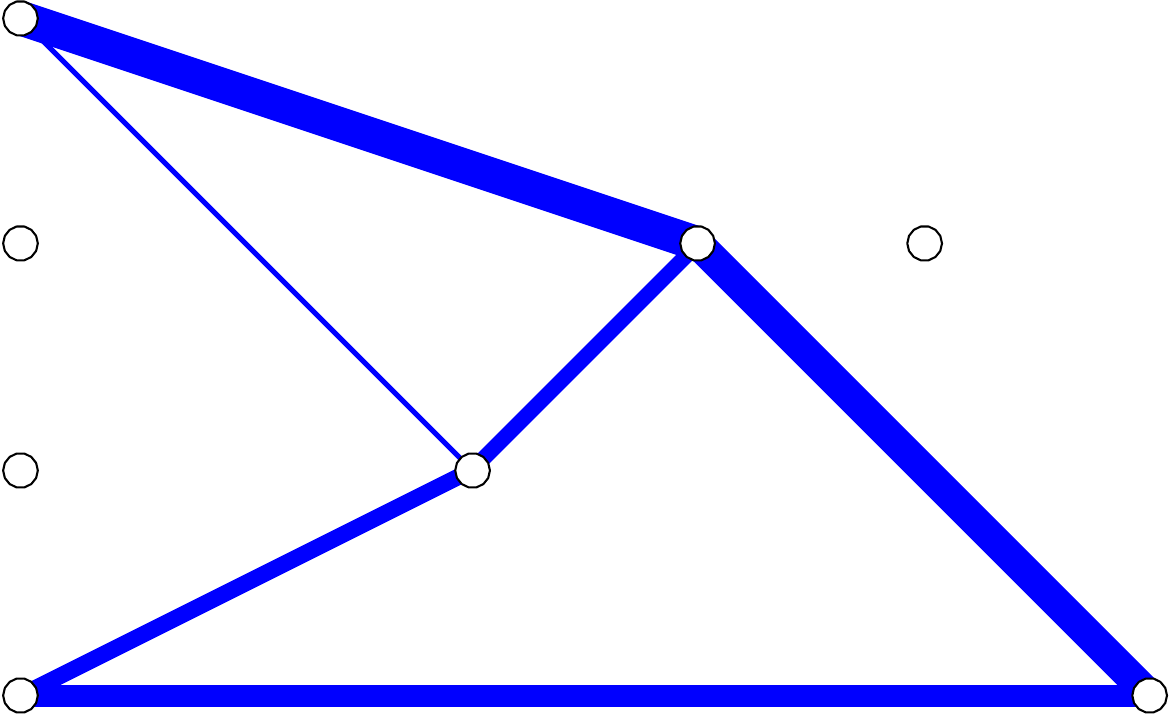}
    \caption{}
    \label{fig:x5_y3_robust}
  \end{subfigure}
  \hfill
  \begin{subfigure}[b]{0.32\textwidth}
    \centering
    \includegraphics[scale=0.25]{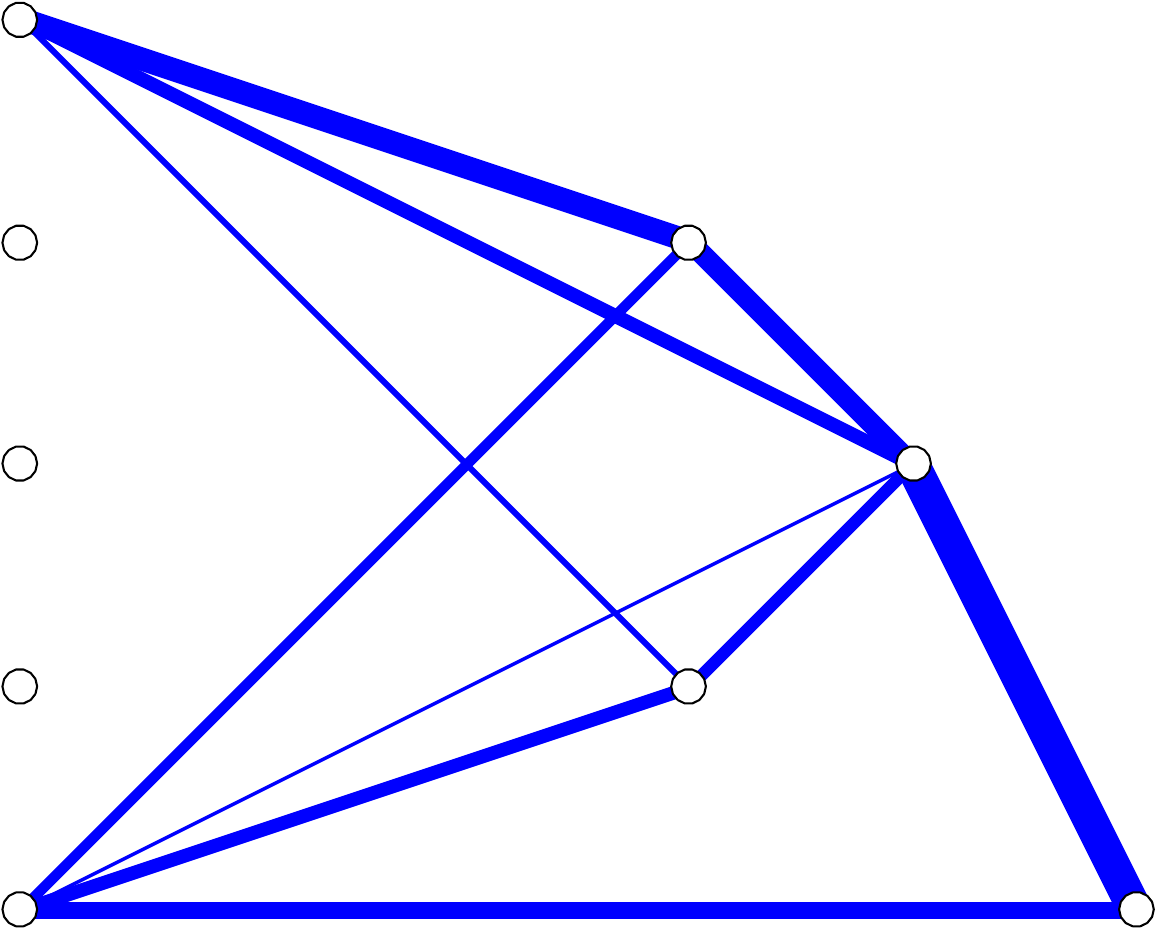}
    \caption{}
    \label{fig:x5_y4_robust}
  \end{subfigure}
  \par\medskip
  \begin{subfigure}[b]{0.45\textwidth}
    \centering
    \includegraphics[scale=0.40]{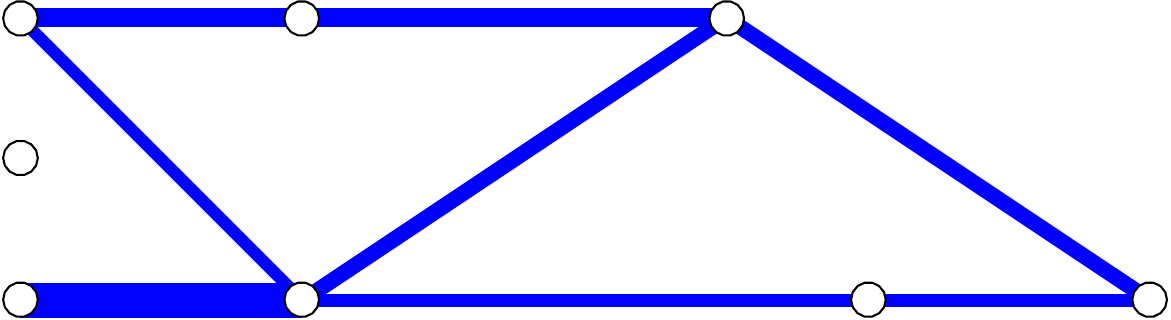}
    \caption{}
    \label{fig:x8_y2_robust}
  \end{subfigure}
  \hfill
  \begin{subfigure}[b]{0.45\textwidth}
    \centering
    \includegraphics[scale=0.40]{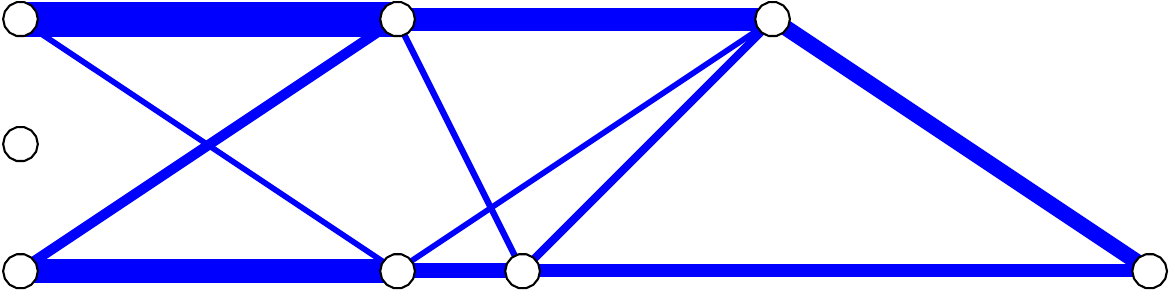}
    \caption{}
    \label{fig:x9_y2_robust}
  \end{subfigure}
  \par\medskip
  \begin{subfigure}[b]{0.45\textwidth}
    \centering
    \includegraphics[scale=0.40]{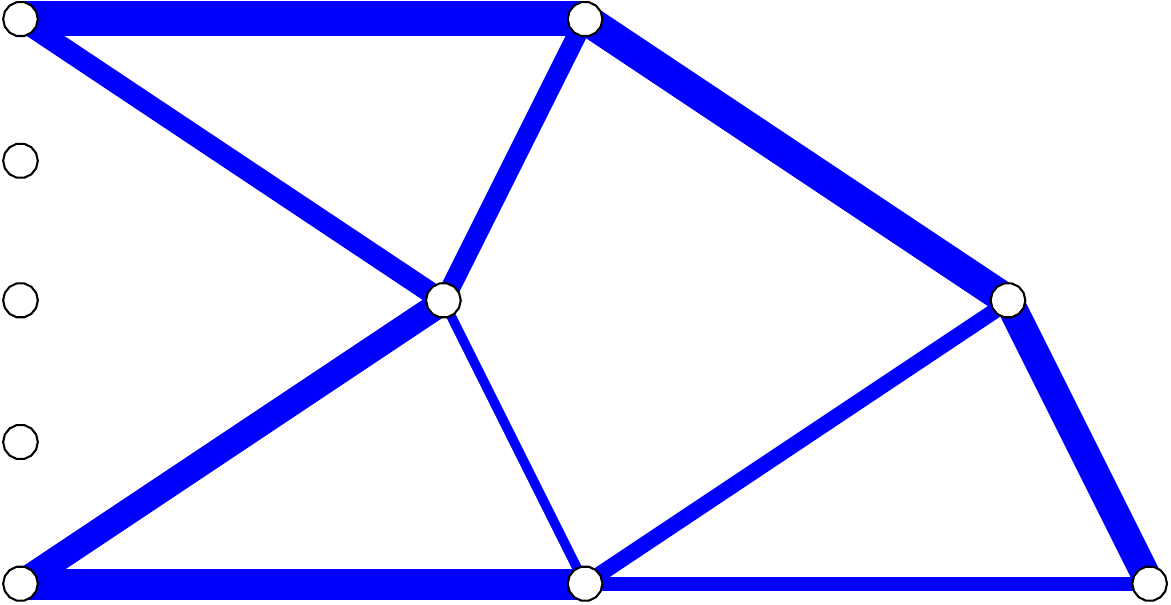}
    \caption{}
    \label{fig:x8_y4_robust}
  \end{subfigure}
  \hfill
  \begin{subfigure}[b]{0.45\textwidth}
    \centering
    \includegraphics[scale=0.40]{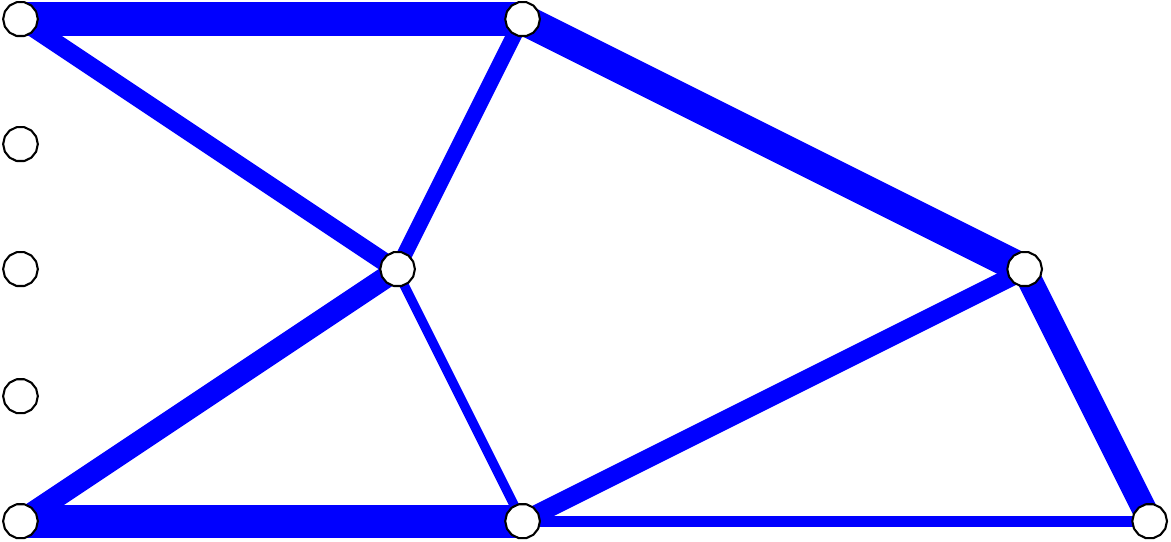}
    \caption{}
    \label{fig:x9_y4_robust}
  \end{subfigure}
  \par\medskip
  \begin{subfigure}[b]{0.45\textwidth}
    \centering
    \includegraphics[scale=0.40]{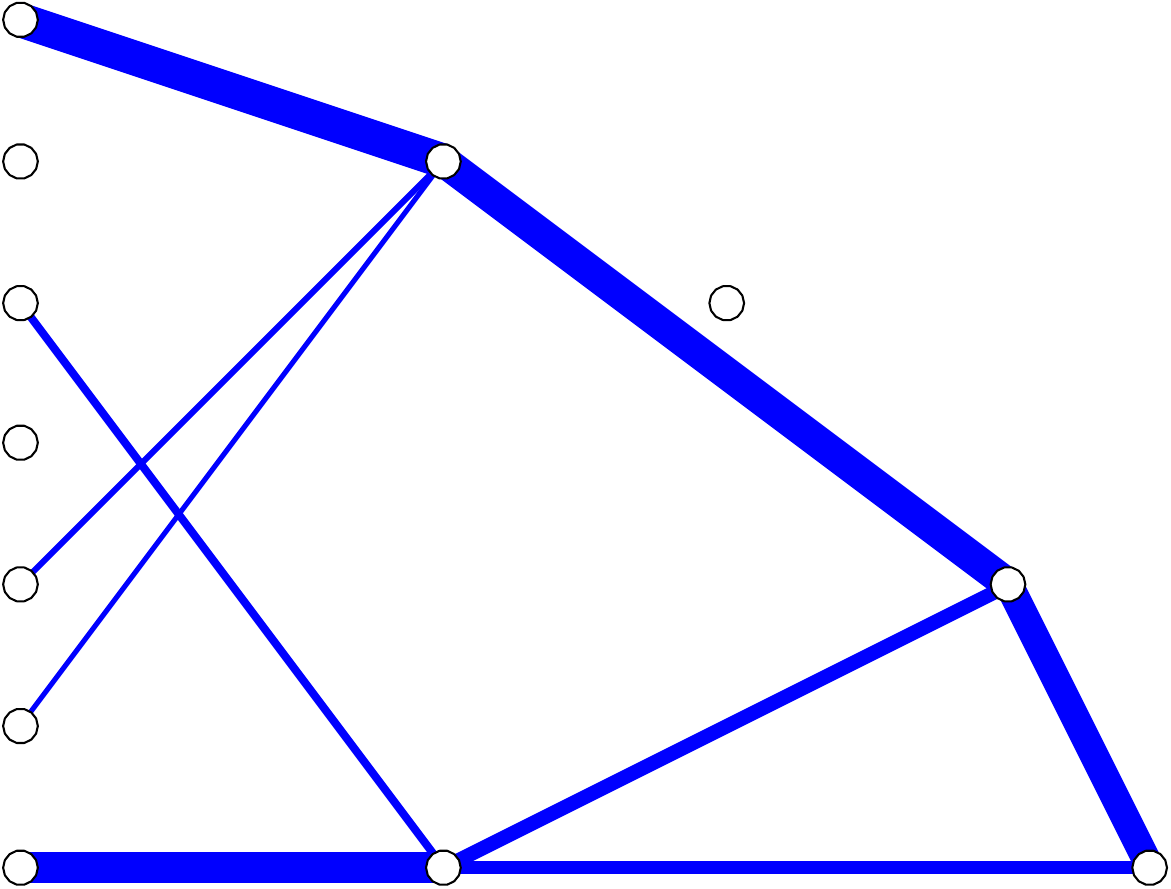}
    \caption{}
    \label{fig:x8_y6_robust}
  \end{subfigure}
  \hfill
  \begin{subfigure}[b]{0.45\textwidth}
    \centering
    \includegraphics[scale=0.40]{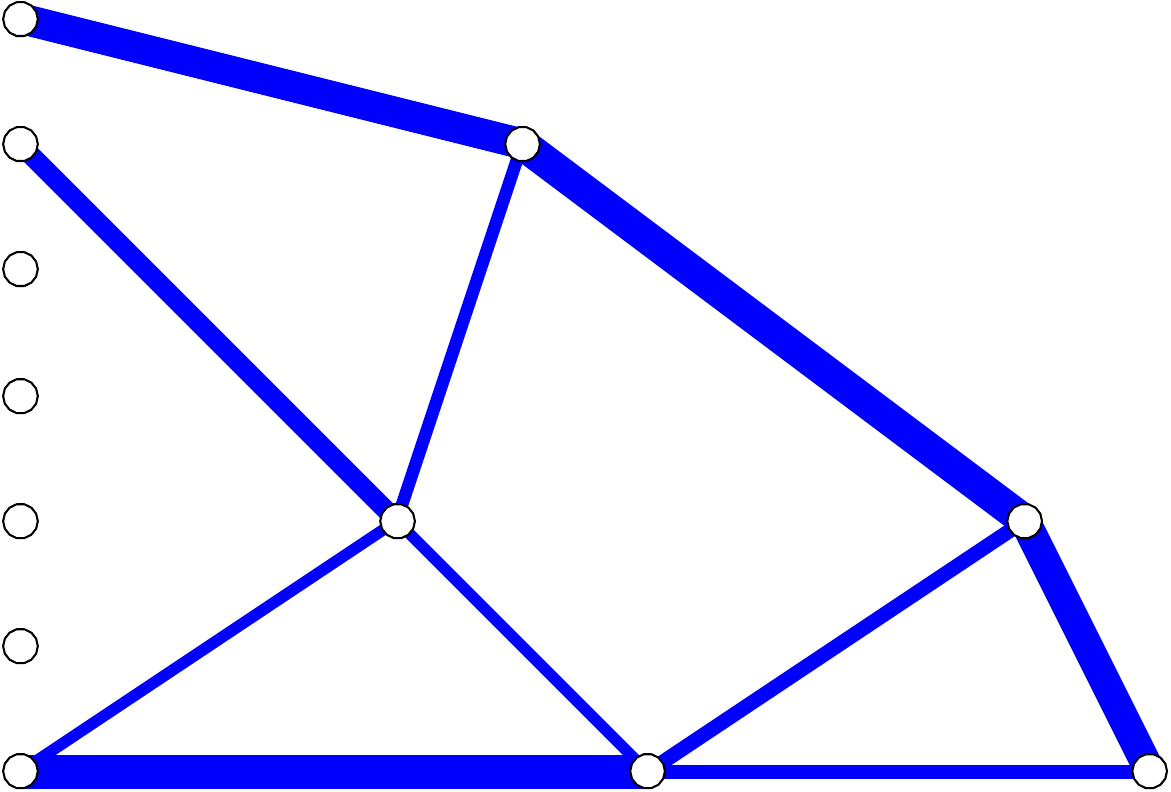}
    \caption{}
    \label{fig:x9_y6_robust}
  \end{subfigure}
  \caption{The solutions obtained by the ADMM applied to the robust 
  optimization. 
  \subref{fig:x5_y2_robust} $(N_{X},N_{Y},n)=(5,2,4)$; 
  \subref{fig:x5_y3_robust} $(N_{X},N_{Y},n)=(5,3,4)$; 
  \subref{fig:x5_y4_robust} $(N_{X},N_{Y},n)=(5,4,4)$; 
  \subref{fig:x8_y2_robust} $(N_{X},N_{Y},n)=(8,2,5)$; 
  \subref{fig:x9_y2_robust} $(N_{X},N_{Y},n)=(9,2,5)$; 
  \subref{fig:x8_y4_robust} $(N_{X},N_{Y},n)=(8,4,5)$; 
  \subref{fig:x9_y4_robust} $(N_{X},N_{Y},n)=(9,4,5)$; 
  \subref{fig:x8_y6_robust} $(N_{X},N_{Y},n)=(8,6,5)$; and 
  \subref{fig:x9_y6_robust} $(N_{X},N_{Y},n)=(9,6,5)$. 
  }
  \label{fig:robust.ADMM}
\end{figure}

\begin{figure}[tp]
  \centering
  \begin{subfigure}[b]{0.45\textwidth}
    \centering
    \includegraphics[scale=0.40]{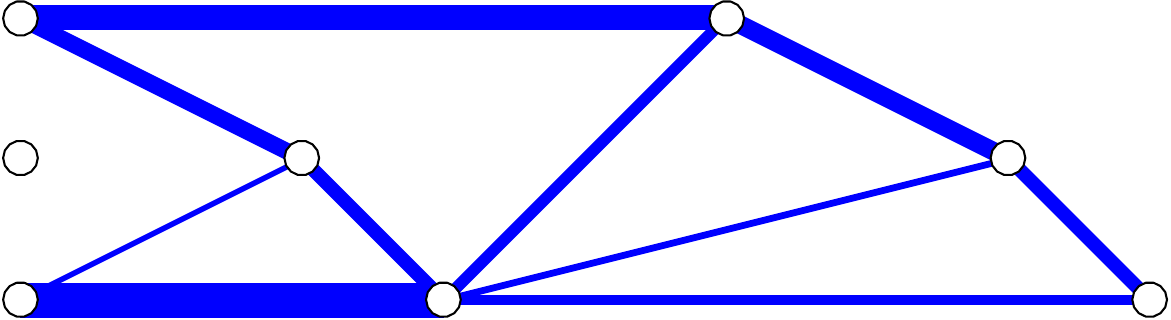}
    \caption{}
    \label{fig:x8_y2_bnb}
  \end{subfigure}
  \hfill
  \begin{subfigure}[b]{0.45\textwidth}
    \centering
    \includegraphics[scale=0.40]{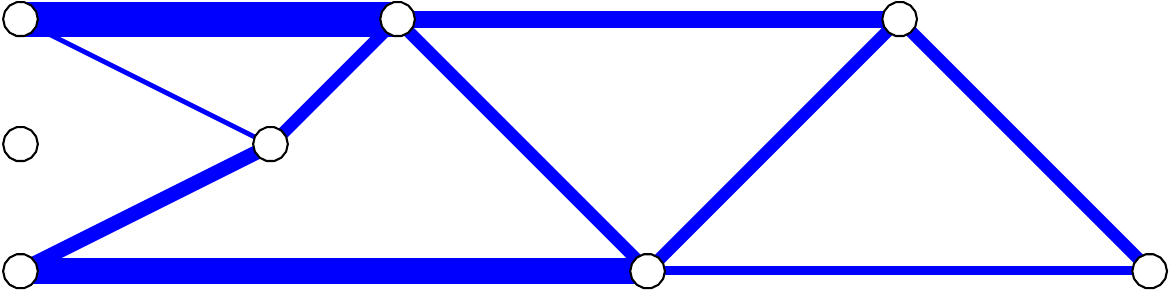}
    \caption{}
    \label{fig:x9_y2_bnb}
  \end{subfigure}
  \par\medskip
  \begin{subfigure}[b]{0.45\textwidth}
    \centering
    \includegraphics[scale=0.40]{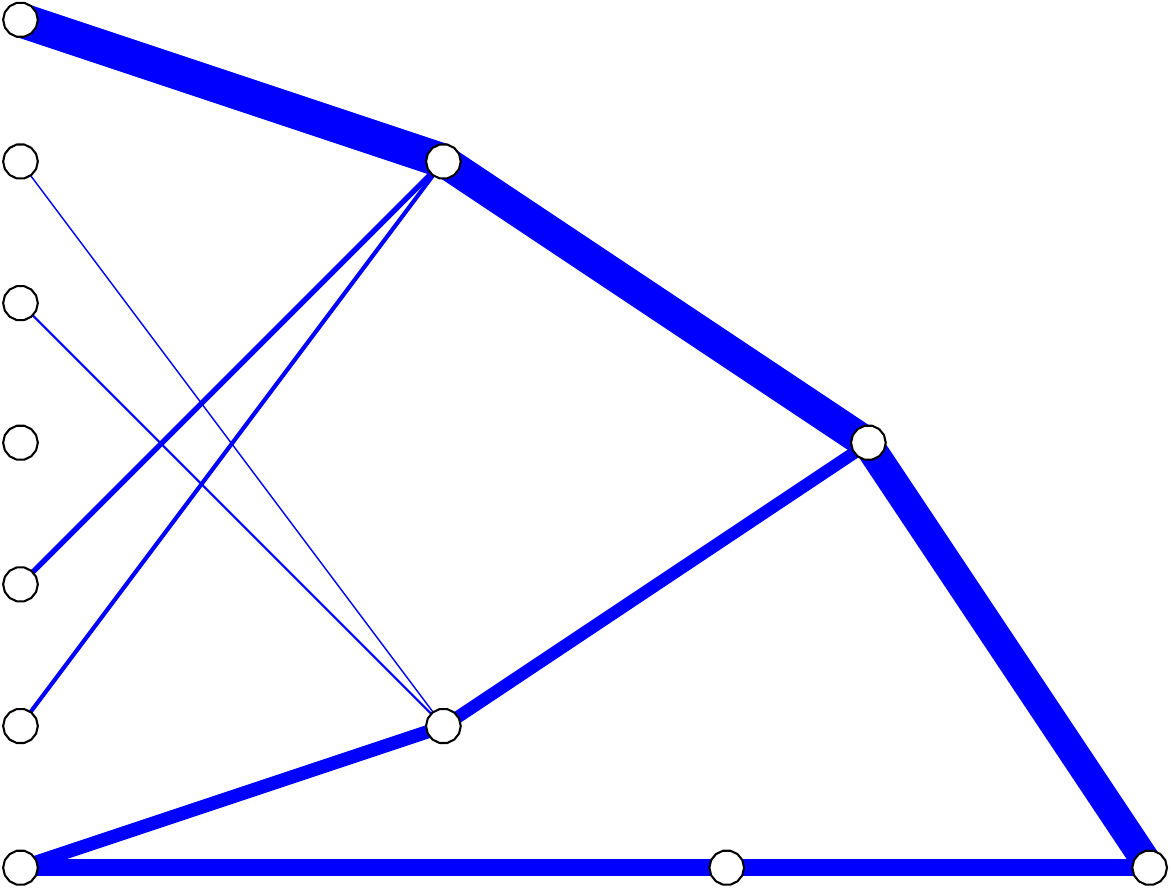}
    \caption{}
    \label{fig:x8_y6_bnb}
  \end{subfigure}
  \caption{The optimal solutions of the robust optimization obtained by YALMIP. 
  \subref{fig:x8_y2_bnb} $(N_{X},N_{Y},n)=(8,2,5)$; 
  \subref{fig:x9_y2_bnb} $(N_{X},N_{Y},n)=(9,2,5)$; and 
  \subref{fig:x8_y6_bnb} $(N_{X},N_{Y},n)=(8,6,5)$. 
  }
  \label{fig:robust.bnb}
\end{figure}

\reftab{tab:ex.I.SDP} reports the computational results. 
For five instances, $(N_{X},N_{Y},n)=(5,2,4)$, $(5,3,4)$, $(5,4,4)$, 
$(8,4,5)$, and $(9,4,5)$, the ADMM approach found the global optimal 
solutions. 
In all these cases, the computational cost of the ADMM is smaller than 
that of YALMIP. 
The difference of computational cost increases as the problem size increases. 
For $(N_{X},N_{Y},n)=(9,6,5)$, YALMIP did not terminate after $300$ 
iterations. 
The best solution is same as the one found by the ADMM with initial 
point (A), but a better solution was found by the ADMM with initial 
point (B). 
For three instances, $(N_{X},N_{Y},n)=(8,2,5)$, $(9,2,5)$, and $(8,6,5)$, 
the solutions found by the ADMM are not optimal. 
The difference between the obtained objective value and the optimal value 
is 7\% or less, like in the cases in section~\ref{sec:ex.n5}. 
The computational time required by YALMIP is more than $10$ times (in 
some cases, more than $100$ times) larger than that of the ADMM. 
\reffig{fig:robust.ADMM} collects the solutions obtained by the ADMM. 
The global optimal solutions that could not be obtained by the ADMM are 
shown in \reffig{fig:robust.bnb}. 
The set of nodes in \reffig{fig:x8_y2_robust} is much different from 
that in \reffig{fig:x8_y2_bnb}. 
The solution in \reffig{fig:x9_y2_robust} has only one node that is not 
included in the solution in \reffig{fig:x9_y2_bnb}. 
Similarly, the difference between the solutions in 
\reffig{fig:x8_y6_robust} and \reffig{fig:x8_y6_bnb} is the location of 
one node. 

\begin{figure}[tp]
  \centering
  \begin{subfigure}[b]{0.45\textwidth}
    \centering
    \includegraphics[scale=0.40]{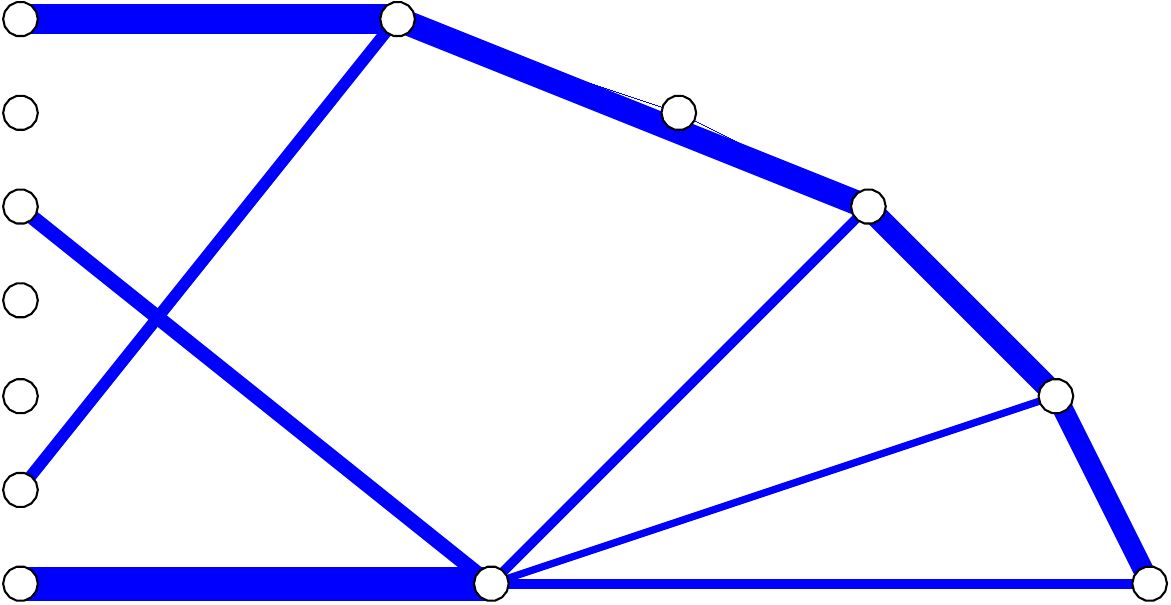}
    \caption{}
    \label{fig:x12_y6_robust}
  \end{subfigure}
  \hfill
  \begin{subfigure}[b]{0.45\textwidth}
    \centering
    \includegraphics[scale=0.40]{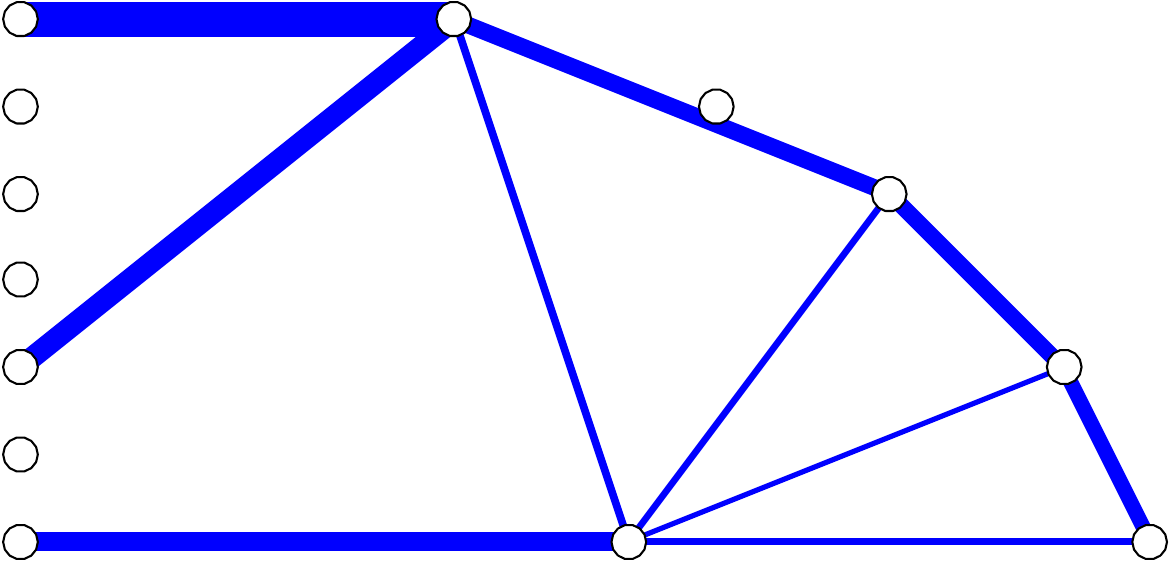}
    \caption{}
    \label{fig:x13_y6_robust}
  \end{subfigure}
  \par\medskip
  \begin{subfigure}[b]{0.45\textwidth}
    \centering
    \includegraphics[scale=0.40]{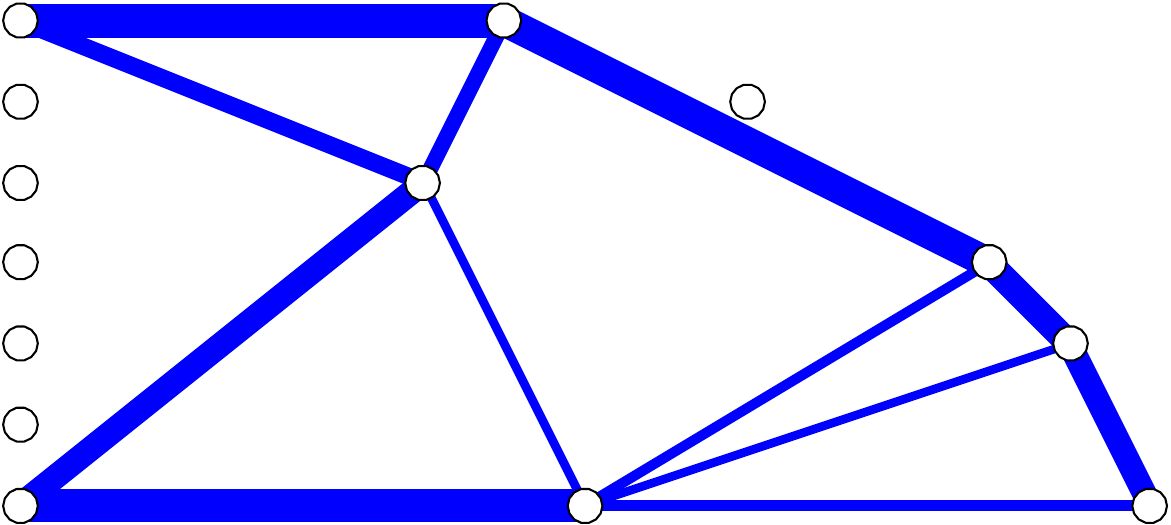}
    \caption{}
    \label{fig:x14_y6_robust}
  \end{subfigure}
  \caption{
  The solutions obtained by ADMM for the large-scale robust optimization 
  problems. 
  \subref{fig:x12_y6_robust} $(N_{X},N_{Y},n)=(12,6,6)$; 
  \subref{fig:x13_y6_robust} $(N_{X},N_{Y},n)=(13,6,6)$; 
  \subref{fig:x14_y6_robust} $(N_{X},N_{Y},n)=(14,6,7)$; 
  }
  \label{fig:x12_x_13_x14_robust}
\end{figure}

\reffig{fig:x12_x_13_x14_robust} shows the solutions obtained by the 
ADMM for problem instances with larger sizes. 
The computational results are listed in \reftab{tab:ex.I.SDP.large}. 
A global optimization method, YALMIP, cannot solve these problems within 
reasonable computational cost.

\begin{table}[bp]
  \centering
  \caption{Computational results of the ADMM approach applied to the 
  robust optimization with the cardinality constraint.}
  \label{tab:ex.I.SDP.large}
  \begin{tabular}{lrrrr}
    \toprule
    $(N_{X},N_{Y},n)$ &  Init.\ sol.\ & $w^{*}$ (J) & {\#}iter & Time (s) \\
    \midrule
    $(12,6,6)$ 
    & $*$ (A) & $7862.67$ & $10$ & $771.7$ \\
    &     (B) & $9572.47$ & $10$ & $730.5$ \\
    \midrule
    $(13,6,6)$ 
    & $*$ (A) & $10191.63$ & $10$ & $972.2$ \\
    &     (B) & $13680.44$ & $12$ & $1108.4$ \\
    \midrule
    $(14,6,7)$ 
    & $*$ (A) & $9772.59$ & $9$ & $1033.9$ \\
    &     (B) & $14980.99$ & $14$ & $1614.1$ \\
    \bottomrule
  \end{tabular}
\end{table}

\section{Conclusions}
\label{sec:conclude}

In this paper we have studied the compliance minimization of a truss 
with the limited number of nodes. 
It has been shown that this optimization problem can be formulated as 
the cardinality-constrained SOCP. 
We have proposed a simple and efficient heuristic based on ADMM. 

The problem considered in this paper can also be formulated as MISOCP 
involving the so-called big-M. 
In the numerical experiments, we have compared the proposed ADMM 
approach with a 
global optimization approach using the MISOCP formulation. 
For small-size problem instances, it has been confirmed that ADMM finds 
a global optimal solution. 
For middle-size instances, the objective value of the solution obtained 
by ADMM is often close to the optimal value. 
The number of iterations of ADMM is almost same for instances with 
different sizes. 
In contrast, the computational cost required by a standard MISOCP solver 
highly depends on instances, even if the instances have similar sizes. 

In the numerical experiments, it has also been illustrated that, for 
some problem instances, the compliance minimization problem of a truss 
has some different optimal solutions, and the number of nodes can be 
decreased without losing the optimality. 
In most of the other cases, the number of nodes can be decreased at the 
expense of only small increase of the compliance.

\paragraph{Acknowledgments}

The work of the first author is partially supported by 
JSPS KAKENHI 15KT0109 and 17K06633.

\end{document}